\documentclass{IEEEtran}
\usepackage{latexsym}
\usepackage{cite}
\usepackage{amssymb}
\usepackage{bbm}
\usepackage{textpos}
\usepackage{enumerate}
\usepackage{bm}
\usepackage{mathtools}
\usepackage[usenames,dvipsnames]{xcolor}
\usepackage{stmaryrd}
\usepackage{amsmath, amssymb,accents}
\usepackage{flushend}
\usepackage{dsfont}
\usepackage{algorithmicx}
\usepackage{algorithm} 
\usepackage[noend]{algpseudocode} %
\usepackage{comment}
\usepackage{accents}
\DeclareMathOperator{\supp}{supp}
\usepackage[pdftex]{graphicx}
\usepackage{psfrag}
\usepackage{units}
\usepackage{subfig}
\usepackage{setspace}
\usepackage{theorem}
\usepackage{mathtools}
\usepackage{fancyhdr}
\allowdisplaybreaks

\newtheorem{definition}{Definition}

\newtheorem{theorem}{Theorem}
\newtheorem{lemma}{Lemma}
\newtheorem{remark}{Remark}

\newtheorem{corollary}{Corollary}
\newtheorem{proposition}{Proposition}

\newcommand\numberthis{\addtocounter{equation}{1}\tag{\theequation}}
{} 

\begin{document}

\makeatletter

\DeclareRobustCommand{\rchi}{{\mathpalette\irchi\relax}}
\newcommand{\irchi}[2]{\raisebox{\depth}{$#1\chi$}} 
\newcommand{\cA}{\mathcal{A}}
\newcommand{\cB}{\mathcal{B}}
\newcommand{\cC}{\mathcal{C}}
\newcommand{\cD}{\mathcal{D}}
\newcommand{\cE}{\mathcal{E}}
\newcommand{\cF}{\mathcal{F}}
\newcommand{\cG}{\mathcal{G}}
\newcommand{\cH}{\mathcal{H}}
\newcommand{\cI}{\mathcal{I}}
\newcommand{\cJ}{\mathcal{J}}
\newcommand{\cK}{\mathcal{K}}
\newcommand{\cL}{\mathcal{L}}
\newcommand{\cM}{\mathcal{M}}
\newcommand{\cN}{\mathcal{N}}
\newcommand{\cO}{\mathcal{O}}
\newcommand{\cP}{\mathcal{P}}
\newcommand{\cQ}{\mathcal{Q}}
\newcommand{\cR}{\mathcal{R}}
\newcommand{\cS}{\mathcal{S}}
\newcommand{\cT}{\mathcal{T}}
\newcommand{\cU}{\mathcal{U}}
\newcommand{\cV}{\mathcal{V}}
\newcommand{\cW}{\mathcal{W}}
\newcommand{\cX}{\mathcal{X}}
\newcommand{\cY}{\mathcal{Y}}
\newcommand{\cZ}{\mathcal{Z}}
\newcommand{\FSG}{\mathcal{F}^{(\mathsf{SG})}_{d,K}}
\newcommand{\FSGO}{\mathcal{F}^{(\mathsf{SG})}_{1,K}}
\newcommand{\BB}{\mathbb{B}}
\newcommand{\CC}{\mathbb{C}}
\newcommand{\DD}{\mathbb{D}}
\newcommand{\EE}{\mathbb{E}}
\newcommand{\FF}{\mathbb{F}}
\newcommand{\GG}{\mathbb{G}}
\newcommand{\HH}{\mathbb{H}}
\newcommand{\II}{\mathbb{I}}
\newcommand{\JJ}{\mathbb{J}}
\newcommand{\KK}{\mathbb{K}}
\newcommand{\LL}{\mathbb{L}}
\newcommand{\MM}{\mathbb{M}}
\newcommand{\NN}{\mathbb{N}}
\newcommand{\OO}{\mathbb{O}}
\newcommand{\PP}{\mathbb{P}}
\newcommand{\QQ}{\mathbb{Q}}
\newcommand{\RR}{\mathbb{R}}
\newcommand{\TT}{\mathbb{T}}
\newcommand{\UU}{\mathbb{U}}
\newcommand{\VV}{\mathbb{V}}
\newcommand{\WW}{\mathbb{W}}
\newcommand{\XX}{\mathbb{X}}
\newcommand{\YY}{\mathbb{Y}}
\newcommand{\ZZ}{\mathbb{Z}}
\newcommand*{\dd}{\, \mathsf{d}}
\newcommand{\vasti}{\bBigg@{3.5 }}
\newcommand{\vast}{\bBigg@{4}}
\newcommand{\Vast}{\bBigg@{5}}
\newcommand{\Vastt}{\bBigg@{7}}
\makeatother
\newcommand{\be}{\begin{equation}}
\newcommand{\ee}{\end{equation}}
\newcommand{\ba}{\begin{align}}
\newcommand{\ea}{\end{align}}
\newcommand{\baa}{\begin{align*}}
\newcommand{\eaa}{\end{align*}}
\newcommand{\ber}{$\ \mbox{Ber}$}
\newcommand{\argmin}{\mathop{\mathrm{argmin}}}
\newcommand{\argmax}{\mathop{\mathrm{argmax}}}
\newcommand{\ubar}[1]{\underaccent{\bar}{#1}}
\newcommand*{\gauss}{\varphi_\sigma}
\newcommand*{\gausss}{\varphi_{\frac{\sigma}{\sqrt{2}}}}
\newcommand*{\gaussI}{\varphi_1}
\newcommand*{\gausssI}{\varphi_{\frac{1}{\sqrt{2}}}}
\newcommand*{\Gauss}{\mathcal{N}_\sigma}
\newcommand*{\Gausss}{\mathcal{N}_{\frac{\sigma}{\sqrt{2}}}}
\newcommand{\p}[1]{\left(#1\right)}
\newcommand{\s}[1]{\left[#1\right]}
\newcommand{\cp}[1]{\left\{#1\right\}}
\newcommand{\abs}[1]{\left|#1\right|}
\newcommand{\nmc}{n_{\mathsf{MC}}}
\newcommand{\var}{\mathsf{var}}
\newcommand*{\E}{\mathbb E}
\newcommand*{\Chi}{\rchi^2}
\newcommand{\fillater}[1] {{\Large \color{red}[#1]}}
\makeatletter
\def\Ddots{\mathinner{\mkern1mu\raise\p@
\vbox{\kern7\p@\hbox{.}}\mkern2mu
\raise4\p@\hbox{.}\mkern2mu\raise7\p@\hbox{.}\mkern1mu}}
\makeatother

\title{Convergence of Smoothed Empirical Measures with Applications to Entropy Estimation}

\author{Ziv Goldfeld, Kristjan Greenewald, Jonathan Niles-Weed and Yury Polyanskiy

		\thanks{This work was partially supported by the MIT-IBM Watson AI Lab. The work of Z. Goldfeld and Y. Polyanskiy was also supported in part by the National Science Foundation CAREER award under grant agreement CCF-12-53205, by the Center for Science of Information (CSoI), an NSF Science and Technology Center under grant agreement CCF-09-39370, and a grant from Skoltech--MIT Joint Next Generation Program (NGP). The work of Z. Goldfeld was also supported by the Rothschild postdoc fellowship. The work of J. Niles-Weed was partially supported by the Josephine de K\'arm\'an Fellowship.
		\newline This work was presented in part at the 2018 IEEE International Conference on the Science of Electrical Engineering (ICSEE-2018), Eilat, Israel, and in part at the 2019 IEEE International Symposium on Information Theory (ISIT-2019), Paris, France.
		\newline Z. Goldfeld is with the School of Electrical and Computer Engineering, Cornell University, Ithaca, NY 14850 US (e-mail: goldfeld@cornell.edu). K. Greenewald is with the MIT-IBM Watson AI Lab, Cambridge, MA 02142 US (email: kristjan.h.greenewald@ibm.com). J. Niles-Weed is with the Courant Institute of Mathematical Sciences and Center for Data Science, New York University, New York, NY 10003 US (email: jnw@cims.nyu.edu). Y. Polyanskiy is with the Department of Electrical Engineering and Computer Science, Massachusetts Institute of Technology, Cambridge, MA 02139 US (e-mail: yp@mit.edu).}
}
\maketitle

\begin{abstract}

This paper studies convergence of empirical measures smoothed by a Gaussian kernel. Specifically, consider approximating $P\ast\Gauss$, for $\Gauss\triangleq\mathcal{N}(0,\sigma^2 \mathrm{I}_d)$, by $\hat{P}_n\ast\Gauss$ under different statistical distances, where $\hat{P}_n$ is the empirical measure. We examine the convergence in terms of the Wasserstein distance, total variation (TV), Kullback-Leibler (KL) divergence, and $\Chi$-divergence. We show that the approximation error under the TV distance and 1-Wasserstein distance ($\mathsf{W}_1$) converges at the rate $e^{O(d)} n^{-1/2}$ in remarkable contrast to a (typical) $n^{-\frac{1}{d}}$ rate for unsmoothed $\mathsf{W}_1$ (and $d\ge 3$). Similarly, for the KL divergence, squared 2-Wasserstein distance ($\mathsf{W}_2^2$), and $\Chi$-divergence, the convergence rate is $e^{O(d)} n^{-1}$, but only if $P$ achieves finite input-output $\Chi$ mutual information across the additive white Gaussian noise (AWGN) channel. If the latter condition is not met, the rate changes to $\omega\left(n^{-1}\right)$ for the KL divergence and $\mathsf{W}_2^2$, while the $\Chi$-divergence becomes infinite -- a curious dichotomy. 

As an application we consider estimating the differential entropy $h(S+Z)$, where $S\sim P$ and $Z\sim\Gauss$ are independent $d$-dimensional random variables. The distribution $P$ is unknown and belongs to some nonparametric class, but $n$ independently and identically distributed (i.i.d) samples from it are available. Despite the regularizing effect of noise, we first show that any good estimator (within an additive gap) for this problem must have a sample complexity that is exponential in $d$. We then leverage the above empirical approximation results to show that the absolute-error risk of the plug-in estimator converges as $e^{O(d)} n^{-1/2}$, thus attaining the parametric rate in $n$. This establishes the~plug-in estimator as minimax rate-optimal for the considered problem, with sharp dependence of the convergence rate both in $n$ and $d$. We provide numerical results comparing the performance of the plug-in estimator to that of general-purpose (unstructured) differential entropy estimators (based on kernel density estimation (KDE) or $k$ nearest neighbors (kNN) techniques) applied to samples of $S+Z$. These results reveal a significant empirical superiority of the plug-in to state-of-the-art KDE and kNN methods. As a motivating utilization of the plug-in approach, we estimate information~flows in deep neural networks and discuss Tishby's Information Bottleneck and the compression conjecture, among others.



\end{abstract}



\section{Introduction}

This work is motivated by a new nonparametric and high-dimensional functional estimation problem, which we refer to as `differential entropy estimation under Gaussian convolutions.' The goal of this problem is to estimate the differential entropy $h(S+Z)$ based on samples of $S$, while knowing the distribution of $Z\sim\Gauss\triangleq\mathcal{N}(0,\sigma^2\mathrm{I}_d)$ which is independent of $S$. The analysis of the estimation risk reduces to evaluating the expected 1-Wasserstein distance or $\Chi$-divergence between $P\ast\Gauss$ and $\hat{P}_{S^n}\ast\Gauss$, where $S^n\triangleq(S_1,\ldots,S_n)$ are i.i.d. samples from $P$ and $\hat{P}_{S^n}=\frac{1}{n}\sum_{i=1}^n\delta_{S_i}$ is the empirical measure.\footnote{Here, $\delta_{S_i}$ stands for the Dirac measure at $S_i$.} Due to the popularity of the additive Gaussian noise model, we start by exploring this smoothed empirical approximation problem in detail, under several additional statistical distances. 

\subsection{Convergence of Smooth Empirical Measures}

Consider the empirical approximation error $\mathbb{E}\delta(\hat{P}_{S^n}\ast\Gauss,P\ast\Gauss)$ under some statistical distance $\delta$. We consider various choices of $\delta$, including the 1-Wasserstein and (squared) 2-Wasserstein distances, total variation (TV), Kullback-Leibler (KL) divergence, and $\Chi$-divergence. We show that, when $P$ is subgaussian, the approximation error under the 1-Wasserstein and TV distances drops at the \emph{parametric} rate for \emph{any} dimension $d$. The exact rate is $c^dn^{-1/2}$, for a constant $c$. The parametric convergence rate is also attained by the squared 2-Wasserstein distance, KL divergence, and $\Chi$-divergence, so long as $P$ achieves finite input-output $\Chi$ mutual information across the additive white Gaussian noise (AWGN) channel. We show that this condition is always met for subgaussian $P$ in the low signal-to-noise ratio (SNR) regime. This fast convergence is remarkable since classical (unconvolved) empirical approximation suffers from the so-called curse of dimensionality. For instance, the empirical 1-Wasserstein distance $\mathbb{E}\mathsf{W}_1(\hat{P}_{S^n},P)$ is known to decay at most as $n^{-\frac{1}{d}}$ \cite{dudley1969speed}, which is sharp for all $d>2$ (see \cite{dobric1995asymptotics} and \cite{FournierGuillin2015} for sharper results). Convolving $P$ and $\hat{P}_{S^n}$ with $\Gauss$ improves the rate from $n^{-\frac{1}{d}}$ to $c^d n^{-1/2}$ (or $c^d n^{-1}$ for squared distances). The latter has a milder dependence on $d$ and can be dominated with practical choices of $n$, even in relatively high~dimensions.

The $\Chi$-divergence $\mathbb{E}\Chi\left(\hat{P}_{S^n}\ast\Gauss\middle\|P\ast\Gauss\right)$ presents a particularly curious behavior: it converges as\footnote{Recall $\Chi$ is a squared distance.} $\frac{1}{n}$ for low SNR and possibly diverges when SNR is high. 
To demonstrate a diverging scenario we construct a subgaussian $P$ for which $\mathbb{E}\Chi\left(\hat{P}_{S^n}\ast\Gauss\middle\|P\ast\Gauss\right) = \infty$ whenever the subgaussian constant is greater than or equal to $\sqrt{2}\sigma$. We also show that, when the $\Chi$-divergence is infinite, the convergence rate of the squared 2-Wasserstein distance and KL divergence are strictly slower than the parametric rate. All of these empirical approximation results are gathered in Section \ref{SEC:distances}. 

\subsection{Differential Entropy Estimation Under Smoothing}

We then apply these empirical approximation results to study the estimation of $h(S+Z)$, where $S\sim P$ is an arbitrary (continuous, discrete, or mixed) $\mathbb{R}^d$-valued random variable and $Z\sim \mathcal{N}_\sigma$ is an isotropic Gaussian. The differential entropy is estimated using $n$ i.i.d. samples $S^n$ from $P$ and assuming $\sigma$ is known. To investigate the decision-theoretic fundamental limit, we consider the minimax absolute-error risk
\begin{equation}
    \mathcal{R}^\star(n,\sigma,\mathcal{F}_d)\triangleq\inf\limits_{\hat{h}}\sup\limits_{P\in\mathcal{F}_d}\mathbb{E}\left|h(P\ast\Gauss)-\hat{h}(S^n,\sigma)\right|,\label{EQ:risk}
\end{equation}
where $\mathcal{F}_d$ is a nonparametric class of distributions and $\hat{h}$ is the estimator. The sample complexity $n^\star(\eta,\sigma,\mathcal{F}_d)$ is the smallest number of samples needed for estimating $h(P\ast\Gauss)$ within an additive gap $\eta$. We aim to understand whether having access~to `clean' samples of $S$ can improve estimation performance (theoretically and empirically) compared to when only `noisy' samples of $S+Z$ are available and the distribution of $Z$ is unknown.

Our results establish the plug-in estimator as minimax rate-optimal for the considered problem. Defining $\mathsf{T}_\sigma(P)\triangleq h(P\ast\Gauss)$ as the functional of interest, the plug-in estimator is $\mathsf{T}_\sigma\big(\hat{P}_{S^n}\big)=h\big(\hat{P}_{S^n}\ast\Gauss\big)$. Plug-in techniques are suboptimal for vanilla discrete (Shannon) and differential entropy estimation (see \cite{Paninski03EstMI_neuro_biology} and \cite{han2017optimal}, respectively). Nonetheless, we show that $h\big(\hat{P}_{S^n}\ast\Gauss\big)$ attains the parametric estimation rate of $O_{\sigma,d}(n^{-1/2})$ when $P$ is subgaussian, establishing the optimality of the plug-in. 

We use the $\Chi$ empirical approximation result to prove the parametric risk convergence rate when $P$ has bounded support. The result is then extended to (unbounded) subgaussian $P$ via a separate argument. Specifically, we first bound the risk by a \emph{weighted} TV distance between $P\ast\Gauss$ and $\hat{P}_{S^n}\ast\Gauss$. This bound is derived by linking the two measures via the maximal TV-coupling. The subgaussianity of $P$ and the smoothing introduced by the Gaussian convolution are used to bound the weighted TV distance by a $e^{O(d)} n^{-1/2}$ term, with all constants explicitly characterized. Notably, while the convergence with $n$ is parametric, the derived rates still depends exponentially on $d$ though the prefactor.



A natural next question is whether the exponential dependence on $d$ is necessary. Answering in the affirmative, we prove that any good estimator of $h(P\ast\Gauss)$, within an additive gap $\eta$, has a sample complexity $n^\star(\eta,\sigma,\mathcal{F}_d)=\Omega\left(\frac{2^{\gamma(\sigma)d}}{\eta d}\right)$, where $\gamma(\sigma)$ is positive and monotonically decreasing in $\sigma$. The proof relates the estimation of $h(P\ast\Gauss)$ to estimating the discrete entropy of a distribution supported on a capacity-achieving codebook for an AWGN channel. Existing literature (e.g., \cite{valiant2010clt,discrete_entropy_est_Wu2016}) implies that the discrete problem has sample complexity exponential in $d$ (because this the growth rate of the codebook's size), which is then carried over to the original problem to establish the result. 

Finally, we focus on the practical estimation of $h(P\ast\Gauss)$. While the above results give necessary and sufficient conditions on the number of samples needed to drive the estimation error below a desired threshold, these are worst-case bounds. In practice, the unknown distribution $P$ may not follow the minimax rates, and the resulting estimation error could be smaller. As a guideline for setting $n$ in practice, we derive a lower bound on the bias of the plug-in estimator that scales as $\log\left(2^d n^{-1}\right)$. Our last step is to propose an efficient implementation of the plug-in estimator based on Monte Carlo (MC) integration. As the estimator amounts to computing the differential entropy of a \emph{known} Gaussian mixture, MC integration allows a simple and efficient computation. We bound the mean squared error (MSE) of the computed value by $c^{(\mathsf{MC})}_{\sigma,d}(n\cdot \nmc)^{-1}$, where $n$ is the number of centers in the mixture\footnote{The number of centers is the number of samples used for estimation.}, $\nmc$ is the number of MC samples, and $c^{(\mathsf{MC})}_{\sigma,d}=\Theta(d)$ is explicitly characterized. The proof uses the Gaussian Poincar{\'e} inequality to reduce the analysis to that of the log-mixture distribution gradient. Several simulations (including an estimation experiment over a small deep neural network (DNN) for a 3-dimensional spiral classification task) illustrate the superiority of the ad hoc plug-in approach over existing general-purpose estimators, both in the rate of error decay and scalability with dimension. 


\subsection{Related Differential Entropy Estimation Results}

General-purpose differential entropy estimators can be used in the considered setup by estimating $h(S+Z)$ using `noisy' samples of $S+Z$ (generated from the available samples of $S$). There are two prevailing approaches for estimating the nonsmooth differential entropy functional: (i) based on kernel density estimators (KDEs) \cite{kandasamy2015nonparametric,moonEntropy,moon2016improving}; and (ii) using $k$ nearest neighbor (kNN) techniques \cite{kozachenko1987sample,kraskov2004estimating,tsybakov1996root,sricharan2012estimation,kumar2013ensemble,singh2016finite,delattre2017kozachenko,jiao2017nearest,berrett2019efficient} (see also \cite{beirlant1997nonparametric,biau2015lectures} for surveys). Many performance analyses of such estimators restrict attention to nonparametric classes of smooth and compactly supported densities that are bounded away from zero (although the support may be unknown \cite{moonEntropy,moon2016improving}). Since the density associated with $P\ast\Gauss$ violates these assumptions, such results do not apply in our setup. The work of Tsybakov and van der Meulen \cite{tsybakov1996root} accounted for densities with unbounded support and exponentially decaying tails for $d=1$, but we are interested in the high-dimensional scenario.

Two recent works weakened or dropped the boundedness from below assumption in the high-dimensional setting, providing general-purpose estimators whose risk bounds are valid in our setup. In \cite{han2017optimal}, a KDE-based differential entropy estimator that also combines best polynomial approximation techniques was proposed. Assuming subgaussian densities with unbounded support, Theorem 2 of \cite{han2017optimal} bounded the estimation risk by\footnote{Multiplicative polylogarithmic factors are overlooked in this restatement.} $O\left(n^{-\frac{s}{s+d}}\right)$, where $s$ is a Lipschitz smoothness parameter assumed to satisfy $0<s\leq 2$. While the result is applicable for our setup when $P$ is compactly supported or subgaussian, the convergence rate for large $d$ is roughly $n^{-\frac{1}{d}}$. This rate deteriorates quickly with dimension and is unable to exploit the smoothness of $P\ast\Gauss$ due to the $s\leq 2$ restriction.\footnote{Such convergence rates are typical in estimating $h(p)$ under boundedness or smoothness conditions on $p$. Indeed, the results cited above (applicable in our framework or otherwise) as well as many others bound the estimation~risk as $O\mspace{-3mu}\left(\mspace{-2mu}n^{\mspace{-3mu}-\mspace{-1mu}\frac{\alpha}{\beta\mspace{-1mu}+\mspace{-1mu}d}}\mspace{-3mu}\right)$, where $\alpha,\beta$ are constants that may~depend on $s$ and~$d$.} This is to be expected because the results of \cite{han2017optimal} account for a wide class of density functions, including highly non-smooth~ones. 

In \cite{berrett2019efficient}, a weighted-KL (wKL) estimator (in the spirit of \cite{kumar2013ensemble}) was studied for smooth densities. In a major breakthrough, the authors proved that with a careful choice of weights the estimator is asymptotically efficient,\footnote{in the sense of, e.g., \cite[p. 367]{van2000asymptotic}.} under certain assumptions on the densities' speed of decay to zero (which captures $P\ast\Gauss$ when, e.g., $P$ is compactly supported). Despite its $O(n^{-1/2})$ risk convergence rate, however, the empirical performance of the estimator seems lacking (at least in our experiments, which use the code provided in \cite{berrett2019efficient}). As shown in Section \ref{SEC:simulations}, the plug-in estimator achieves superior performance even in rather simple scenarios of moderate dimension. 
The empirical performance of the estimator from \cite{berrett2019efficient} may originate from the dependence of its estimation risk on the dimension $d$, which was not characterized therein.


\subsection{Relation to Information Flows in Deep Neural Networks}

The considered differential entropy estimation problem is closely related to that of estimating information flows in DNNs. There has been a recent surge of interest in measuring the mutual information between selected groups of neurons in a DNN \cite{DNNs_ICLR2018,jacobsen2018revnet,liu2018understanding,gabrie2018entropy,reeves2017additivity,ICML_Info_flow2019}, partially driven by the Information Bottleneck (IB) theory \cite{tishby_DNN2015,DNNs_Tishby2017}. Much of the focus centers on the mutual information $I(X;T)$ between the input feature $X$ and a hidden activity vector $T$. However, as explained in \cite{ICML_Info_flow2019}, this quantity is vacuous in deterministic DNNs\footnote{i.e., DNNs that, for fixed parameters, define a deterministic mapping from input to output.} and becomes meaningful only when a mechanism for discarding information (e.g., noise) is integrated into the system. Such a noisy DNN framework was proposed in \cite{ICML_Info_flow2019}, where each neuron adds a small amount of Gaussian noise (i.i.d. across neurons) after applying the activation function. While the injection of noise alleviates the degeneracy of $I(X;T)$, the concatenation of Gaussian noises and nonlinearities makes this mutual information impossible to compute analytically or even evaluate numerically. Specifically, the distribution of $T$ (marginal or conditioned on $X$) is highly convoluted and thus the appropriate mode of operation becomes treating it as unknown, belonging to some nonparametric class.

Herein, we lay the groundwork for estimating $I(X;T)$ (or any other mutual information between layers) over real-world DNN classifiers. To achieve this, the estimation of $I(X;T)$ is reduced to the problem of differential entropy estimation under Gaussian convolutions described above. Specifically, in a noisy DNN each hidden layer can be written as $T=S+Z$, where $S$ is a deterministic function of the previous layer and $Z$ is a centered isotropic Gaussian vector. The DNN's generative model enables sampling $S$, while the distribution of $Z$ is known since the noise is injected by design. Estimating mutual information over noisy DNNs thus boils down to estimating $h(T) = h(S+Z)$ from samples of $S$, which is a main focus in this work.




\textbf{Outline:} The remainder of this paper is organized as follows. Section \ref{SEC:distances} analyzes the convergence of various statistical distances between $P\ast \Gauss$ and its Gaussian-smoothed empirical approximation. In Section \ref{SEC:Results} we set up the differential entropy estimation problem and state our main results. Section \ref{SEC:applications} presents applications of the considered estimation problem, focusing on mutual information estimation over DNNs. Simulation results are given in Section \ref{SEC:simulations}, and proofs are provided in Section \ref{SEC:proofs}. The main insights from this work and potential future directions are discussed in Section~\ref{SEC:summary}.

\textbf{Notation:} Logarithms are with respect to (w.r.t.) base $e$. For an integer $k\geq1$, we set $[k]\triangleq\big\{i\in\mathbb{Z}\big| 1\leq i \leq k\big\}$. $\|x\|$ is the Euclidean norm in $\mathbb{R}^d$, and $\mathrm{I}_d$ is the $d\times d$ identity matrix. We use $\mathbb{E}_P$ for an expectation w.r.t. a distribution $P$, omitting the subscript when $P$ is clear from the context. For a continuous $X\sim P$ with PDF $p$, we interchangeably use $h(X)$, $h(P)$ and $h(p)$ for its differential entropy. The $n$-fold product extension of $P$ is denoted by $P^{\otimes n}$. The convolution of two distributions $P$ and $Q$ on $\mathbb{R}^d$ is $(P\ast Q)(\mathcal{A})=\int\int\mathds{1}_{\mathcal{A}}(x+y)\dd P(x)\dd Q(y)$, where $\mathds{1}_\mathcal{A}$ is the indicator of the Borel set~$\mathcal{A}$. We use $\mathcal{N}_\sigma$ for the isotropic Gaussian measure of parameter $\sigma$, and denotes its PDF by $\gauss$. 

\section{Smooth Empirical Approximation}\label{SEC:distances}

This section studies the convergence rate of $\delta(\hat{P}_{S^n}\ast\Gauss,P\ast\Gauss)$ for different statistical distances $\delta(\cdot,\cdot)$, when $P$ is a $K$-subgaussian distribution, as defined next, and $\hat{P}_{S^n}$ is the empirical measure associated with $S^n\sim P^{\otimes n}$.

\begin{definition}[Subgaussian Distribution]\label{DEF:SG}
A $d$-dimensional distribution $P$ is $K$-subgaussian, for $K>0$, if $X\sim P$ satisfies
\begin{equation}
\E\Big[\exp\big(\alpha^T(X \mspace{-2mu}-\mspace{-2mu} \E X)\big)\Big] \leq \exp\big(0.5 K^2 \|\alpha\|^2\big),\quad \forall\alpha \in \RR^d.\label{EQ:SG}
\end{equation}
\end{definition}
In words, the above requires that every one-dimensional projection of $X$ be subgaussian in the traditional scalar sense. When $\big(X-\E X\big)\in\mathcal{B}(0,R)\triangleq\big\{x\in\RR^d\big|\|x\|\leq R\}$, \eqref{EQ:SG} holds with $K=R$.

When $\delta(\cdot, \cdot)$ is the TV or the 1-Wasserstein distance, the distance between the convolved distributions converges at the rate $n^{-1/2}$ for all $K$ and $d$. However, when $\delta(\cdot, \cdot)$ is the KL divergence or squared 2-Wasserstein distance (both are squared distances), convergence at rate $n^{-1}$ happens only when $K$ is sufficiently small or $P$ has bounded support. Interestingly, when $\delta(\cdot,\cdot)$ is the $\Chi$-divergence, two very different behaviors are observed. For low SNR (and in particular when $P$ has bounded support), $\E_{P^{\otimes n}}\Chi\left(\hat{P}_{S^n}\ast\Gauss\middle\|P\ast\Gauss\right)$ converges as $n^{-1}$. However, when SNR is high, we construct a subgaussian $P$ for which the expected $\Chi$-divergence is infinite.

Another way to summarize the results of this section is through the following curious dichotomy. This dichotomy highlights the central role of the $\Chi$-divergence in Gaussian-smoothed empirical approximation. To state it, let $Y=S+Z$, where $S\sim P$ and $Z\sim\Gauss$ are independent. Denote the joint distribution of $(S,Y)$ by $P_{S,Y}$, and let $P_S=P$ and $P_Y=P\ast\Gauss$ be their marginals. Setting $I_{\Chi}(S;Y)\triangleq\Chi\left(P_{S,Y}\middle\|P_S\otimes P_Y\right)$ as the $\Chi$ mutual information between $S$ and $Y$, we have:
\begin{enumerate}
\item If $P_S$ is $K$-subgaussian with $K<\frac{\sigma}{2}$ then $I_{\Chi}(S;Y)< \infty$. If $K> \sqrt{2}\sigma$ then $I_{\Chi}(S;Y)= \infty$ for some distributions $P_S$.
\item Assume $I_{\Chi}(S;Y) < \infty$. Then $\E_{P^{\otimes
n}}\delta(\hat{P}_{S^n}\ast\Gauss,P\ast\Gauss) = O\left(n^{-1}\right)$ if $\delta$ is the KL or $\Chi$-divergence. If $\delta$ is the TV or the 1-Wasserstein distance, the convergence rate is
$O(n^{-1/2})$. If, in addition, $P_S$ is subgaussian (with any constant), then the rate is also $O(n^{-1})$ if $\delta$ is the square of the 2-Wasserstein distance.
\item Assume $I_{\Chi}(S;Y) = \infty$. Then $\E_{P^{\otimes n}}\delta(\hat{P}_{S^n}\ast\Gauss,P\ast\Gauss) = \omega(n^{-1})$ if $\delta$ is the KL divergence or the squared 2-Wasserstein distance. If $\delta$ is the $\Chi$-divergence, it is infinite. 
\end{enumerate}
All the above are stated or immediately implied by the results to follow.


\subsection{1-Wasserstein Distance}\label{SUBSEC:1Wasserstein}

The 1-Wasserstein distance between $\mu$ and $\nu$ is given by $\mathsf{W}_1(\mu, \nu) \triangleq \inf \mathbb{E}\|X - Y\|$, where the infimum is taken over all couplings of $\mu$ and
$\nu$, i.e., joint distributions $P_{X,Y}$ whose marginals satisfy $P_X = \mu$ and $P_Y = \nu$. 

\begin{proposition}
[Smooth $\bm{\mathsf{W}_1}$ Approximation]\label{PROP:W1_parametric}
Fix $d\geq 1$, $\sigma>0$ and $K>0$. For any $K$-subgaussian distribution $P$, we have
\begin{equation}
\E_{P^{\otimes n}} \mathsf{W}_1(\hat{P}_{S^n} \ast\Gauss, P\ast\Gauss) \leq c^{(\mathsf{W_1})}_{\sigma,d,K} \frac{1}{\sqrt n},\label{EQ:W1_parametric_equation}
\end{equation}
where $c^{(\mathsf{W_1})}_{\sigma,d,K}$ is given in \eqref{EQ:W1_bound_constant}.
\end{proposition}

\begin{IEEEproof}
We can assume without loss of generality that $P_S$ has mean $0$. We start with the following upper bound \cite[Theorem~6.15]{villani2008optimal}:
\begin{equation}
\mathsf{W}_1(\hat{P}_{S^n}\ast\Gauss, P\ast\Gauss) \leq \int_{\RR^d} \|z\| \big|r_{S^n}(z) - q(z)\big| \dd z,
\end{equation}
where $r_{S^n}$ and $q$ are the densities associated with $\hat{P}_{S^n}\ast\Gauss$ and $P\ast\Gauss$, respectively. This inequality follows by coupling $\hat{P}_{S^n}\ast\Gauss$ and $P\ast\Gauss$ via the maximal TV-coupling. 

Let $f_a:\RR^d \to \RR$ be the PDF of $\mathcal{N}\left(0,\frac{1}{2a} \mathrm{I}_d\right)$, for $a > 0$ specified later. The Cauchy-Schwarz inequality implies 
\begin{align}
&\mathbb{E}_{P^{\otimes n}}\int_{\RR^d} \|z\| \big|r_{S^n}(z) - q(z)\big| \dd z\label{EQ:W1_CS_bound}\\\nonumber
&\leq \left(\int_{\RR^d} \|z\|^2f_a(z)\dd z\right)^{\frac{1}{2}}\!\left(\mathbb{E}_{P^{\otimes n}}\!\int_{\RR^d}\! \frac{\big(q(z)-r_{S^n}(z)\big)^2}{f_a(z)} \dd z\right)^{\frac{1}{2}}\!\!.
\end{align}
The first term is the expected  squared Euclidean norm of $\mathcal{N}\left(0,\frac{1}{2a} \mathrm{I}_d\right)$, which equals $\frac{d}{2a}$. 

For the second integral, note that $r_{S^n}(z)$ is a sum of i.i.d. terms with expectation $q(z)$. This implies 
\begin{align*}
    \mathbb{E}_{P^{\otimes n}}\big(q(z)-&r_{S^n}(z)\big)^2
    =\var_{P^{\otimes n}}\left(\frac{1}{n}\sum_{i=1}^n\gauss(z-S_i)\right)\\
    &=\frac{1}{n}\var_P\big(\gauss(z-S)\big)\leq\frac{c_1^2}{n}\E_P e^{-\frac{1}{\sigma^2} \|z - S\|^2},
\end{align*}
where $c_1 = (2 \pi \sigma^2)^{-d/2}$. Consequently, we have
\begin{equation}
\int_{\RR^d}\E_{P^{\otimes n}} \frac{\big(q(z)-r_{S^n}(z)\big)^2}{f_a(z)} \dd z \leq \frac{c_1}{n2^{d/2}} \E \frac{1}{f_a(S+ Z/\sqrt{2} )},\label{EQ:W1_1/f_bound}
\end{equation}
where $S\sim P$ and $Z\sim\cN(0, \sigma^2\mathrm{I_d})$ are independent.

Setting $c_2\triangleq\left(\frac{\pi}{a}\right)^{\frac{d}{2}}$, we have $\big(f_a(z)\big)^{-1}= c_2 \exp\big(a \|z\|^2\big)$. Since $S$ is $K$-subgaussian and $Z$ is $\sigma$-subgaussian, $S + Z/\sqrt{2}$ is $(K+ \sigma/\sqrt{2})$-subgaussian. Following \eqref{EQ:W1_1/f_bound}, for any $0<a<\frac{1}{2(K + \sigma/\sqrt{2})^2}$, we have~\cite[Remark 2.3]{hsu2012tail}
\begin{align*}
&\frac{c_1}{n2^{d/2}}\E \frac{1}{f_a(S + Z/\sqrt{2})}=\frac{c_1 c_2}{n2^{d/2}} \E \exp\Big(a \big\|S + Z/\sqrt{2}\big\|^2\Big)\\
&\leq \frac{c_1 c_2}{n2^{d/2}} \exp\left(\big(K+\sigma/ \sqrt 2\big)^2 a d+\frac{(K+\sigma/ \sqrt 2)^4 a^2d}{1 -2(K+\sigma/ \sqrt 2)^2 a}\right),\numberthis\label{EQ:W1_int2_bound}
\end{align*}
where the last inequality uses the subgaussianity of  $S+Z/\sqrt{2}$ and Definition \ref{DEF:SG}. Setting $a=\frac{1}{4(K+\sigma/\sqrt{2})^2}$, we combine \eqref{EQ:W1_CS_bound}-\eqref{EQ:W1_int2_bound} to obtain the result
\begin{align}
\E_{P^{\otimes n}}\mathsf{W}_1(\hat{P}_{S^n}\ast\Gauss &, P\ast\Gauss)\nonumber\\ &\leq \sigma\sqrt{2d}\left(\frac{1}{\sqrt 2}+\frac{K}{\sigma}\right)^{\frac{d}{2}+1}e^{\frac{3d}{16}}\frac{1}{\sqrt{n}}.\label{EQ:W1_bound_constant}
\end{align}
\end{IEEEproof}

\begin{remark}[Smooth $\bm{\mathsf{W}_1}$ for Bounded Support]
A better constant is attainable if attention is restricted to the bounded support case. It was shown in \cite{Weed_bdd_Support2018} that analyzing $\E_{P^{\otimes n}}\mathsf{W}_1(\hat{P}_{S^n}\ast\Gauss, P\ast\Gauss)$ directly for $\supp(P)\subseteq[-1,1]^d$, one can obtain the bound $\frac{2^{d+2}\sqrt{d}}{\min\{1,\sigma^{d}\}}n^{-1/2}$.
\end{remark}

\subsection{Total Variation Distance}\label{SUBSEC:TV}

The TV distance between $\mu$ and $\nu$ is $\|\mu-\nu\|_{\mathsf{TV}} \triangleq \sup_{A\in\mathcal{F}}|\mu(A)-\nu(A)|$, where $\mathcal{F}$ is the sigma-algebra. When $\mu$ and $\nu$ have densities, say $f$ and $g$, respectively, the TV distance reduces to $\frac{1}{2}\int|f(z)-g(z)|\dd z$.

\begin{proposition}
[Smooth Total Variation Approximation]\label{PROP:TV_parametric}
Fix $d\geq 1$, $\sigma>0$ and $K>0$. For any $K$-subgaussian distribution $P$, we have
\begin{equation}
\E_{P^{\otimes n}} \left\|\hat{P}_{S^n} \ast\Gauss-P\ast\Gauss\right\|_{\mathsf{TV}} \leq c^{(\mathsf{TV})}_{\sigma,d,K} \frac{1}{\sqrt n},\label{EQ:TV_parametric_equation}
\end{equation}
where $c^{(\mathsf{TV})}_{\sigma,d,K}$ is given in \eqref{EQ:TV_bound_constant}.
\end{proposition}

\begin{IEEEproof}
Noting that $\E_{P^{\otimes n}} \left\|\hat{P}_{S^n} \ast\Gauss-P\ast\Gauss\right\|_{\mathsf{TV}}=\frac{1}{2}\E_{P^{\otimes n}} \int\big|r_{S^n}(z)-q(z)\big|\dd z$, we may apply the Cauchy-Schwarz inequality similarly to \eqref{EQ:W1_CS_bound}. The only difference now is that the first integral sums up to 1 (rather than being a Gaussian moment). Repeating steps \eqref{EQ:W1_1/f_bound}-\eqref{EQ:W1_int2_bound} we obtain
\begin{equation}
\E_{P^{\otimes n}} \left\|\hat{P}_{S^n} \ast\Gauss-P\ast\Gauss\right\|_{\mathsf{TV}} \leq \left(\frac 1{\sqrt 2}+\frac{K}{\sigma}\right)^{\frac{d}{2}}e^{\frac{3d}{16}}\frac{1}{\sqrt{n}},\label{EQ:TV_bound_constant}
\end{equation}
as desired.\end{IEEEproof}


\subsection{$\Chi$-Divergence}\label{SUBSEC:Chi}

The $\Chi$-divergence $\Chi(\mu\|\nu)\triangleq\int \left(\frac{\dd\mu}{\dd\nu}-1\right)^2\dd\nu$ presents perhaps the most surprising behavior of all the considered distances. When the signal-to-noise ratio (SNR) $\frac{K}{\sigma}<\frac{1}{2}$, we prove that $\E_{P^{\otimes n}} \Chi\left(\hat{P}_{S^n}\ast\Gauss\middle\|P\ast\Gauss\right)$ converges as $n^{-1}$ for all dimensions. However, if $K\geq\sqrt{2} \sigma$, then there exists $K$-subgaussian distributions $P$ such that $\E_{P^{\otimes n}} \Chi\left(\hat{P}_{S^n}\ast\Gauss\middle\|P\ast\Gauss\right)=\infty$ even in $d=1$. 
Our results rely on the following identity.
\begin{lemma}[$\bm{\Chi}$-Divergence and Mutual Information]\label{lem:chi2-information}
Let $S \sim P$ and $Y = S + Z$, with $Z \sim \mathcal{N}_\sigma$ independent of~$S$.
Then
\begin{equation*}
\E_{P^{\otimes n}} \Chi\left(\hat{P}_{S^n} \ast\Gauss\middle\|P\ast\Gauss\right) = \frac 1n I_{\Chi}(S; Y)\,.
\end{equation*}
\end{lemma}
\begin{IEEEproof}
Recall that $r_{S^n}(z)$ is a sum of i.i.d.\ terms with expectation $q(z)$.
This yields
\begin{align*}
    \E_{P^{\otimes n}} &\Chi\left(\hat{P}_{S^n} \ast\Gauss\middle\|P\ast\Gauss\right)\\&=\E_{P^{\otimes n}} \int_{\RR^d} \frac{\big(r_{S^n}(z) - q(z)\big)^2}{q(z)} \dd z\\
    &= \frac{1}{n}\left( \int_{\RR^d}\E_P \frac{(\gauss(z - S) - q(z))^2}{q(z)}\dd z\right)\\
    &=\frac{1}{n}I_{\Chi}(S; Y)\,.\numberthis\label{EQ:Chi_empirical_to_MI}
\end{align*}
\end{IEEEproof}

Lemma~\ref{lem:chi2-information} implies that $\E_{P^{\otimes n}} \Chi\left(\hat{P}_{S^n} \ast\Gauss\middle\|P\ast\Gauss\right) = O(n^{-1})$ if and only if $I_{\Chi}(S; Y) < \infty$.
When $I_{\Chi}(S; Y) = \infty$, $\E_{P^{\otimes n}} \Chi\left(\hat{P}_{S^n} \ast\Gauss\middle\|P\ast\Gauss\right)$ diverges for all $n$. It therefore suffices to examine conditions under which $I_{\Chi}(S; Y)$ is finite.


\vspace{3mm}
\subsubsection{\underline{Convergence at Low SNR and Bounded Support}} We start by stating and proving the convergence results.

\begin{proposition}[Smooth $\bm{\Chi}$ Approximation]\label{PROP:Chi_parametric}
Fix $d\geq 1$ and $\sigma>0$. If $P$ is $K$-subgaussian with $K<\frac{\sigma}{2}$, then $I_{\Chi}(S;Y)\leq c^{(\Chi)}_{\sigma,d,K}<\infty$, where $c^{(\Chi)}_{\sigma,d,K}$ is given in \eqref{EQ:Chi_bound_constant}.
Consequently,
\begin{equation*}
\E_{P^{\otimes n}} \Chi\left(\hat{P}_{S^n} \ast\Gauss\middle\|P\ast\Gauss\right) \leq c^{(\Chi)}_{\sigma,d,K} \frac 1n\,.
\end{equation*}
\end{proposition}

\begin{IEEEproof}
Denote by $\mathcal N(x, \sigma^2 \mathrm{I}_d)$ an isotropic Gaussian of entrywise variance $\sigma^2$ centered at $x$.
Then by the convexity of the $\Chi$-divergence,
\begin{align*}
    I_{\Chi}(S; Y)
    &= \E_P\Chi\left(\mathcal{N}(S,\sigma^2\mathrm{I}_d)\middle\|\E_P\mathcal{N}(\tilde{S},\sigma^2\mathrm{I}_d)\right) \\& \leq \E_{P^{\otimes 2}} \Chi\left(\mathcal{N}(S,\sigma^2\mathrm{I}_d)\middle\| \mathcal{N}(\tilde{S},\sigma^2\mathrm{I}_d)\right) \\
    & = \E_{P^{\otimes 2}} e^{\frac{1}{\sigma^2}\|S-\tilde{S}\|^2}\,,\numberthis\label{EQ:Chi_bound_last}
\end{align*}
where the last equality follows from the closed form expression for the $\Chi$-divergence between Gaussians~\cite{nielsen2014chi}.

Since $S-\tilde{S}$ is $\sqrt 2 K$-subgaussian, the RHS above converges if $K<\frac{\sigma}{2}$ and gives the following bound~\cite[Remark 2.3]{hsu2012tail}
\begin{equation}
     I_{\Chi}(S; Y)\leq\exp\left(2d\left(\frac{K}{\sigma}\right)^2\frac{\sigma^2-2K^2}{\sigma^2-4K^2}\right)\,.\label{EQ:Chi_bound_constant}
\end{equation}
     The second claim follows from Lemma~\ref{lem:chi2-information}.
\end{IEEEproof}

The proof of Proposition \ref{PROP:Chi_parametric} immediately implies $\Chi$ convergence for any compactly supported $P$.

\begin{corollary}[Smooth $\bm{\Chi}$ for Bounded Support]\label{CORR:Chi_parametric_bdd}
If $P$ has a bounded support with diameter $D\triangleq\sup_{x,y\in\supp(P)}\|x-y\|$, then $I_{\Chi}(S; Y) \leq \exp\left(\frac{D^2}{\sigma^2}\right)$.
Consequently,
\begin{equation}
\E_{P^{\otimes n}} \Chi\left(\hat{P}_{S^n} \ast\Gauss\middle\|P\ast\Gauss\right) \leq \exp\left(\frac{D^2}{\sigma^2}\right)\cdot\frac{1}{n},\label{EQ:Chi_parametric_equation_bdd}
\end{equation}
\end{corollary}

\begin{IEEEproof}
Follows by inserting $\|S-\tilde{S}\|\leq D$ into \eqref{EQ:Chi_bound_last}.
\end{IEEEproof}


\ \\
\subsubsection{\underline{Diverging Example}}
This section shows that for $K \geq \sqrt{2}\sigma$, there exist $K$-subgaussian distributions $P$ for which $I_{\Chi}(S; Y) = \infty$.

Let $d=1$ and without loss of generality assume $\sigma = 1$.
Furthermore, for simplicity of the proof we set $K = \sqrt{2}\sigma$, since the resulting counterexample will apply for any $K \geq \sqrt{2}\sigma$ (recalling that any $\sqrt{2}\sigma$-subgaussian distribution is $K$-subgaussian for any $K \geq \sqrt{2}\sigma$).

Let $\epsilon = \frac{1}{2K^2} = \frac{1}{4}$ and define the sequence $\{r_k\}_{k=0}^\infty$ by $r_0 = 0$, $r_1 = 1$ and $r_k = \frac{r_{k-1}}{1-\sqrt{2\epsilon}}$, for $k\geq 2$. Let $P$ be  discrete distribution with $\supp(P)=\{r_k\}_{k=0}^\infty$ given by
\begin{subequations}\label{eq:Pcounter}
\begin{equation}
P = \sum_{k = 0}^\infty p_k \delta_{r_k},
\end{equation}
where $\delta_x$ is the Dirac measure at $x$. We make $P$ $K$-subgaussian by setting
\begin{equation}
p_k = \left\{\begin{array}{ll} 2\sqrt{\frac{{\epsilon}}{{\pi}}} e^{-\epsilon r_k^2}& k \geq 1\\ 1 - \sum_{k = 1}^\infty p_k & k = 0.\end{array} \right. 
\end{equation}
\end{subequations}%
Note then that since $\min_{k \geq 1} |r_k - r_{k-1}| = 1$ and $p_k=2\varphi_K(r_k)$, we get that $\sum_{k = 1}^\infty p_k < 1$; the remainder of the probability is allocated to $r_0=0$. As stated in the next proposition, $I_{\Chi}(S;Y)$ diverges when $S\sim P$ as constructed above (which, in turn implies that the $\Chi$ smoothed empirical approximation is also infinite). This stands in contrast to the classic KL mutual information, which is always finite over an AWGN channel for inputs with a bounded second moment. 

\begin{proposition}[$\bm{\Chi}$ Diverging Example]\label{prop:div}
For $P$ as in \eqref{eq:Pcounter} and $\epsilon = 1/4$, we have that \begin{equation}
    I_{\Chi}(S;Y)=\infty\,.
\end{equation}
Consequently
\begin{equation*}
\E_{P^{\otimes n}} \Chi\left(\hat{P}_{S^n} \ast\mathcal{N}_1\middle\|P\ast\mathcal{N}_1\right)= \infty,\quad \forall n\in\mathbb{N}.
\end{equation*}
\end{proposition}
The proof is given in Appendix \ref{app:counterex}.
Intuitively, the constructed $P$ has infinitely many atoms at sufficiently large distance from each other such that the tail contribution at $r_k$ from any $j\neq k$ component of the mixture $P\ast\mathcal{N}_1$ is negligible. Note that we grow $r_k$ exponentially to counter the exponentially shrinking $p_k$ weights. Since $\supp(P)=\mathbb{N}\cup\{0\}$, for any finite $n$ there are infinitely many atoms which were not sampled in $S^n$. Since they are sufficiently well-separated, each of these unsampled atoms contributes a constant value to the considered $\Chi$-divergence, which consequently becomes infinite.


\subsection{2-Wasserstein Distance and Kullback-Leibler Divergence}
One can leverage the above $\Chi$ results to obtain analogous bounds for KL divergence and for $\mathsf{W}_2$. The squared 2-Wasserstein distance between $\mu$ and $\nu$ is $ \mathsf{W}^2_2(\mu,\nu) \triangleq \inf \E\|X - Y\|^2$, where the infimum is taken over all couplings of $\mu$ and $\nu$. The KL divergence is given by $\mathsf{D}_{\mathsf{KL}}(\mu\|\nu)\triangleq \int\log\left(\frac{\dd\mu}{\dd\nu}\right)\dd\mu$.

The behavior of the 2-Wasserstein distance and KL divergence are governed by $I_{\Chi}(S;Y)$. If $I_{\Chi}(S; Y) < \infty$, then the KL divergence converges at the rate $O(n^{-1})$. If additionally $P$ is $K$-subgaussian (for any $K < \infty$), then the squared 2-Wasserstein distance also converges as $O(n^{-1})$. On the other hand, if $I_{\Chi}(S; Y) = \infty$, then both the KL divergence and the squared 2-Wasserstein distance are $\omega(n^{-1})$ in expectation.
%
\vspace{3mm}
 
\subsubsection{\underline{Parametric convergence when $I_{\Chi}(S; Y) < \infty$}}
Our bounds on the $\Chi$-divergence immediately imply analogous bounds for the KL divergence.
\begin{proposition}[Smooth KL Divergence Approximation]\label{prop:kl_convergence}
Fix $d \geq 1$ and $\sigma > 0$. We have
\begin{equation}\label{eq:kl_subg}
\E_{P^{\otimes n}} \mathsf{D}_{\mathsf{KL}} \left(\hat{P}_{S^n} \ast\Gauss\middle\|P\ast\Gauss\right) \leq \frac 1 {n} I_{\Chi}(S; Y) \,.
\end{equation}
In particular, if $P$ is $K$-subgaussian with $K<\frac{\sigma}{2}$, then
\begin{subequations}
\begin{equation}
    \E_{P^{\otimes n}} \mathsf{D}_{\mathsf{KL}} \left(\hat{P}_{S^n} \ast\Gauss\middle\|P\ast\Gauss\right) \leq c_{\sigma, d,K}^{(\Chi)}  \frac 1n
\end{equation}
and if $P$ is supported on a set of diameter $D$, then 
\begin{equation}
    \E_{P^{\otimes n}} \mathsf{D}_{\mathsf{KL}} \left(\hat{P}_{S^n} \ast\Gauss\middle\|P\ast\Gauss\right) \leq \exp\left(\frac{D^2}{\sigma^2}\right)\cdot\frac{1}{n}.
\end{equation}
\end{subequations}
\end{proposition}
\begin{IEEEproof}
The first claim follows directly from Lemma~\ref{lem:chi2-information} because $D_{\mathsf{KL}}(\mu\|\nu)\leq \log\big(1+\Chi(\mu\|\nu)\big)$ for any two probability measures $\mu$ and $\nu$.
Proposition~\ref{PROP:Chi_parametric} and Corollary~\ref{CORR:Chi_parametric_bdd} then imply the subsequent claims.
\end{IEEEproof}

To obtain bounds for the 2-Wasserstein distance, we leverage a transport-entropy inequality that connects KL divergence and $\mathsf{W}_2$. We have the following.
\begin{proposition}[Smoothed $\bm{\mathsf{W}_2^2}$ Approximation]\label{prop:Wass2}
Let $d \geq 1$, $\sigma > 0$ and $K < \infty$. If $P$ is a $K$-subgaussian distribution and $I_{\Chi}(S; Y) < \infty$, then
\begin{equation}
\E_{P^{\otimes n}} \mathsf{W}_2^2(\hat{P}_{S^n} \ast\Gauss, P\ast\Gauss) = O\left(\frac 1 n\right)\,.
\end{equation}
In particular, this holds if $K < \frac \sigma 2$.
More explicitly, if $P$ is supported on a set of diameter $D$, then
\begin{equation}
    \E_{P^{\otimes n}} \mathsf{W}_2^2(\hat{P}_{S^n} \ast\Gauss, P\ast\Gauss) \leq c_{\sigma, d, D}^{(\mathsf{W}_2)} \cdot \frac 1 n\,,
\end{equation}
where $c_{\sigma, d, D}^{(\mathsf{W}_2)} $ is given in~\eqref{EQ:w2_bound}.
\end{proposition}
\begin{IEEEproof}
The subgaussianity of $P$ implies~\cite[Lemma~5.5]{vershynin2010introduction} that  $\E_P e^{\varepsilon \|S\|^2} < \infty$ for $\varepsilon > 0$ sufficiently small.
Therefore, \cite[Theorem 1.2]{WanWan16} implies that $P \ast \Gauss$ satisfies a log-Sobolev inequality with some constant $C_{P, \sigma}$, depending on $P$ and $\sigma$. This further means \cite{OttVil00,BobGenLed01} that $P \ast \Gauss$ satisfies the transport-entropy inequality
\begin{equation}
        \mathsf{W}_2^2(Q, P \ast \Gauss) \leq C_{P, \sigma} \mathsf{D}_{\mathsf{KL}} (Q\|P\ast\Gauss) 
\end{equation}
for all probability measures $Q$.
Combining this inequality with Proposition~\ref{prop:kl_convergence} yields the first claim.

If $P$ is supported on a set of diameter $D$, we have the following more explicit bound:
\begin{equation}
    \mathsf{W}_2^2(Q, P \ast \Gauss) \leq c' \sqrt{d} \sigma^2\left(1+\frac{D^2}{4 \sigma^2}\right)e^{\frac{D^2}{\sigma^2}} \mathsf{D}_{\mathsf{KL}} (Q\|P\ast\Gauss) 
\end{equation}
for an absolute constant $c'$ and any probability measure $Q$ on $\RR^d$.
Applying Proposition~\ref{prop:kl_convergence} yields
\begin{align}
    \E_{P^{\otimes n}} \mathsf{W}_2^2&(\hat{P}_{S^n} \ast\Gauss, P\ast\Gauss) \nonumber\\&\leq c'\frac{\sqrt{d} \sigma^2}{n}\left(1+\frac{D^2}{4 \sigma^2}\right) \exp\left(\frac{2D^2}{\sigma^2}\right) .\label{EQ:w2_bound}
\end{align}
\end{IEEEproof}

\begin{remark}[$\bm{\mathsf{W}_2^2}$ Speedup in One Dimension]
The convergence of (unsmoothed) empirical measures in the squared 2-Wasserstein distance suffers from the curse of dimensionality, converging at rate $n^{-2/d}$ when $d$ is large \cite{boissard2014mean,dudley1969speed}. Proposition~\ref{prop:Wass2}, however, shows that when smoothed with Gaussian kernels the convergence rate improves to $n^{-1}$ for all $d$ and low SNR ($K<\frac{\sigma}{2}$). Interestingly, the Gaussian smoothing speeds up the convergence rate even for $d=1$. For instance, Theorem 7.11 from \cite{bobkov2014one} shows that $\mathbb{E}_{P^{\otimes n}}\mathsf{W}^2_2(\hat{P}_{S^n},P)\asymp n^{-1/2}$ whenever the support of $P$ is not an interval in $\mathbb{R}$, which is slower than the $n^{-1}$ attained under Gaussian smoothing. Even for the canonical case when $P$ is Gaussian, Corollary 6.14 from \cite{bobkov2014one} shows that $\mathbb{E}_{P^{\otimes n}}\mathsf{W}^2_2(\hat{P}_{S^n},P)\asymp {\frac{\log\log n}{n}}$.


\end{remark}

\subsubsection{\underline{Slower convergence when $I_{\Chi}(S; Y) = \infty$}}
Unlike the $\Chi$-divergence, it is easy to see that the 2-Wasserstein distance and KL divergence between the convolved measures are always finite when $P$ is subgaussian.
However, when $I_{\Chi}(S; Y) = \infty$, the rate of convergence of the KL divergence and the squared 2-Wasserstein distance is strictly slower than parametric.

\begin{proposition}[Smooth KL Divergence vs. Parametric]\label{prop:kl_slow}
If $I_{\Chi}(S;Y)=\infty$, for $S\sim P$ and $Y=S+Z$, with $Z\sim\mathcal{N}_\sigma$ independent of $S$, then
\begin{equation}
\E_{P^{\otimes n}} \mathsf{D}_{\mathsf{KL}} \left(\hat{P}_{S^n} \ast\Gauss\middle\|P\ast\Gauss\right) = \omega\left(\frac 1n \right).
\end{equation}
For examples of $K$-subgaussian $S\sim P$ distributions with 
$I_{\Chi}(S;Y)=\infty$ see Proposition \ref{prop:div}.
\end{proposition}
\begin{IEEEproof}
By rescaling, we assume that $\sigma = 1$. Let $S^n\sim P^{\otimes n}$, $Z\sim \mathcal{N}_1$ and $W \sim \mathsf{Unif}([n])$ be independent random variable. Defining $V = S_W + Z$, 
the proof of Lemma~\ref{LEMMA:SP_bias_MI} in Appendix~\ref{APPEN:SP_bias_MI_proof} establishes that $V$ has law $P\ast\mathcal{N}_1$ and that the conditional distribution of $V$ given $S^n = s^n$ is $\hat{P}_{s^n} \ast\mathcal{N}_1$.
This implies
\begin{align*}
\E_{P^{\otimes n}} \mathsf{D}_{\mathsf{KL}}& \left(\hat{P}_{S^n}\ast\mathcal{N}_1\middle\|P\ast\mathcal{N}_1\right) = I(S^n; V)\\
&\stackrel{(a)}= \sum_{i=1}^n h(S_i|S^{i-1}) - h(S_i|V,S^{i-1})\\
&\stackrel{(b)}\geq \sum_{i=1}^n h(S_i) - h(S_i|V)\\
&\stackrel{(c)}= n I(S_1; V),\numberthis
\end{align*}
where (a) uses $I(S^n;V)=h(S^n)-h(S^n|V)$ and the entropy chain rule, (b) is since $S_i$ are independent (first term) and because conditioning cannot increase entropy (second term), while (c) follows because $S_i$ are identically distributed and since $P_{V|S_i}$ does not depend on $i$.

Conditioned on $S_1 = s_1$, the random variable $V$ has law $\frac 1 n \delta_{s_1} \ast \Gauss + \frac{n-1}{n} P \ast \mathcal{N}_1$.
Therefore
\begin{equation}
I(S_1; V) = \E_{P} \mathsf{D}_{\mathsf{KL}} \left(\frac 1 n \delta_{S_1} \ast \mathcal{N}_1 + \frac{n-1}{n} P \ast \mathcal{N}_1 \middle\| P \ast \mathcal{N}_1\right).
\end{equation}

Let $S \sim P$ and $Y = S + Z$, and denote by $P_S$ and $P_Y$ the distributions of $S$ and $Y$, respectively.
By Fatou's~lemma,
\begin{align*}
&\liminf_{n \to \infty} n \E_{P^{\otimes n}} \mathsf{D}_{\mathsf{KL}} \left(\hat{P}_{S^n} \ast\mathcal{N}_1\middle\|P\ast\mathcal{N}_1\right) \\& \geq \liminf_{n \to \infty} \E_{P_S}  n^2  \mathsf{D}_{\mathsf{KL}} \left(\frac 1 n \delta_{S} \ast \mathcal{N}_1 + \frac{n-1}{n} P \ast \mathcal{N}_1 \middle \| P \ast \mathcal{N}_1\right) \\
& = \liminf_{n \to \infty}  \E_{P_S}  n^2 \mathsf{D}_{\mathsf{KL}} \left(\frac 1 n P_{Y | S} + \frac{n-1}{n} P_Y \middle \| P_Y \right) \\
& = \liminf_{n \to \infty}  n^2 \mathsf{D}_{\mathsf{KL}} \left(\frac 1 n P_{S,Y} + \frac{n-1}{n} P_S \otimes P_Y \middle \| P_S \otimes P_Y \right) \\
&\qquad \stackrel{(a)}\geq \E_{P_S \otimes P_Y} \left(\frac{\mathrm{d}P_{S,Y}}{\mathrm{d}P_{S} \otimes P_Y} - 1\right)^2 \\
&\qquad = I_{\Chi}(S; Y)\\
&\qquad= \infty\numberthis,
\end{align*}
where (a) follows by \cite[Proposition 4.2]{PolyWu-LecNotes}. We conclude that $\E_{P^{\otimes n}} \mathsf{D}_{\mathsf{KL}} \left(\hat{P}_{S^n} \ast\mathcal{N}_1\middle\|P\ast\mathcal{N}_1\right) = \omega(n^{-1})$.
\end{IEEEproof}
\ \\
We likewise obtain an analogous claim for $\mathsf{W}_2$, showing that the squared 2-Wasserstein distance converges as $\omega(n^{-1})$ when smoothed by any Gaussian with strictly smaller variance.

\begin{corollary}[Smooth $\bm{\mathsf{W}_2^2}$ Slower than Parametric]
If $I_{\Chi}(S;Y)=\infty$, for $S\sim P$ and $Y=S+Z$, with $Z\sim\mathcal{N}_\sigma$ independent of $S$, then for any $\tau < \sigma$,
\begin{equation}
\E_{P^{\otimes n}} \mathsf{W}_2^2(\hat{P}_{S^n} \ast\mathcal{N}_\tau, P\ast\mathcal{N}_\tau) = \omega \left(\frac 1n \right)\,.
\end{equation}
For examples of $K$-subgaussian $S\sim P$ distributions with 
$I_{\Chi}(S;Y)=\infty$ see Proposition \ref{prop:div}.
\end{corollary}
\begin{IEEEproof}
We assume as above that $\sigma = 1$. 
Let $S^n \sim P^{\otimes n}$, define $T_i \triangleq \frac{K}{\sqrt 2} S_i$, and Let $\lambda \triangleq \sqrt{1 - \tau^{2}} > 0$.
If $P_{R_n, R}$ is any coupling between $\hat P_{S^n} \ast \mathcal{N}_{\tau}$ and $P \ast \mathcal{N}_{\tau}$, then the joint convexity of the KL divergence implies that
\begin{align*}
\mathsf{D}_{\mathsf{KL}} &\left(\hat{P}_{S^n} \ast\mathcal{N}_\sigma\middle\|P\ast\mathcal{N}_\sigma\right) \\& = \mathsf{D}_{\mathsf{KL}} \left(\E_{P_{R_n}}\mathcal{N}(R_n,\lambda^2 \mathrm{I}_d) \middle \| \E_{P_R}\mathcal{N}(R,\lambda^2 \mathrm{I}_d)\right) \\
&  \leq \E_{P_{R_n, R}} \mathsf{D}_{\mathsf{KL}} \left(\mathcal{N}(R_n,\lambda^2 \mathrm{I}_d) \middle\| \mathcal{N}(R,\lambda^2 \mathrm{I}_d)\right) \\
& = \E_{P_{R_n, R}} \frac{1}{2\lambda^2} \|R_n - R\|^2\numberthis,
\end{align*}
where the last step uses the explicit expression for the KL divergence between isotropic Gaussians. Taking an infimum over all valid couplings yields
\begin{equation}
\mathsf{D}_{\mathsf{KL}} \left(\hat{P}_{S^n} \ast\mathcal{N}_1\middle\|P\ast\mathcal{N}_1\right) \leq \frac{1}{2 \lambda^2} \mathsf{W}_2^2(\hat{P}_{S^n} \ast\mathcal{N}_\tau, P\ast\mathcal{N}_\tau)\,.
\end{equation}
Proposition~\ref{prop:kl_slow} then implies
\begin{align*}
\E_{P^{\otimes n}} \mathsf{W}_2^2(\hat{P}_{S^n} \ast\mathcal{N}_\tau,& P\ast\mathcal{N}_\tau) \\&\geq 2 \lambda^2 \E_{P^{\otimes S^n}} \mathsf{D}_{\mathsf{KL}} \left(\hat{P}_{S^n}\ast\mathcal{N}_1\middle\|P\ast\mathcal{N}_1\right) \\&= \omega\left(\frac 1n \right),
\end{align*}
which concludes the proof.
\end{IEEEproof}

\subsection{Open Questions}\label{SUBSEC:open_questions}
We list here some open questions that remain unanswered by the above. First, recall we focused on bounds of the form:
	$$ \EE_{P^{\otimes n}}\delta\left(\hat P_{S^n} \ast \mathcal{N}_\sigma, P \ast \mathcal{N}_\sigma\right) \le C(d) \frac{1}{n}\,.$$
What is the correct dependence of $C(d)$ on dimension? For $\delta=\mathsf{W}_1$ we proved a bound with $C(d) =e^{O(d)}$ for all subgaussian $P$. Similarly, for $\delta=\mathsf{W}_2^2$ we have shown bounds with $C(d)=e^{O(d)}$ (for small-variance subgaussian $P$) and $C(d)=\sqrt{d}e^{O(D^2)}$ (for $P$ supported on a set of diameter $D$). What is the sharp dependence on dimension? Does it change as a function of the subgaussian constant?

A second, and perhaps more interesting, direction is to understand the rate of convergence of $\mathsf{W}_2^2$ in cases when~it is $\omega(n^{-1})$. A proof similar to Proposition~\ref{PROP:W1_parametric} can be used to show that for subgaussian $P$ (with any constant), we have
	$$ \EE_{P^{\otimes n}}\mathsf{W}_2^2\left(\hat P_{S^n} \ast \mathcal{N}_\sigma, P \ast \mathcal{N}_\sigma\right) = O\left(
		\frac{1}{\sqrt{n}}\right)\,.$$
Is this ever sharp? What rates are possible in the range between $\omega(n^{-1})$ and $O(n^{-1/2})$?

Heuristically, one may think that since $\mathsf{W}_1$ converges as $n^{-1/2}$, then some truncation argument should be able to recover $\mathsf{W}_2^2 \lesssim {\frac{\mathrm{polylog}(n)}{n}}$. Rigorizing this reasoning requires, however, analyzing the distance distribution between $A\sim \hat P_{S^n} \ast \Gauss$ and $B\sim P \ast \Gauss$ under the optimal $\mathsf{W}_1$-coupling. The TV-coupling that was used in Proposition~\ref{PROP:W1_parametric} will not work here because under it we have $\PP\big(\|A-B\|>\Omega(1)\big)=\Omega(n^{-1/2})$, which results in the $n^{-1/2}$ rate for~$\mathsf{W}_2^2$.

Finally, as we saw, the finiteness of $I_{\Chi}(S;Y)$ is a sufficient condition for many of the above empirical measure convergence results. When $S\sim P$ is $K$-subgaussian with $K<\frac{\sigma}{2}$, Proposition \ref{PROP:Chi_parametric} shows that $I_{\Chi}(S;Y)<\infty$ always holds. However, for $K\geq\sqrt{2}\sigma$, there exist $K$-subgaussian distribution for which $I_{\Chi}(S;Y)=\infty$ (Proposition \ref{prop:div}). Characterizing the sharp threshold at which $I_{\Chi}(S;Y)$ may diverge is another open question.


\section{Differential Entropy Estimation}\label{SEC:Results}


Our main application of the Gaussian smoothed empirical approximation questions is the estimation of $h(P\ast\mathcal{N}_\sigma)$, based on samples $S^n\sim P^{\otimes n}$ and knowledge of $\sigma$. The $d$-dimensional distribution $P$ is unknown and belongs to some nonparametric class. We first consider the class $\mathcal{F}_d$ of all distributions $P$ with $\supp(P)\subseteq[-1,1]^d$.\footnote{One may consider any other class of compactly supported distributions.} The second class of interest is $\mathcal{F}^{(\mathsf{SG})}_{d,K}$, which contains all $K$-subgaussian distributions (see Definition \ref{DEF:SG}). 

\subsection{Lower Bounds on Risk}

The sample complexity for estimating $h(P\ast\Gauss)$ over the class $\mathcal{F}_d$ is $n^\star(\eta,\sigma,\mathcal{F}_d)\triangleq\min\big\{n\in\NN: \mathcal{R}^\star(n,\sigma,\mathcal{F}_d)\leq \eta\big\}$, where $\mathcal{R}^\star(n,\sigma,\mathcal{F}_d)$ is defined in \eqref{EQ:risk}. As claimed next, the sample complexity of any estimator is exponential in $d$. 

\begin{theorem}[Exponential Sample Complexity]\label{TM:sample_complex_worstcase_asymp}  
The following claims hold:
\begin{enumerate}
    \item Fix $\sigma > 0$. There exist $d_0(\sigma)\in\mathbb{N}$, $\eta_0(\sigma)>0$ and $\gamma(\sigma)>0$ (monotonically decreasing in $\sigma$), such that for all $d\geq d_0(\sigma)$ and $\eta<\eta_0(\sigma)$, we have $n^\star(\eta,\sigma,\mathcal{F}_d) =\Omega\left(\frac{2^{\gamma(\sigma) d}}{\eta d}\right)$.
    
    \item Fix $d\in\mathbb{N}$. There exist $\sigma_0(d),\eta_0(d)>0$, such that for~all $\sigma<\sigma_0(d)$ and $\eta<\eta_0(d)$, we have $n^\star(\eta,\sigma,\mathcal{F}_d)=\Omega\left(\frac{2^d}{\eta d}\right)$.
\end{enumerate}
\end{theorem}

Theorem \ref{TM:sample_complex_worstcase_asymp} is proven in Section \ref{SUPPSEC:sample_complex_worstcase_proof_asymp}, based on channel coding arguments. For instance, the proof of Part 1 relates the estimation of $h(P\ast\Gauss)$ to discrete entropy estimation of a distribution supported on a capacity-achieving codebook for a peak-constrained AWGN channel. Since the codebook size is exponential in $d$, discrete entropy estimation over the codebook within a small gap $\eta>0$ is impossible with~less than order of $\frac{2^{\gamma(\sigma) d}}{\eta d}$ samples \cite{valiant2010clt,discrete_entropy_est_Wu2016}. The exponent $\gamma(\sigma)$ is monotonically decreasing in $\sigma$, implying that larger $\sigma$ values are favorable for estimation. The 2nd part of the theorem relies on a similar argument but for a $d$-dimensional AWGN channel and an input constellation that comprises the vertices of the $d$-dimensional hypercube $[-1,1]^d$.


\begin{remark}[Complexity for Restricted Distribution Classes]
Restricting $\mathcal{F}_d$ by imposing smoothness or lower-boundedness assumptions on the distributions in the class would not alleviate the exponential dependence on $d$ from Theorem \ref{TM:sample_complex_worstcase_asymp}. For instance, consider convolving any $P\in\mathcal{F}_d$ with $\Gausss$, i.e., replacing each $P$ with $Q=P\ast\Gausss$. These $Q$ distributions are smooth, but if one could accurately estimate $h\big(Q\ast\Gausss\big)$ over the convolved class, then $h(P\ast\Gauss)$ over $\mathcal{F}_d$ could have been estimated as well. Therefore, Theorem \ref{TM:sample_complex_worstcase_asymp} applies also for the class of such smooth $Q$ distributions. 
\end{remark}

The next propositions shows that the absolute-error risk attained by any estimator of $h(P\ast\Gauss)$ decays no faster than~$n^{-1/2}$.

\begin{proposition}[Risk Lower Bound]\label{PROP:minimaxRate}
For any $\sigma>0,\ d\geq1$, we have
\begin{equation}
    \mathcal{R}^\star\left(n,\sigma,\FSG\right)= \Omega\left(\frac{1}{\sqrt{n}}\right).
\end{equation}
\end{proposition}
That proposition states the so-called parametric lower bound on the absolute-error estimation risk. Under (the square root of the) quadratic loss, the result trivially follows from the Cram{\'e}r-Rao lower bound. For completeness, Appendix \ref{app:rateVerification} provides a simple proof for the absolute-error loss considered herein, based on the Hellinger modulus \cite{chen1997general}.

\subsection{Upper Bound on Risk}

We establish the minimax-rate optimality of the plug-in estimator by showing its risk converges as $n^{-1/2}$. Our risk bounds provide explicit constants (in terms of $\sigma$, $K$ and $d$). These constants depend exponentially on dimension, in accordance to the results of Theorem \ref{TM:sample_complex_worstcase_asymp}. Recall that given a collection of samples $S^n\sim P^{\otimes n}$, the estimator is $h(\hat{P}_{S^n}\ast\Gauss)$, where $\hat{P}_{S^n}=\frac{1}{n}\sum_{i=1}^n\delta_{S_i}$ is the empirical measure. 

A risk bound for the bounded support case is presented first. Although a special case of Theorem \ref{TM:SP_sample_complex_new}, where the subgaussian class $\FSG$ is considered, we state the bounded support result separately since it gives a cleaner bound with a better constant. 

\begin{theorem}[Plug-in Risk Bound - Bounded Support Class]\label{TM:SP_sample_complex_new_bdd}
Fix $\sigma>0$ and $d\geq1$. For any $n$, we have
\begin{equation}
    \sup_{P \in \mathcal{F}_d}\mathbb{E}_{P^{\otimes n}}\left|h(P\ast\Gauss)-h(\hat{P}_{S^n}\ast\Gauss)\right| \leq \tilde{C}_{\sigma,d} \frac{1}{\sqrt{n}},\label{EQ:Plugin_risk_bound_bdd}
\end{equation}
where $\tilde{C}_{\sigma,d}=O_{\sigma}(c^d)$, for a numerical constant $c$, is explicitly characterized in \eqref{EQ:Plugin_risk_bound_constant_bdd}.
\end{theorem}

The proof (given in Section \ref{SUBSEC:SP_sample_complex_new_bdd_proof}) relies on the $\Chi$-divergence $n^{-1/2}$ convergence rate from Corollary \ref{CORR:Chi_parametric_bdd}. Specifically, we relate the differential entropy estimation error to $\E_{P^{\otimes n}}\Chi\left(\hat{P}_{S^n}\ast\Gauss\middle\|P\ast\Gauss\right)$ using $\Chi$ variational representation. The result then follows by controlling certain variance terms and using Corollary~\ref{CORR:Chi_parametric_bdd}.


We next bound the estimation risk when $P\in\FSG$.

\begin{theorem}[Plug-in Risk Bound - Subgaussian Class]\label{TM:SP_sample_complex_new}
Fix $\sigma>0$ and $d\geq1$. For any $n$, we have
\begin{equation}
    \sup_{P \in \mathcal{F}_{d,K}^{(\mathsf{SG})}}\mathbb{E}_{P^{\otimes n}}\left|h(P\ast\Gauss)-h(\hat{P}_{S^n}\ast\Gauss)\right| \leq C_{\sigma,d,K} \frac{1}{\sqrt{n}},\label{EQ:Plugin_risk_bound}
\end{equation}
where $C_{\sigma,d,K}=O_{\sigma,K}(c^d)$, for a numerical constant $c$, is explicitly characterized in \eqref{EQ:Plugin_risk_bound_constant}.
\end{theorem}

The proof of Theorem \ref{TM:SP_sample_complex_new} is given in Section \ref{SUBSEC:plug-in_risk_bound_proof}. While the result follows via arguments similar to the bounded support case (namely, through the $\Chi$ subgaussian bound from Proposition \ref{PROP:Chi_parametric}), this method only covers the regime $\frac{K}{\sigma}<\frac{1}{2}$. To prove Theorem \ref{TM:SP_sample_complex_new} without restricting $\sigma$ and $K$, we resort to a different argument. Using the maximal TV-coupling, we bound the estimation risk by a weighted TV distance between $P\ast\Gauss$ and $\hat{P}_{S^n}\ast\Gauss$. The smoothing induced by the Gaussian convolutions allows us to control this TV distance by a $e^{O(d)} n^{-1/2}$. Several things to note about the result are the following:
\begin{enumerate}
    \item The theorem does not require any smoothness conditions on the distributions in $\mathcal{F}_{d,K}^{(\mathsf{SG})}$. This is achievable due to the inherent smoothing introduced by the convolution with the Gaussian density. Specifically, while the differential entropy $h(q)$ is not a smooth functional of the underlying density $q$ in general, our functional is $\mathsf{T}_\sigma(P)\triangleq h(P\ast\Gauss)$, which is smooth. 
    
    \item The above smoothness also allows us to avoid any assumptions on $P$ being bounded away from zero. So long as $P$ has subgaussian tails, the distribution may be arbitrary.
\end{enumerate}

\begin{remark}[Knowledge of Noise Parameter]\label{REM:unknown_sigma}
Our original motivation for this work is the noisy DNN setting, where additive Gaussian noise is injected into the system to enable tracking ``information flows'' during training (see \cite{ICML_Info_flow2019}). In this setting, the parameter $\sigma$ is known and the considered observation model reflects this. However, an interesting scenario is when $\sigma$ is unknown. To address this, first note that samples from $P$ contain no information about $\sigma$. Hence, in the setting where $\sigma$ is unknown, presumably samples of both $S\sim P$ and $S+Z\sim P\ast\mathcal{N}_\sigma$ would be available. Under this alternative model, estimating $\sigma$ can be done immediately by comparing the empirical variance of $S$ and $S+Z$. This empirical proxy would converge as $O\left((n d)^{-1/2}\right)$, implying that for large enough dimension, the empirical $\sigma$ can be substituted into our entropy estimator (in place of the true $\sigma$) without affecting the $O\left(c^d n^{-1/2}\right)$ convergence rate.
\end{remark}


\subsection{Bias Lower Bound}\label{SUBSEC:SP_estimator_implement}

To have a guideline as to the smallest number of samples needed to avoid biased estimation, we present the following lower bound on the estimator's bias $\sup_{P\in\mathcal{F}_d}\big|h(P\ast\Gauss)-\mathbb{E}_{P^{\otimes n}}h(\hat{P}_{S^n}\ast\Gauss)\big|$. 

\begin{theorem}[Bias Lower Bound]\label{TM:SP_Bias_LB}
Fix $d\geq 1$ and $\sigma>0$, and let $\epsilon\in\left(1-\left(1-2Q\left(\frac{1}{2\sigma}\right)\right)^d,1\right]$, where $Q$ is the Q-function.\footnote{The Q-function is defined as $Q(x)\triangleq\frac{1}{\sqrt{2\pi}}\int_x^\infty e^{-\frac{t^2}{2}}dt$.} Set $k_\star\triangleq\bigg\lfloor\frac{1}{\sigma Q^{-1}\left(\frac{1}{2}\left(1-(1-\epsilon)^{\frac{1}{d}}\right)\right)}\bigg\rfloor$, where $Q^{-1}$ is the inverse of the Q-function. By the choice of $\epsilon$, clearly $k_\star\geq 2$, and we have
\begin{align*}
 \sup_{P\in\mathcal{F}_d}\big|h(P\ast\Gauss)&-\mathbb{E}_{P^{\otimes n}}h(\hat{P}_{S^n}\ast\Gauss)\big|\\&\geq \log\left(\frac{k_\star^{d(1-\epsilon)}}{n}\right)-H_b(\epsilon). \numberthis\label{EQ:SP_Bias_LB}
\end{align*}
Consequently, the bias cannot be less than a given $\delta>0$ so long as $n\leq k_\star^{d(1-\epsilon)}\cdot e^{-(\delta+H_b(\epsilon))}$.
\end{theorem}
The theorem is proven in Section \ref{SUBSEC:SP_Bias_LB}. Since $H_b(\epsilon)$ shrinks with $\epsilon$, for sufficiently small $\epsilon$ values, the lower bound from \eqref{EQ:SP_Bias_LB} essentially shows that the our estimator will not have negligible bias unless $n> k_\star^{d(1-\epsilon)}$ is satisfied. The condition $\epsilon>1-\left(1-2Q\left(\frac{1}{2\sigma}\right)\right)^d$ is non-restrictive in any relevant regime of $d$ and $\sigma$. For the latter, values we have in mind are inspired by \cite{ICML_Info_flow2019}, where noisy DNNs with parameter $\sigma$ are studied. In that work,  $\sigma$ values are around $0.1$, for which the lower bound on $\epsilon$ is at most 0.0057 for all dimensions up to at least $d=10^4$. For example, when setting $\epsilon=0.01$ (for which $H_b(0.01)\approx 0.056$), the corresponding $k_\star$ equals 3 for $d\leq 11$ and 2 for $12\leq d\leq 10^4$. Thus, with these parameters, a negligible bias requires $n$ to be at least $2^{0.99d}$.

\subsection{Computing the Estimator}\label{SUBSEC:SP_estimator_compute}

Evaluating the plug-in estimator $h(\hat{P}_{S^n}\ast\Gauss)$ requires computing the differential entropy of a $d$-dimensional $n$-mode Gaussian mixture ($\hat{P}_{S^n}\ast\Gauss$). Although it cannot be computed in closed form, this section presents a method for computing an arbitrarily accurate approximation via MC integration \cite{robert2004montecarlo}. To simplify the presentation, we present the method for an arbitrary Gaussian mixture without referring to the notation of the estimation setup. 

Let $g(t)\triangleq\frac{1}{n}\sum_{i=1}^n \gauss(t-\mu_i)$ be a $d$-dimensional, $n$-mode Gaussian mixture, with centers $\{\mu_i\}_{i=1}^n\subset\mathbb{R}^d$. Let $C\sim\mathsf{Unif}\p{\{\mu_i\}_{i=1}^n}$ be independent of $Z\sim\Gauss$ and note that $V\triangleq C+Z\sim g$. First, rewrite $h(g)$ as follows:
\begin{align*}
    h(g)=-\mathbb{E}\log g(V)&=-\frac{1}{n}\sum_{i=1}^n\mathbb{E}\Big[\log g(\mu_i+Z)\Big|C=\mu_i\Big]\\&=-\frac{1}{n}\sum_{i=1}^n\mathbb{E}\log g(\mu_i+Z),\numberthis\label{EQ:MC_expansion}
\end{align*}
where the last step uses the independence of $Z$ and $C$. Let $\cp{Z_j^{(i)}}_{\substack{i\in[n]\\j\in[\nmc]}}$ be $n\times \nmc$ i.i.d. samples from $\gauss$. For each $i\in[n]$, we estimate the $i$-th summand on the RHS of \eqref{EQ:MC_expansion} by
\begin{subequations}
\begin{equation}
    \hat{I}_\mathsf{MC}^{(i)}\triangleq \frac{1}{\nmc}\sum_{j=1}^{\nmc}\log g\left(\mu_i+Z_j^{(i)}\right),\label{EQ:MC_per_i_est}
\end{equation}
which produces 
\begin{equation}
    \hat{h}_\mathsf{MC}\triangleq \frac{1}{n}\sum_{i=1}^n\hat{I}_\mathsf{MC}^{(i)}\label{EQ:MC_full_est}
\end{equation}\label{EQ:MCI}%
\end{subequations}
as the approximation of $h(g)$. Note that since $g$ is a mixture of $n$ Gaussians, it can be efficiently evaluated using off-the-shelf KDE software packages, many of which require only $O(\log n)$ operations on average per evaluation of $g$.

Define the MSE of $\hat{h}_\mathsf{MC}$ as
\begin{equation}
    \mathsf{MSE}\left(\hat{h}_\mathsf{MC}\right)\triangleq \mathbb{E}\bigg[\Big(\hat{h}_\mathsf{MC}-h(g)\Big)^2\bigg].
\end{equation}

We have the following bounds on the MSE. 
\begin{theorem}[MSE Bounds for the MC Estimator]\label{TM:MC_MSE}\ \\

\begin{enumerate}[(i)]
    \item \underline{Bounded support:} Assume $C\in[-1,1]^d$ almost surely, then 
    \begin{equation}
        \mathsf{MSE}\left(\hat{h}_\mathsf{MC}\right)\leq \frac{2d(2+\sigma^2)}{\sigma^2}\frac{1}{n\cdot \nmc}.\label{EQ:MC_MSE_Tanh}
    \end{equation}
    \item \underline{Bounded moment:} Assume $m\triangleq\mathbb{E}\|C\|_2^2<\infty$, then
    \begin{align*}
        &\mathsf{MSE}\left(\hat{h}_\mathsf{MC}\right)\leq\numberthis\label{EQ:MC_MSE_ReLU}\\& \frac{9d\sigma^2+8(2+\sigma\sqrt{d})m+3(11\sigma\sqrt{d}+1)\sqrt{m}}{\sigma^2}\frac{1}{n\cdot \nmc}.
    \end{align*}
\end{enumerate}
\end{theorem}
The proof is given in Section \ref{SUBSEC:MC_MSE_proof}. The bounds on the MSE scale only linearly with the dimension $d$, making $\sigma^2$ in the denominator often the dominating factor experimentally. 

\section{Information Flow in Deep Neural Networks}\label{SEC:applications}

A utilization of the developed theory is estimating the mutual information between selected groups of neurons in DNNs. Much attention was recently devoted to this task \cite{DNNs_ICLR2018,jacobsen2018revnet,liu2018understanding,reeves2017additivity,gabrie2018entropy,ICML_Info_flow2019}, partly motivated by the Information Bottleneck (IB) theory for DNNs \cite{tishby_DNN2015,DNNs_Tishby2017}. The theory tracks the mutual information pair $\big(I(X;T),I(Y;T)\big)$, where $X$ is the DNN's input (i.e., feature), $Y$ is the true label and $T$ is the hidden representation vector. An interesting claim from \cite{DNNs_Tishby2017} is that the mutual information $I(X;T)$ undergoes a so-called `compression' phase during training. Namely, after an initial short `fitting' phase (where $I(Y;T)$ and $I(X;T)$ both grow), $I(X;T)$ exhibits a slow long-term decrease, which is termed the `compression' phase. According to \cite{DNNs_Tishby2017}, this phase explains the excellent generalization performance of DNNs. 

The main caveat in the supporting empirical results from \cite{DNNs_Tishby2017} (and the partially opposing results from the followup work \cite{DNNs_ICLR2018}) is that in deterministic networks, where $T=f(X)$ with strictly monotone activations, $I(X;T)$ is either infinite (when the data distribution $P_X$ is continuous) or a constant (when $P_X$ is discrete\footnote{The mapping from the discrete values of $X$ to $T$ is almost always (except for a measure-zero set of weights) injective whenever the nonlinearities are, thereby causing $I(X;T)=H(X)$ for any hidden layer $T$, even if $T$ consists of a single neuron.}). As explained in \cite{ICML_Info_flow2019}, the reason \cite{DNNs_Tishby2017} and \cite{DNNs_ICLR2018} miss this fact stems from an inadequate application of a binning-based mutual information estimator for $I(X;T)$. 


To fix this constant/infinite mutual information issue, \cite{ICML_Info_flow2019} proposed the framework of noisy DNNs, in which each neuron adds a small amount of Gaussian noise (i.i.d. across all neurons) after applying the activation function. The injected noise makes the map $X \mapsto T$ a stochastic parameterized channel, and as a consequence, $I(X;T)$ is a finite quantity that depends on the network's parameters. Although the primary purpose of the noise injection in \cite{ICML_Info_flow2019} was to ensure that $I(X;T)$ depends on the system parameters, experimentally it was found that the network's performance is optimized at non-zero noise variance, thus providing a natural to select this parameter. In the following, we first define noisy DNNs and then show that estimating $I(X;T)$, $I(Y;T)$ or any other mutual information term between layers of a noisy DNN can be reduced to differential entropy estimation under Gaussian convolutions. The reduction relies on a sampling procedure that leverages the DNN's generative model.




\subsection{Noisy DNNs and Mutual Information between Layers}\label{SUBSEC:MI_Noisy_DNNs}

We start by describing the noisy DNN setup from \cite{ICML_Info_flow2019}. Let $(X,Y)\sim P_{X,Y}$ be a feature-label pair, where $P_{X,Y}$ is the (unknown) true distribution of $(X,Y)$, and $\big\{(X_i,Y_i)\big\}_{i=1}^n$ be $n$ i.i.d. samples from~$P_{X,Y}$. 

\begin{figure}[t!]
	\begin{center}
		\begin{psfrags}
			\psfragscanon
			\psfrag{A}[][][1]{\ \ \ \ \ \ $T_{\ell-1}$}
			\psfrag{B}[][][0.8]{$a\mspace{-3mu}\left(\mathrm{W}_\ell^{(k)}T_{\ell-1}\mspace{-3mu}+\mspace{-3mu}b_\ell(k)\right)$}
			\psfrag{C}[][][1]{\ \ \ \ \ \ \ $S_\ell(k)$}
			\psfrag{D}[][][1]{\ \ \ \ \ \ \ \ \ \ \ \ \ \ \ \ \ \ \ \ $Z_\ell(k)\sim\mathcal{N}(0,\sigma^2)$}
			\psfrag{E}[][][1]{$\mspace{-5mu}T_\ell(k)$}
			\includegraphics[scale = .71]{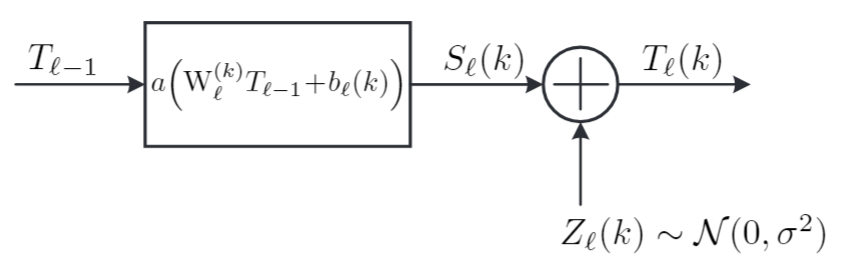}
			\caption{$k$-th noisy neuron in a fully connected or a convolutional layer $\ell$ with activation function $a$; $\mathrm{W}_\ell^{(k)}$ and $b_\ell(k)$ are the $k$-th row and the $k$-th entry of the weight matrix and the bias vector,  respectively.}\label{FIG:noisy_neuron}
			\psfragscanoff
		\end{psfrags}
	\end{center}
\end{figure}

Consider an $(L+1)$-layered (fixed / trained) noisy DNN with layers $T_0,T_1,\ldots,T_L$, input $T_0=X$ and output $T_L=\hat{Y}$ (i.e., an estimate of $Y$). For each $\ell\in[L-1]$, the $\ell$-th hidden layer is given by $T_\ell=S_\ell+Z_\ell$, where $S_\ell\triangleq f_\ell(T_{\ell-1})$ with $f_\ell:\mathbb{R}^{d_{\ell-1}}\to\mathbb{R}^{d_\ell}$ being a deterministic function of the previous layer and $Z_\ell\sim\mathcal{N}\big(0,\sigma^2 \mathrm{I}_{d_\ell}\big)$ being the noise injected at layer $\ell$. The functions $f_1,f_2,\ldots,f_L$ can represent any type of layer (fully connected, convolutional, max-pooling, etc.). For instance, $f_\ell(t)=a(\mathrm{W}_\ell t+b_\ell)$ for a fully connected or a convolutional layer, where $a$ is the activation function which operates on a vector component-wise, $\mathrm{W}_\ell\in\mathbb{R}^{d_\ell\times d_{\ell-1}}$ is the weight matrix and $b_\ell\in\mathbb{R}^{d_\ell}$ is the bias. For fully connected layers $\mathrm{W}_\ell$ is arbitrary, while for convolutional layers $\mathrm{W}_\ell$ is Toeplitz. Fig.~\ref{FIG:noisy_neuron} shows a neuron in a noisy DNN.



The noisy DNN induces a stochastic map from $X$ to the rest of the network, described by the conditional distribution $P_{T_1,\ldots,T_L|X}$. The joint distribution of the tuple $(X,Y,T_1,\ldots,T_L)$ is $P_{X,Y,T_1,\ldots,T_L}\triangleq P_{X,Y}P_{T_1,\ldots,T_L|X}$ under which $Y-X-T_1-\ldots-T_L$ forms a Markov chain. 
For any $\ell\in[L-1]$, consider the mutual information between the hidden layer and the input (see Remark \ref{REM:MI_label} for an account of $I(Y;T_\ell)$):
\begin{align}
    I(X;T_\ell)&=h(T_\ell)-h(T_\ell|X)\nonumber\\&=h(P_{T_\ell})-\int dP_X(x) h(P_{T_\ell|X=x}).\label{EQ:DNN_mutual_information_input}
\end{align}
Since $P_{T_\ell}$ and $P_{T_\ell|X}$ have a highly complicated structure (due to the composition of Gaussian noises and nonlinearities), this mutual information cannot be computed analytically and must be estimated. Based on the expansion from \eqref{EQ:DNN_mutual_information_input}, an estimator of $I(X;T_\ell)$ is constructed by estimating the unconditional and each of the conditional differential entropy terms, while approximating the expectation by an empirical average. As explained next, all these entropy estimation tasks are instances of our framework of estimating $h(P\ast\Gauss)$ based on samples from $P$ and knowledge of $\sigma$.


\subsection{From Differential Entropy to Mutual Information}\label{SUBSEC:MI_Noisy_DNNs_estimation}

Recall that $T_\ell=S_\ell+Z_\ell$, where $S_\ell\sim P_{S_\ell}=P_{f_\ell(T_{\ell-1})}$ and $Z_\ell\sim\mathcal{N}(0,\sigma^2\mathrm{I}_{d_\ell})$ are independent. Thus, 
\begin{subequations}
\begin{equation}
h(P_{T_\ell})=h(P_{S_\ell}\ast\Gauss)\label{EQ:uncond_ent}
\end{equation}
and
\begin{equation}
h(P_{T_\ell|X=x_i})=h(P_{S_\ell|X=x_i}\ast\Gauss).\label{EQ:cond_ent}
\end{equation}
\end{subequations}
The DNN's forward pass enables sampling from $P_{S_\ell}$ and $P_{S_\ell|X}$ as follows: 
\begin{enumerate}
    \item \underline{Unconditional Sampling:} To generate the sample set from $P_{S_\ell}$, feed each $X_i$, for $i\in[n]$, into the DNN and collect the outputs it produces at the $(\ell-1)$-th layer. The function $f_\ell$ is then applied to each collected output to obtain $S_\ell^n\triangleq\cp{S_{\ell,1},S_{\ell,2},\ldots,S_{\ell,n}}$, which is the a set of $n$ i.i.d. samples from $P_{S_\ell}$. 

    \item \underline{Conditional Sampling Given $X$:} 
    To generate i.i.d. samples from $P_{S_\ell|X=x_i}$, for $i\in[n]$, we feed $x_i$ into the DNN $n$ times, collect outputs from $T_{\ell-1}$ corresponding to different noise realizations, and apply $f_\ell$ on each. Denote the obtained samples by $S_\ell^n(X_i)$.
    \footnote{The described sampling procedure is valid for any layer $\ell\geq 2$. For $\ell=1$, $S_1$ coincides with $f_1(X)$ but the conditional samples are undefined. Nonetheless, noting that for the first layer $h(T_1|X)=h(Z)=\frac{d}{2}\log(2\pi e \sigma^2)$, we see that no estimation of the conditional entropy is needed. The mutual information estimator given in  \eqref{EQ:MI_CONT} is modified by replacing the subtracted term with $h(Z)$.}


\end{enumerate}
The knowledge of $\sigma$ and together with the samples $S_\ell^n$ and $S_\ell^n(X_i)$ can be used to estimate the unconditional and the conditional entropies, from \eqref{EQ:uncond_ent} and \eqref{EQ:cond_ent}, respectively.

For notational simplicity, we henceforth omit the layer index $\ell$. Based on the above sampling procedure we construct an estimator $\hat{I}\big(X^n,\hat{h}\big)$ of $I(X;T)$ using a given estimator $\hat{h}(A^n,\sigma)$ of $h(P\ast\Gauss)$ for $P$ supported inside $[-1,1]^d$ (i.e., a tanh / sigmoid network), based on i.i.d. samples $A^n=\{A_1,\ldots,A_n\}$ from $P$ and knowledge of $\sigma$. Assume that $\hat{h}$~attains
\begin{equation}
    \sup_{P\in\mathcal{F}_d} \mathbb{E}_{P^{\otimes n}}\left|h(P\ast\Gauss)-\hat{h}(A^n,\sigma)\right|\leq \Delta_{\sigma,d}(n).\label{EQ:Data_Dist_Assumption}
\end{equation}
An example of such an $\hat{h}$ is the estimator $h(\hat{P}_{A^n}\ast\Gauss)$. The corresponding $\Delta_{\sigma,d}(n)$ term is given in Theorem~\ref{TM:SP_sample_complex_new}. Our estimator for the mutual information is
\begin{equation}
\hat{I}_{\mathsf{Input}}\p{X^n,\hat{h},\sigma}\triangleq \hat{h}(S^n,\sigma) - \frac{1}{n} \sum_{i=1}^n \hat{h}\big(S^n(X_i),\sigma\big).\label{EQ:MI_CONT}
\end{equation}
The absolute-error estimation risk of $\hat{I}_{\mathsf{Input}}\p{X^n,\hat{h},\sigma}$ is bounded in the following proposition, proven in Section~\ref{SUBSEC:MI_True_Data_Dist_proof}.


\begin{proposition}[Input--Hidden Layer Mutual Information]\label{PROP:MI_True_Data_Dist}
For the above described estimation setting, we have
\begin{align*}
    \sup_{P_X}\mathbb{E} &\left|I(X;T)-\hat{I}_{\mathsf{Input}}\p{X^n,\hat{h},\sigma}\right|\\&\qquad\qquad\leq 2\Delta_{\sigma,d}(n)+\frac{d\log\left(1+\frac{1}{\sigma^2}\right)}{4\sqrt{n}}. 
\end{align*}
\end{proposition}
The quantity $\frac{1}{\sigma^2}$ is the SNR between $S$ and $Z$. The larger $\sigma$ is the easier estimation becomes, since the noise smooths out the complicated $P_X$ distribution. Also note that the dimension of the ambient space in which $X$ lies does not appear in the absolute-risk bound. The bound depends only on the dimension of $T$ (through $\Delta_{\sigma,d}$). This happens because the blurring effect caused by the noise enables uniformly lower bounding $\inf_x h(T|X=x)$ and thereby controlling the variance of the estimator for each conditional entropy. This reduces the impact of $X$ on the estimation error to that of an empirical average converging to its expected value with rate $n^{-1/2}$.


\begin{remark}[Subgaussian Class and Noisy ReLU DNNs]\label{REM:SG_natural}
We provide performance guarantees for the plug-in estimator also over the more general class  $\mathcal{F}_{d,K}^{(\mathsf{SG})}$ of distributions with subgaussian marginals. This class accounts for the following important cases:
\begin{enumerate}
    \item Distributions with bounded support, which correspond to noisy DNNs with bounded nonlinearities. This case is directly studied through the bounded support class $\mathcal{F}_d$. 
    \item Discrete distributions over a finite set, which is a special case of bounded support.
    \item Distributions $P$ of a random variable $S$ that is a hidden layer of a noisy ReLU DNN, so long as the input $X$ to the network is itself subgaussian. To see this recall that linear combinations of independent subgaussian random variables are also subgaussian. Furthermore, for any (scalar) random variable $A$, we have that $\big|\mathsf{ReLU}(A)\big|=\big|\max\{0,A\}\big|\leq |A|$, almost surely. Each layer in a noisy $\mathsf{ReLU}$ DNN is a coordinate-wise $\mathsf{ReLU}$ applied to a linear transformation of the previous layer plus a Gaussian noise. Consequently, for a $d$-dimensional hidden layer $S$ and any $i\in[d]$, one may upper bound $\big\|S(i)\big\|_{\psi_2}$ by a constant, provided that the input $X$ is coordinate-wise subgaussian. This constant depends on the network's weights and biases, the depth of the hidden layer, the subgaussian norm of the input, and the noise variance.
\end{enumerate}
In the context of estimation of mutual information over DNNs, the input distribution is typically taken as uniform over the dataset \cite{DNNs_Tishby2017,DNNs_ICLR2018,ICML_Info_flow2019}, adhering to case (2).
\end{remark}

\begin{figure*}[t!]
    \centering\includegraphics[width=0.83\linewidth]{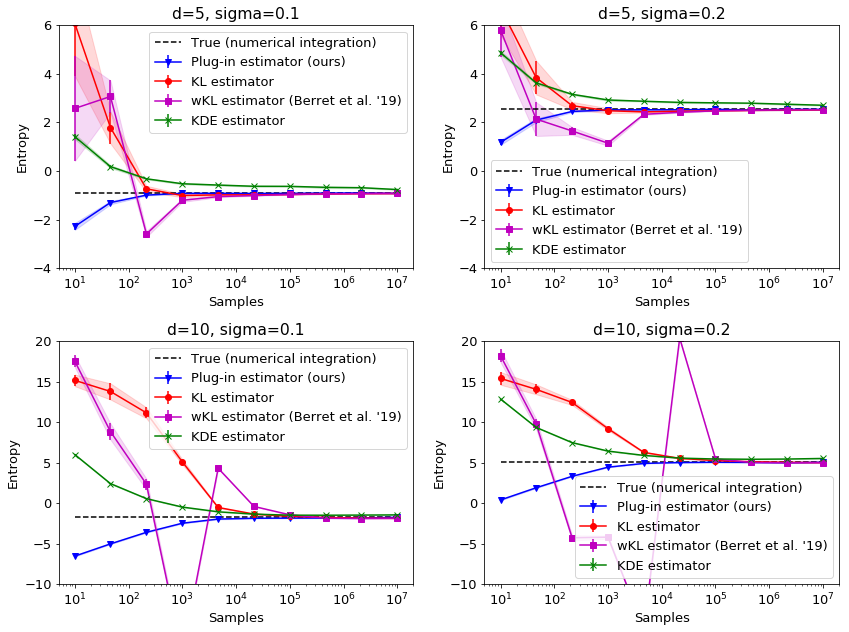}
    \caption{Estimation results comparing the plug-in estimator to: (i) a KDE-based method \cite{kandasamy2015nonparametric}; (ii) the KL estimator \cite{kozachenko1987sample}; and (iii) the wKL estimator \cite{berrett2019efficient}. The differential entropy of $S+Z$ is estimated, where $S$ is a truncated $d$-dimensional mixture of $2^d$ Gaussians and $Z\sim\mathcal{N}(0,\sigma^2\mathrm{I}_d)$. Results are shown as a function of $n$, for $d = 5,10$ and $\sigma=0.1,0.5$. Error bars are one standard deviation over 20 random trials. The $h(P\ast\Gauss)$ estimator presents faster convergence rates, improved stability and better scalability with dimension compared to the two competing methods.} \label{Fig:CornersGaussian}
\end{figure*}

\begin{remark}[Hidden Layer--Label Mutual Information]\label{REM:MI_label} Another quantity of interest is the mutual information between the hidden layer and the true label (see, e.g., \cite{DNNs_Tishby2017}). For $(X,Y)\sim P_{X,Y}$, and a hidden layer $T$ in a noisy DNN with input $X$, the joint distribution of $(X,Y,S,T)$ is $P_{X,Y}P_{S,T|X}$, under which $Y-X-(S,T)$ forms a Markov chain.\footnote{In fact, the Markov chain is $Y-X-S-T$ since $T=S+Z$, but this is inconsequential here.} The mutual information of interest is then 
\begin{equation}
    I(Y;T)=h(P_S\ast\Gauss)-\sum_{y\in\mathcal{Y}}P_Y(y)h(P_{S|Y=y}\ast\Gauss),\label{EQ:MI_label_rem}
\end{equation}
where $\mathcal{Y}$ is the (known and) finite set of labels. Just like for $I(X;T)$, estimating $I(Y;T)$ reduces to differential entropy estimation under Gaussian convolutions. Namely, an estimator for $I(Y;T)$ can be constructed by estimating the unconditional and each of the conditional differential entropy terms in \eqref{EQ:MI_label_rem}, while approximating the expectation by an empirical average. There are several required modifications for estimating $I(Y;T)$ as compared to $I(X;T)$. Most notably is the procedure for sampling from $P_{S|Y=y}$, which results in a sample set whose size is random (Binomial). In appendix \ref{APPEN:label_MI}, the estimation of $I(Y;T)$ is described in detail and a corresponding risk bound is derived. 
\end{remark}

This section shows that the performance in estimating mutual information depends on our ability to estimate $h(P\ast\Gauss)$. In Section \ref{SEC:simulations} we present experimental results for $h(P\ast\Gauss)$, when $P$ is induced by a DNN.

\section{Simulations}\label{SEC:simulations}

We present empirical results illustrating the convergence of the plug-in estimator compared to several competing methods: (i) the KDE-based estimator of \cite{kandasamy2015nonparametric}; (ii) and kNN Kozachenko-Leonenko (KL) estimator \cite{kozachenko1987sample}; and (iii) the recently developed wKL estimator from \cite{berrett2019efficient}. These competing methods are general-purpose estimators of the differential entropy $h(Q)$ based on i.i.d. samples from $Q$. Such methods are applicable for estimating $h(P\ast\Gauss)$ by sampling $\Gauss$ and adding the noise values to the samples from $P$.

\subsection{Simulations for Differential Entropy Estimation}

\subsubsection{\underline{$P$ with Bounded Support}} Convergence rates in the bounded support regime are illustrated first. We set $P$ as a mixture of Gaussians truncated to have support in $[-1,1]^d$. Before truncation, the mixture consists of $2^d$ Gaussian components with means at the $2^d$ corners of $[-1,1]^d$ and standard deviations 0.02. This produces a distribution that is, on one hand, complicated ($2^d$ mixtures) while, on the other hand, is still simple to implement. The entropy $h(P\ast \Gauss)$ is estimated for various values of $\sigma$.

\begin{figure*}[t!]
    \centering\includegraphics[width=6in]{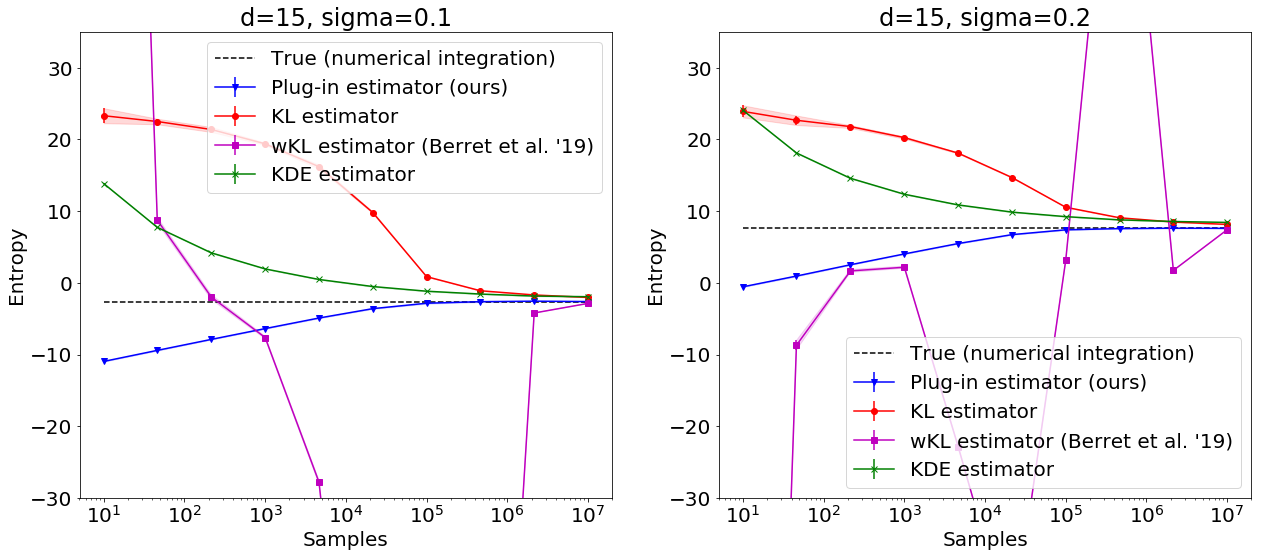}
    \caption{Estimation results comparing the plug-in estimator to: (i) a KDE-based method \cite{kandasamy2015nonparametric}; (ii) the KL estimator \cite{kozachenko1987sample}; and (iii) the wKL estimator \cite{berrett2019efficient}. Here $P$ is an untruncated $d$-dimensional mixture of $2^d$ Gaussians and $Z\sim\mathcal{N}(0,\sigma^2\mathrm{I}_d)$. Results are shown as a function of $n$, for $d = 5,10$ and $\sigma=0.01,0.1,0.5$. Error bars are one standard deviation over 20 random trials.  
    }
    \label{Fig:CornersGaussianNoTrunc}
\end{figure*}

Fig. \ref{Fig:CornersGaussian} shows estimation results as a function of $n$, for $d=5,10$ and $\sigma=0.1,0.2$. The KL and plug-in estimators require no tuning parameters; for wKL we used the default weight setting in the publicly available software. We stress that the KDE estimate is highly unstable and, while not shown here, the estimated value is very sensitive to the chosen kernel width. The kernel width (varying with both $d$ and $n$) for the KDE estimate was chosen by generating a variety of different Gaussian mixture constellations of moderately different
cardinalities and optimizing the kernel width for good performance across regimes (evaluated by comparing finite sample  estimates to the large-sample entropy estimate).\footnote{This gives an unrealistic advantage to KDE method, but we prefer this to uncertainty about baseline performance that stems from inferior kernel width selection.} As seen in Fig. \ref{Fig:CornersGaussian}, the KDE, KL and wKL estimators converge slowly, at a rate that degrades with increased $d$, underperforming the plug-in estimator. Finally, we note that in accordance to the explicit risk bound from \eqref{EQ:Plugin_risk_bound_constant}, the absolute error increases with larger $d$ and smaller $\sigma$.


\vspace{2mm}

\subsubsection{\underline{$P$ with Unbounded Support}} In Fig.~\ref{Fig:CornersGaussianNoTrunc}, we show the convergence rates in the unbounded support regime by considering the same setting with $d = 15$ but without truncating the $2^d$-mode Gaussian mixture. The fast convergence of the plug-in estimator is preserved, outperforming the competing methods. Notice that the performance of the wKL estimator from \cite{berrett2019efficient} (whose asymptotic efficiency was established therein) deteriorates in this relatively high-dimensional setup. This may be a result of the dependence of its estimation error on $d$, which was not characterized in \cite{berrett2019efficient}.


\subsection{Monte Carlo Integration}

\begin{figure}[!t]
    \centering
    \includegraphics[width=\columnwidth]{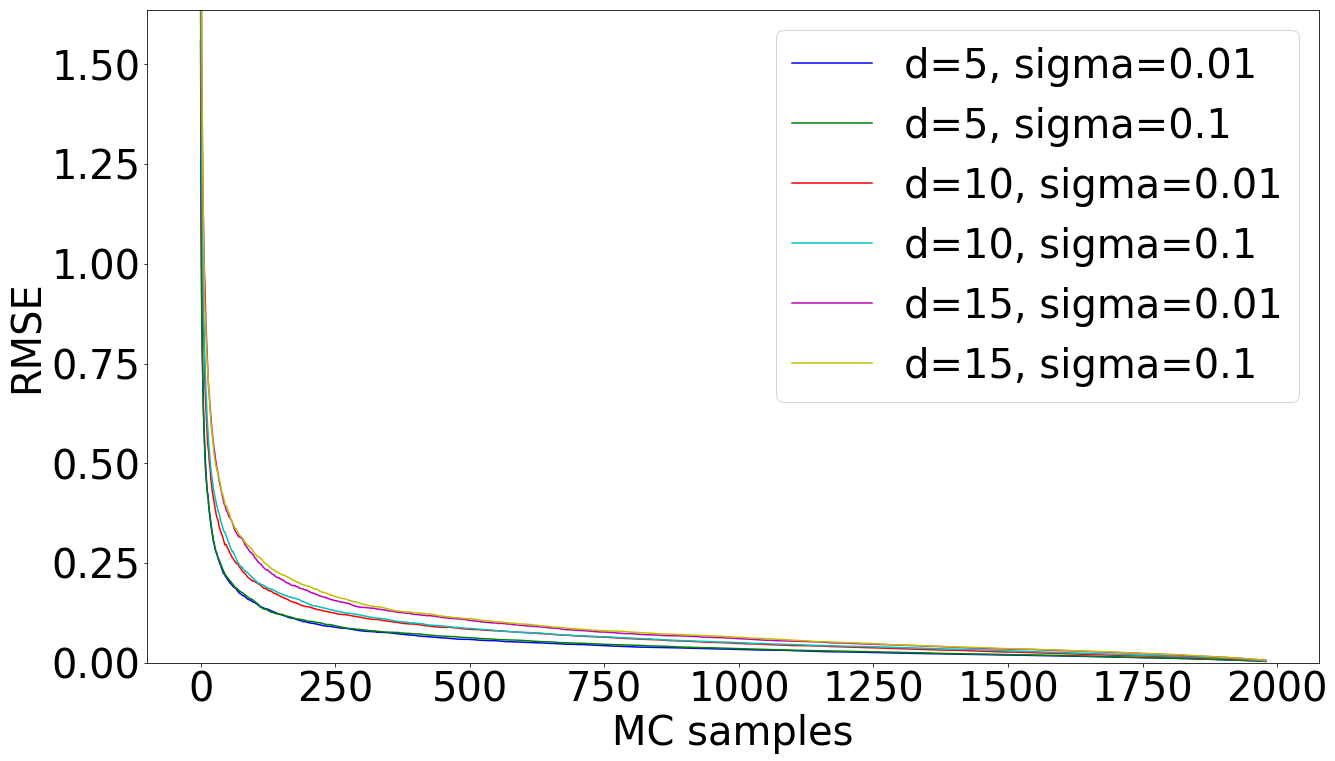}
    \caption{Convergence of the Monte Carlo integrator computation of the proposed estimator. Shown is the decay of the RMSE as the number of Monte Carlo samples increases, for a variety of $\sigma$ and $d$ values. The MC integrator is computing the $h(P\ast\Gauss)$ estimate of the entropy of $S+Z$ where $S$ is a truncated $d$-dimensional mixture of $2^d$ Gaussians and $Z\sim\mathcal{N}(0,\sigma^2\mathrm{I}_d)$. The number of samples of $S$ used by $h(P\ast\Gauss)$ is $10^4$. }\label{Fig:MC}
\end{figure}

Fig.~\ref{Fig:MC} illustrates the convergence of the MC integration method for computing the plug-in estimator. The figure shows the root-MSE (RMSE) as a function of MC samples $\nmc$, for the truncated $2^d$ Gaussian mixture distribution with $n=10^4$ (which corresponds to the number of modes in the Gaussian mixture $\hat{P}_{S^n}\ast\Gauss$ whose entropy approximates $h(P\ast\Gauss)$), $d = 5,10,15$, and $\sigma = 0.01,0.1$. Note the error decays approximately as $\nmc^{1/2}$ in accordance with Theorem \ref{TM:MC_MSE}, and that the convergence does not vary excessively for different $d$ and $\sigma$~values.

\begin{figure*}[t!]
    \centering
    \subfloat[]{\includegraphics[width=2.7in]{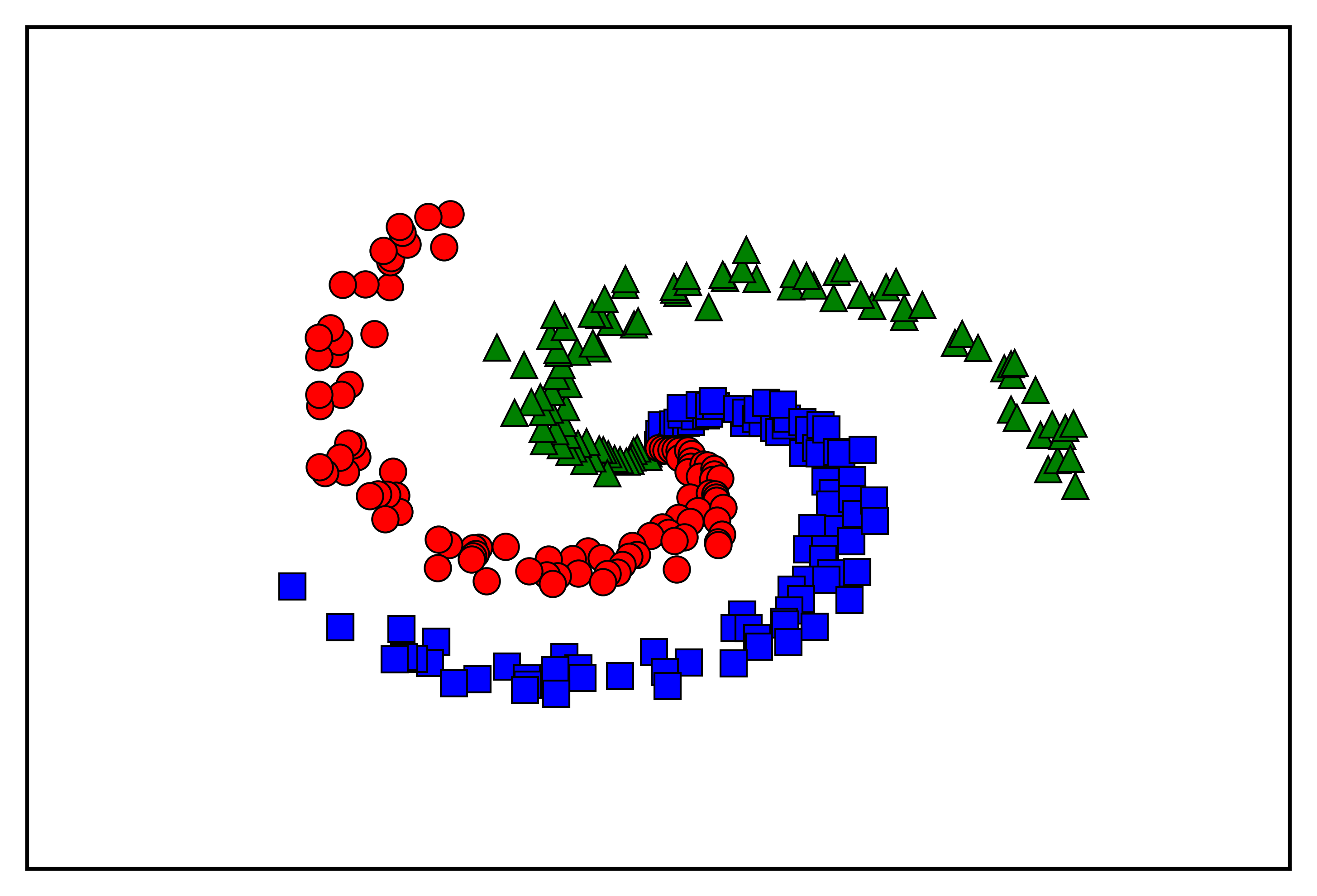}}\ \ \ \ \subfloat[]{\includegraphics[width=3.1in]{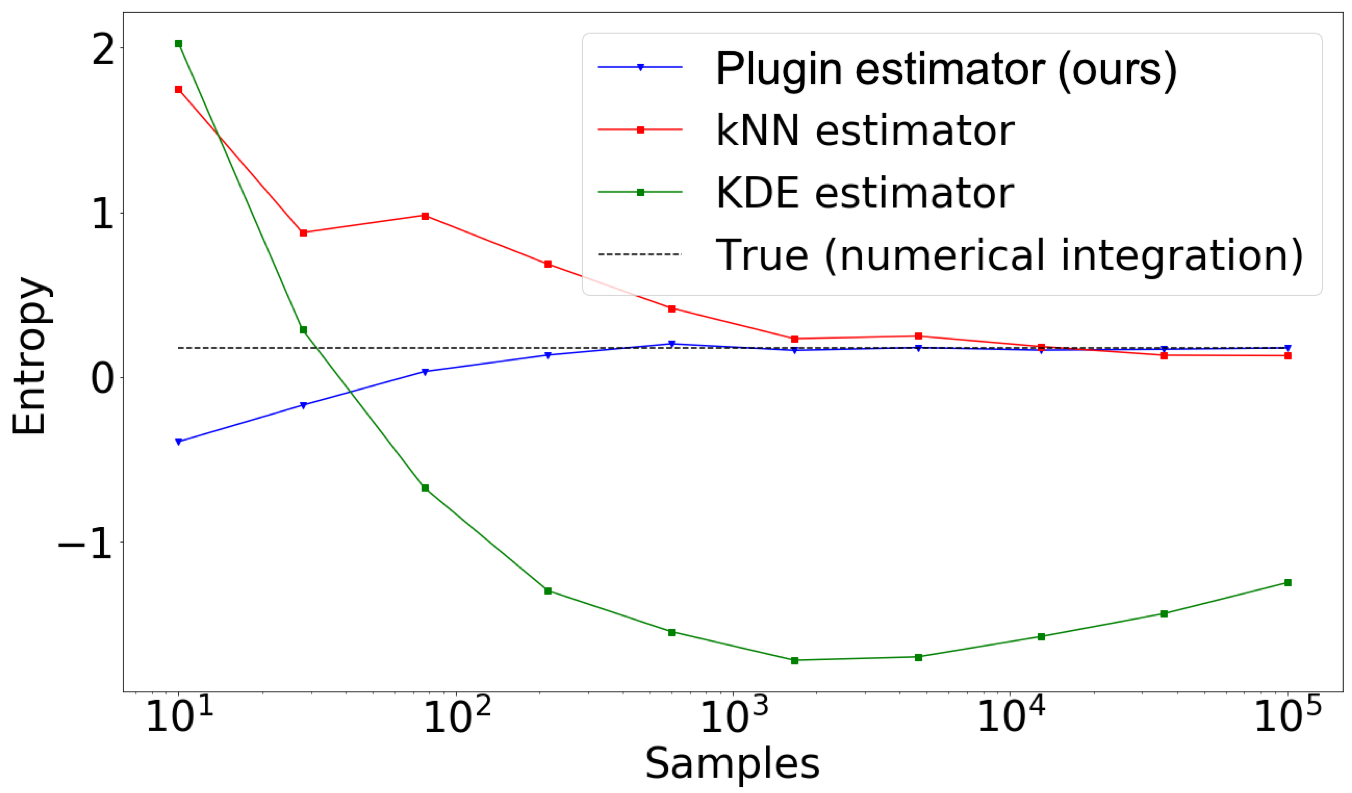}}
    \caption{10-dimensional entropy estimation in a 3-layer neural network trained on the 2-dimensional 3-class spiral dataset shown on the left. Estimation results for the plug-in estimator compared to general-purpose kNN and KDE methods are shown on the right. The differential entropy of $S+Z$ is estimated, where $S$ is the output of the third (10-dimensional) layer. Results are shown as a function of samples $n$ with $\sigma =0.2$.
    }\label{Fig:SpiralNN}
\end{figure*}

    
\subsection{Estimation in a Noisy Deep Neural Network}

We next illustrate entropy estimation in a noisy DNN. The dataset is a $2$-dimensional 3-class spiral (shown in Fig.~\ref{Fig:SpiralNN}(a)). The network has 3 fully connected layers of sizes 8-9-10, with tanh activations and $\mathcal{N}(0,\sigma^2)$ Gaussian noise added to the output of each neuron, where $\sigma = 0.2$. We estimate the entropy of the output of the 10-dimensional third layer in the network trained to achieve 98\% classification accuracy. Estimation results are shown in Fig.~\ref{Fig:SpiralNN}(b), comparing the plug-in estimator to the KDE and KL estimators; the wKL estimator from \cite{berrett2019efficient} is omitted due to its poor performance in this experiment. As before, the plug-in estimate converges faster than the competing methods illustrating its efficiency for entropy and mutual information estimation over noisy DNNs. The KDE estimate, which performed quite well in the synthetic experiments, underperform here. In our companion work \cite{ICML_Info_flow2019}, additional examples of mutual information estimation in DNNs based on the proposed estimator are provided.



\begin{figure*}[t!]
    \centering
    \subfloat[]{\includegraphics[width=3.1in]{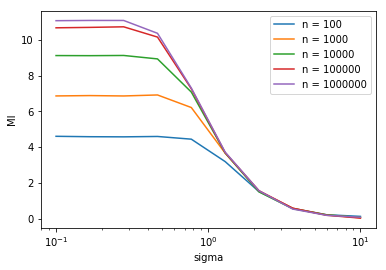}}\ \ \ \ 
    \subfloat[]{\includegraphics[width=3.1in]{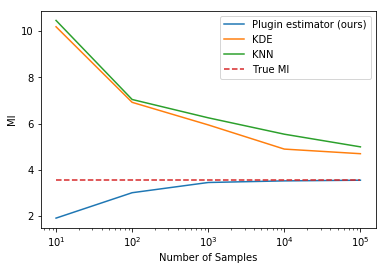}}
    \caption{Estimating $I(S;S+Z)$, where $S$ comes from a BPSK modulated Reed-Muller and $Z\sim\mathcal{N}(0,\sigma^2\mathrm{I}_d)$: (a) Estimated $I(S;S+Z)$ as a function of $\sigma$, for different $n$ values, for the $\mathsf{RM}(4,4)$ code. (b) Plug-in, KDE and KL $I(S;S+Z)$ estimates for the $\mathsf{RM}(5,5)$ code and $\sigma=2$ as a function~of~$n$. Shown for comparison are the curves for the kNN and KDE estimators based on noisy samples of $S+Z$ as well as the true value (dashed).}\label{Fig:ReedMuller}
\end{figure*}

\subsection{Reed-Muller Codes over AWGN Channels}
We next consider data transmission over an AWGN channel using a binary phase-shift keying (BPSK) modulation of a Reed-Muller code. A Reed-Muller code $\mathsf{RM}(r,m)$ of parameters $r,m\in\mathbb{N}$, where $0\leq r \leq m$, encodes messages of length $k=\sum _{i=0}^{r}{\binom {m}{i}}$ into $2^m$-lengthed binary codewords. Let $\mathcal{C}_{\mathsf{RM}(r,m)}$ be set of BPSK modulated sequences corresponding to $\mathsf{RM}(r,m)$ (with $0$ and $1$ mapped to $-1$ and $1$, respectively). The number of bits reliably transmittable over the $2^m$-dimensional AWGN with noise $Z\sim\mathcal{N}(0,\sigma^2\mathrm{I}_{2^m})$ is given by
\begin{equation}
     I(S;S+Z)=h(S+Z)-2^{m-1}\log(2\pi e \sigma^2),
\end{equation}
where $S\sim\mathsf{Unif}(\mathcal{C}_{\mathsf{RM}(r,m)})$ and $Z$ are independent. Despite $I(S;S+Z)$ being a well-behaved function of $\sigma$, an exact computation of this quantity is~infeasible.

Our estimator readily estimates $I(S;S+Z)$ from samples of $S$. Results for the Reed-Muller codes $\mathsf{RM}(4,4)$ and $\mathsf{RM}(5,5)$ (containing $2^{16}$ and $2^{32}$ codewords, respectively) are shown in Fig.~\ref{Fig:ReedMuller} for various values of $\sigma$ and $n$. Fig.~\ref{Fig:ReedMuller}(a) shows our estimate of $I(S;S+Z)$ for an $\mathsf{RM}(4,4)$ code as a function of $\sigma$, for different values of $n$. As expected, the plug-in estimator converges faster when $\sigma$ is larger.  Fig.~\ref{Fig:ReedMuller}(b) shows the estimated $I(S;S+Z)$ for $S\sim\mathsf{Unif}(\mathcal{C}_{\mathsf{RM}(5,5)})$ and $\sigma=2$, with the KDE and KL estimates based on samples of $(S+Z)$ shown for comparison. Our method significantly outperforms the competing general-purpose methods (with the wKL estimator being again omitted due to its instability in this high-dimensional ($d=32$) setting).

\begin{remark}[AWGN with Input Constraint]
When $\supp(P)$ lies inside a ball of radius $\sqrt{d}$, the subgaussian constant $K$ is proportional to $d$, and the bound from \eqref{EQ:Plugin_risk_bound} scales like $d^{d/2}n^{-1/2}$. This corresponds to the popular setup of an AWGN channel with an input constraint.
\end{remark}

\begin{remark}[Calculating the Ground Truth] To compute the true value of $I(S;S+Z)$ in Fig. \ref{Fig:ReedMuller}(b) (dashed red line) we used our MC integrator and the fact the Reed-Muller code is known. Specifically, the distribution of $S+Z$ is a Gaussian mixture, whose differential entropy we compute via the expression from \eqref{EQ:MCI}. Convergence of the computed value was ensured using Theorem \ref{TM:MC_MSE}.
\end{remark}

\section{Proofs for Section \ref{SEC:Results}}\label{SEC:proofs}


\subsection{Proof of Theorem \ref{TM:sample_complex_worstcase_asymp}}\label{SUPPSEC:sample_complex_worstcase_proof_asymp}




\subsubsection{\underline{Part 1}} Consider a AWGN channel $Y=X+N$, where the input $X$ is bound to a peak constraint $X\in[-1,1]$, almost surely, and $N\sim\mathcal{N}(0,\sigma^2)$ is an AWGN independent of $X$. The capacity (in nats) of this channel is 
\begin{equation}
    \mathsf{C}_\mathsf{AWGN}(\sigma)=\max_{X\sim P:\ P\in\mathcal{F}_d}I(X;Y),
\end{equation}
which is positive for any $\sigma<\infty$. The positivity of capacity implies the following \cite{Gallager_book1968}: for any rate $0<R < \mathsf{C}_\mathsf{AWGN}(\sigma)$, there exists a sequence of block codes (with blocklength $d$) of that rate, with an exponentially decaying (in $d$) maximal probability of error. More precisely, for any $\epsilon\in\big(0,\mathsf{C}_\mathsf{AWGN}(\sigma)\big)$, there exists a codebook $\mathcal{C}_d\subset[-1,1]^d$ of size $|\mathcal{C}_d|\doteq e^{d(\mathsf{C}_\mathsf{AWGN}(\sigma)-\epsilon)}$ and a decoding function $\psi_d:\mathbb{R}^d\to [-1,1]^d$ such that 
\begin{equation}
    \mathbb{P}\Big(\psi_d(Y^d)=c\Big| X^d=c\Big)\geq 1-e^{-\epsilon^2d},\quad \forall c\in\mathcal{C}_d,\label{EQ:maximal_error_prob}
\end{equation}
where $X^d\triangleq (X_1,X_2,\ldots,X_d)$ and $Y^d\triangleq(Y_1,Y_2,\ldots,Y_d)$ are the channel input and output sequences, respectively. The sign $\doteq$ stands for equality in the exponential scale, i.e., $a_k\doteq b_k$ means that $\lim_{k\to\infty}\frac{1}{k}\log\frac{a_k}{b_k}=0$. 

Since \eqref{EQ:maximal_error_prob} ensures an exponentially decaying error probability for any $c\in\mathcal{C}_d$, we also have that the error probability induced by a randomly selected codeword is exponentially small. Namely, let $X^d$ be a discrete random variable with any distribution $P$ over the codebook $\mathcal{C}_d$. We have
\begin{align*}
    \mathbb{P}\big(X^d\neq \psi_d(Y^d)\big) &=\sum_{c\in\mathcal{C}_d}P(c)\mathbb{P}\Big(\psi_d(c+N^d)\neq c\Big| X^d=c\Big)\\&\leq e^{-\epsilon^2d}.\numberthis\label{EQ:error_probability}
\end{align*}
Based on \eqref{EQ:error_probability}, Fano's inequality implies 
\begin{equation}
    H\p{X^d\middle|\psi_d(Y^d)}\leq H_b\p{e^{-\epsilon^2d}}+e^{-\epsilon^2d}\log|\mathcal{C}_d|\triangleq\delta^{(1)}_{\sigma,d},\label{EQ:Fano}
\end{equation}
where $H_b(\alpha)=-\alpha\log \alpha -(1-\alpha)\log(1-\alpha)$, for $\alpha\in[0,1]$, is the binary entropy function. Although not explicit in our notation, the dependence of $\delta^{(1)}_{\sigma,d}$ on $\sigma$ is through $\epsilon$. Note that $\lim_{d\to\infty}\delta^{(1)}_{\sigma,d}=0$, for all $\sigma>0$, because $\log|\mathcal{C}_d|$ grows only linearly with $d$ and $\lim_{q\to 0} H_b(q)=0$.


This further gives
\begin{align*}
        I\p{X^d;Y^d}&=H\p{X^d}-H\p{X^d\middle|Y^d}\\&\stackrel{(a)}\geq H\p{X^d}-H\p{X^d\middle|\psi_d(Y^d)}\\&\stackrel{(b)}\geq H\p{X^d}-\delta^{(1)}_{\sigma,d},\numberthis\label{EQ:sample_complex_MIexpansion1}
\end{align*}
where (a) follows because $H(A|B)\leq H\big(A\big|f(B)\big)$ for any pair of random variables $(A,B)$ and any deterministic function $f$, while (b) uses \eqref{EQ:Fano}.

Non-negativity of discrete entropy also implies $I(X^d;Y^d)\leq H(X^d)$, which means that $H(X^d)$ and $I(X^d;Y^d)$ become arbitrarily close as $d$ grows:
\begin{equation}
    \Big|H(X^d)-I(X^d;Y^d)\Big|\leq \delta^{(1)}_{\sigma,d}.\label{EQ:sample_complex_gap}
\end{equation}
This means that any good estimator (within an additive gap) of $H(X^d)$ over the class of distributions $\big\{P\mspace{3mu}\big|\supp(P)=\mathcal{C}_d\big\}\subseteq\mathcal{F}_d$ is also a good estimator of the mutual information. Using the well-known lower bound on the sample complexity of discrete entropy estimation in the large alphabet regime (see, e.g., \cite[Corollary 10]{valiant2010clt} or \cite[Proposition 3]{discrete_entropy_est_Wu2016}), we have that estimating $H(X^d)$ within a small additive gap $\eta>0$ requires at~least
\begin{equation}
\Omega\p{\frac{|\mathcal{C}_d|}{\eta\log|\mathcal{C}_d|}}=\Omega\p{\frac{2^{\gamma(\sigma)d}}{\eta d}},
\end{equation}
where $\gamma(\sigma)\triangleq\mathsf{C}_\mathsf{AWGN}(\sigma)-\epsilon>0$ is independent of $d$.

We relate the above back to the considered differential estimation setup by noting that 
\begin{align*}
    I(X^d;Y^d)&= h(X^d+N^d)-h(N^d)\\&=h(X^d+N^d)-\frac{d}{2}\log_2(2\pi e \sigma^2).\numberthis\label{EQ:sample_complex_MIexpansion2}
\end{align*}
Letting $S\sim P$ and noting that $Z\stackrel{\mathcal{D}}=N^d$, where $\stackrel{\mathcal{D}}=$ denotes equality in distribution, we have $h(X^d+N^d)=h(S+Z)$. Assuming in contradiction that there exists an estimator of $h(S+Z)$ that uses $o\big(2^{\gamma(\sigma)d}/(\eta d)\big)$ samples and achieves an additive gap $\eta>0$ over $\big\{P\mspace{3mu}\big|\supp(P)=\mathcal{C}_d\big\}$, implies that $H(X^d)$ can be estimated from these samples within gap $\eta+\delta^{(1)}_{\sigma,d}$. This follows from \eqref{EQ:sample_complex_gap} by taking the estimator of $h(S+Z)$ and subtracting the constant $\frac{d}{2}\log_2(2\pi e \sigma^2)$. We arrive at a contradiction.

\subsubsection{\underline{Part 2}}

Fix $d\geq 1$ and consider a $d$-dimensional AWGN channel $Y=X+N$, with input $X$ and noise $N\sim\mathcal{N}(0,\sigma^2\mathrm{I}_d)$. Let $\mathcal{C}=\{-1,1\}^d$ and consider the set of all (discrete) distributions $P$ with $\supp(P)=\mathcal{C}$. For $X\sim P$, with $P$ being an arbitrary distribution from the aforementioned set, and any mapping $\psi_\mathcal{C}:\mathbb{R}^d\to\mathcal{C}$, Fano's inequality gives
\begin{equation}
    H(X|Y)\leq H\big(X\big|\psi_\mathcal{C}(Y)\big)\leq H_b\big(\mathsf{P}_\mathsf{e}(\mathcal{C})\big)+\mathsf{P}_\mathsf{e}(\mathcal{C})\cdot \log|\mathcal{C}|,\label{EQ:Fano_smalld}
\end{equation}
where $\mathsf{P}_{\mathsf{e}}(\mathcal{C})\triangleq \mathbb{P}\Big(\psi_\mathcal{C}(Y)\neq X\Big)$ is the error probability. We choose $\psi_\mathcal{C}$ as the maximum likelihood decoder: upon observing $y\in\mathbb{R}^d$ it returns the closest point in $\mathcal{C}$ to $y$. Namely, $\psi_\mathcal{C}$ returns $c\in\mathcal{C}$ if and only if $y$ falls inside the unique orthant that contains $c$. We have:
\begin{align*}
    \mathsf{P}_{\mathsf{e}}(\mathcal{C})&=\sum_{c\in\mathcal{C}}P(c)\mathbb{P}\Big(\psi_\mathcal{C}(c+Z)\neq c\Big|X=c\Big)\\&=1 - \left(1 - Q\left(\frac{1}{\sigma}\right)\right)^d\triangleq\epsilon_{\sigma,d},\numberthis\label{EQ:Error_Prob_UB}
\end{align*}
where $Q$ is the Q-function. Together, \eqref{EQ:Fano_smalld} and \eqref{EQ:Error_Prob_UB} give
$H(X|Y)\leq H_b(\epsilon_{\sigma,d})+\epsilon_{\sigma,d} d\log 2 \triangleq \delta^{(2)}_{\sigma,d}$. Note that for any $d\geq 1$, $\lim_{\sigma\to 0} \delta^{(2)}_{\sigma,d}=0$ exponentially fast in $\frac{1}{\sigma^2}$ (this follows from the large $x$ approximation of $Q(x)$). Similarly to \eqref{EQ:sample_complex_gap}, the above implies that
\begin{equation}
    \Big|H(X)-I(X;Y)\Big|\leq\delta^{(2)}_{\sigma,d}.\label{EQ:sample_complex_gap_smalld}
\end{equation}
Thus, any good estimator (within an additive gap $\eta$) of $H(X)$ within the class of $X$ distributions $P$ with $\supp(P)=\mathcal{C}$, can be used to estimate $I(X;Y)$ within an $\eta+\delta^{(2)}_{\sigma,d}$ gap.

Now, for $\sigma$ small enough $\epsilon_{\sigma,d}$, and consequently $\delta^{(2)}_{\sigma,d}$ are arbitrarily close to zero. Hence we may again use lower bounds on the sample complexity of discrete entropy estimation. Like in the proof of Theorem \ref{TM:sample_complex_worstcase_asymp}, setting $S\sim P$, any estimator of $h(S+Z)$ within a small gap $\eta$ produces an estimator of $H(X)$ (through $H(X)=h(S+Z)-\frac{d}{2}\log(2\pi e\sigma^2)$ and \eqref{EQ:sample_complex_gap_smalld}) within an $\eta+\delta^{(2)}_{\sigma,d}$ gap. Therefore, for sufficiently small $\sigma>0$ and $\eta>0$, any estimator of $h(S+Z)$ within a gap of $\eta$ requires at least
\begin{equation}
\Omega\left(\frac{\supp(P)}{\big(\eta+\delta^{(2)}_{\sigma,d}\big)\log\big(\supp(P)\big)}\right)=\Omega\left(\frac{2^{d}}{\big(\eta+\delta^{(2)}_{\sigma,d}\big) d}\right)
\end{equation}
samples. This concludes the proof.


\subsection{Proof of Theorem \ref{TM:SP_sample_complex_new_bdd}}\label{SUBSEC:SP_sample_complex_new_bdd_proof}

We start with the following lemma. 

\begin{lemma}\label{LEMMA:entropy_bound}
Let $U\sim P_U$ and $V\sim P_V$ be continuous random variables with densities $p_U$ and $p_V$, respectively. If $\big|h(U)\big|,\big|h(V)\big|<\infty$, then
\begin{equation*}
    \big|h(U)-h(V)\big|\leq \max\left\{\left|\E\log\frac{p_V(V)}{p_V(U)}\right|,\left|\E\log\frac{p_U(U)}{p_U(V)}\right|\right\}.
\end{equation*}
\end{lemma}

\begin{IEEEproof}
Recall the identity
\begin{align*}
h(U) - h(V) &\leq h(U)-h(V)+D(P_U||P_V)\\&=\E\log \frac{p_V(V)}{p_V(U)}\leq\left|\E\log \frac{p_V(V)}{p_V(U)}\right|.
\end{align*}
Reversing the roles of $U$ and $V$ in the above derivation establishes the second bound and completes the proof. 
\end{IEEEproof}

Recall now the variational characterization of the $\Chi$-divergence:
\begin{equation}
    \Chi(\mu\|\nu)=\sup_{g:\ \mathsf{var}_\nu(g)\leq 1}\big|\E_\mu g-\E_\nu g\big|^2.\label{EQ:Chi_squared_variational}
\end{equation}
Combining this with Lemma \ref{LEMMA:entropy_bound}, we obtain
\begin{align*}
    \big|h(U)&-h(V)\big|\\&\leq\max\Big\{\sqrt{\mathsf{var}_{P_V}\big(\log p_V(V)\big)\Chi(P_U\|P_V)},\\&\quad\qquad\qquad\sqrt{\mathsf{var}_{P_V}\big(\log p_U(V)\big)\Chi(P_U\|P_V)}\Big\}.\numberthis\label{EQ:ent_Chi_bound}
\end{align*}
Setting $P_V=P\ast\Gauss$ and $P_U=\hat{P}_{S^n}\ast\Gauss$, the next lemma is useful in controlling the variance terms. To state it recall that $q$ and $r_{S^n}$ are the PDFs of $P\ast\Gauss$ and $\hat{P}_{S^n}\ast\Gauss$, respectively, and set 
$\tilde q\triangleq \frac{q}{c_1}$ and $\tilde{r}_{S^n}\triangleq\frac{r}{c_1}$ for $c_1 = (2 \pi \sigma^2)^{-d/2}$.

\begin{lemma}\label{LEMMA:log_bound}
Let $S \sim P$. For all $z \in \RR^d$ it holds that
\begin{subequations}
\begin{align}
\E_{P^{\otimes n}}\big(\log \tilde r_{S^n}(z)\big)^2 & \leq \frac{1}{4 \sigma^4} \E_P \|z - S\|^4\label{EQ:Elog_bound1} \\
\big(\log \tilde q(z)\big)^2 & \leq \frac{1}{4 \sigma^4} \E_P \|z - S\|^4.\label{EQ:Elog_bound2}
\end{align}
\end{subequations}
\end{lemma}

\begin{IEEEproof}
We prove \eqref{EQ:Elog_bound1}; the proof of \eqref{EQ:Elog_bound2} is similar and therefore omitted. The map $x \mapsto (\log x)^2$ is convex on $[0,1]$. For any fixed $s^n$, let $\hat{S}\sim \hat{P}_{s^n}$. Jensen's~inequality gives
\begin{align*}
\big(\log \tilde r_{s^n}(z)\big)^2 &= \left(\log \E_{\hat P_{s^n}} \exp\left(-\frac{\|z-\hat S\|^2}{2 \sigma^2}\right)\right)^2\\&\leq\E_{\hat P_{s^n}}\frac{\|z-\hat S\|^4}{4 \sigma^4}.
\end{align*}
Taking an outer expectation w.r.t. $S^n\sim P^{\otimes n}$ yields
\begin{align*}
\E_{P^{\otimes n}} \big(\mspace{-2mu}\log \tilde r_{S^n}(z)\big)^2\leq \E_{P^{\otimes n}}\E_{\hat P_{S^n}}\frac{\|z - \hat{S}\|^4}{4 \sigma^4}= \frac{\E_P \|z - S\|^4}{4 \sigma^4}.
\end{align*}
\end{IEEEproof}

Let $Y=S'+Z$, where $S'\sim P$ and $Z\sim\Gauss$ are independent. Since variance is translation invariant, we get
\begin{align*}
    \mathsf{var}_{P\ast\Gauss}\big(\log q(Y)\big)&=\mathsf{var}_{P\ast\Gauss}\big(\log \tilde{q}(Y)\big)\\&\leq \frac{1}{4 \sigma^4} \E \|Z + S' - S\|^4\\&\leq\frac{\sigma^2 d(2+d)(2+\sigma^2) + 8 d^2}{4\sigma^4}.\numberthis\label{EQ:var_bound1}
\end{align*}
When combined with Proposition \ref{PROP:Chi_parametric}, the above bound takes care of the first term in \eqref{EQ:ent_Chi_bound}.

For the second term, we apply Cauchy-Schwartz and treat the expected values of $\mathsf{var}_{P_V}\big(\log p_U(V)\big)$ and $\Chi(P_U\|P_V)$ separately. For the variance, using \eqref{EQ:Elog_bound1} and an argument similar to \eqref{EQ:var_bound1} we get the same bound therein. The expected $\Chi$-square divergence in both arguments of the maximum in \eqref{EQ:ent_Chi_bound} is bounded using Corollary \ref{CORR:Chi_parametric_bdd}. Combining the pieces, for any $P\in\mathcal{F}_d$, we obtain
\begin{align*}
    \E_{P^{\otimes n}}&\big|h(P\ast\Gauss)-h(\hat{P}_{S^n}\ast\Gauss)\big|\\&\leq 2 \sqrt{\frac{\sigma^2 d(2+d)(2+\sigma^2) + 8 d^2}{4\sigma^4}}e^{\frac{2d}{\sigma^2}}\cdot\frac{1}{\sqrt{n}}.\numberthis\label{EQ:Plugin_risk_bound_constant_bdd}
\end{align*}

\begin{remark}
An alternative proof of the parametric estimation rate was given in \cite{Weed_bdd_Support2018} using the 1-Wasserstein distance instead of $\Chi$-square. Specifically, one may invoke \cite[Proposition 5]{PolyWu_Wasserstein2016} to reduce the analysis of $\E_{P^{\otimes n}}\big|h(P\ast\Gauss)-h(\hat{P}_{S^n}\ast\Gauss)\big|$ to that of $\E_{P^{\otimes n}}\mathsf{W}_1\left(\hat{P}_{S^n}\ast\Gauss,P\ast\Gauss\right)$. Then, using \cite[Theorem~6.15]{villani2008optimal} and the bounded support assumption, the parametric risk convergence rate follows with the constant $\frac{\sqrt d \cdot 2^{d+2}}{\min\{1, \sigma^d\}}$.
\end{remark}

\subsection{Proof of Theorem \ref{TM:SP_sample_complex_new}}\label{SUBSEC:plug-in_risk_bound_proof}

Starting from Lemma \ref{LEMMA:entropy_bound}, we again focus on bounding the maximum of the two expected log ratios. The following lemma allows converting 
$\left|\E\log\frac{p_V(V)}{p_V(U)}\right|$ and $\left|\E\log\frac{p_U(U)}{p_U(V)}\right|$ into forms that are more convenient to analyze.

\begin{lemma}\label{LEMMA:test-g}
Let $U\sim P_U$ and $V\sim P_V$ be continuous random variables with PDFs $p_U$ and $p_V$, respectively. For any measurable function $g : \mathbb{R}^d  \to \mathbb{R}$
\begin{equation*}
\big|\E g(U)-\E g(V)\big| \leq \int \big|g(z)\big| \cdot \big|p_U(z) - p_V(z)\big| \dd z\ .
\end{equation*}
\end{lemma}

\begin{IEEEproof} We couple $P_U$ and $P_V$ via the maximal TV coupling\footnote{This coupling attains maximal probability of the event $\{U=V\}$.}. Specifically, let $(P_U - P_V)_+$ and $(P_U - P_V)_-$ are the positive and negative parts of the signed measure $(P_U-P_V)$. Define $(P_U\wedge P_V)\triangleq P_U-(P_U-P_V)_+$, and let $(\mathrm{Id},\mathrm{Id})_\sharp (P_U \wedge P_V)$ be the push-forward measure of $P_U \wedge P_V$ through $(\mathrm{Id}, \mathrm{Id}):\RR^d\to\RR^d\times\RR^d$. Letting $\alpha \triangleq \frac{1}{2}\int|p_U(x) - p_V(x)|dx$ (note that $\int \dd(P_U - P_V)_+ = \int \dd(P_U - P_V)_- = \alpha$), the maximal TV coupling is give by 
\begin{equation}
\pi \triangleq (\mathrm{Id}, \mathrm{Id})_\sharp (P_U \wedge P_V) + \frac{1}{\alpha} (P_U - P_V)_+ \otimes (P_U-P_V)_-.\label{EQ:maximal_coupling}
\end{equation}
Jensen's inequality implies $\big|\E g(U)-\E g(V)\big|\leq \E_\pi\big|g(U)-g(V)\big|$ and we proceed as
\begin{align*} 
&\E_\pi\big|g(U)-g(V)\big|\\
& \leq \frac 1 \alpha \int\Bigg( \Big(\mspace{-2mu}\big|g(u)\big|+\big|g(v)\big|\Big)\big(p_U(u)-p_V(u)\big)_+\\&\qquad\qquad\qquad\qquad\qquad\qquad \cdot\big(p_U(v)-p_V(v)\big)_-\Bigg)\dd u\dd v \\
& =\int \big|g(u)\big|\big(p_U(u)-p_V(u)\big)_+\dd u\\&\qquad\qquad+\int\big|g(v)\big|\big(p_U(v)-p_V(v)\mspace{-2mu}\big)_- \dd v \\
& = \int \big|g(z)\big| \Big(\big(p_U(z) - p_V(z)\big)_+ + (p_U(z) - p_V(z)\big)_-\Big) \dd z \\
&= \int \big|g(z)\big|\cdot \big|p_U(z) - p_V(z)\big| \dd z.\numberthis
\end{align*}
\end{IEEEproof}

Fix any $P\in\FSG$ and assume that $\mathbb{E}_PS=0$. This assumption comes with no loss of generality since both the target functional $h(P\ast\Gauss)$ and the plug-in estimator are translation invariant. Note that $\big|h(P\ast\Gauss)\big|,\big|h(\hat{P}_{S^n}\ast\Gauss)\big|<\infty$. Combining Lemmata \ref{LEMMA:entropy_bound}
and \ref{LEMMA:test-g}, we a.s. have
\begin{align*}
\big|h(&P\ast\Gauss)-h(\hat{P}_{S^n}\ast\Gauss)\big| \\&\leq \max\Bigg\{\int\big|\log \tilde r_{S^n}(z)\big|\cdot |q(z) - r_{S^n}(z)| \dd z,\\&\hskip 19mm \int  \big|\log \tilde q(z)\big|\cdot \big|q(z) - r_{S^n}(z)\big| \dd z\Bigg\},\numberthis\label{EQ:ent_bound_lemmata}
\end{align*}
where, as before, $q$ and $r_{S^n}$ are the PDFs of $P\ast\Gauss$ and $\hat{P}_{S^n}\ast\Gauss$, respectively, while $\tilde q \triangleq \frac{q}{c_1}$ and $\tilde r_{S^n} \triangleq \frac{r_{S^n}}{c_1}$, for $c_1 = (2 \pi \sigma^2)^{-d/2}$.

Recalling that $\E[\max\{|X|,|Y|\}]\leq \E |X|+ \E|Y|$, for any random variable $X$, $Y$, we now bound $\int\big|\log \tilde r_{S^n}(z)\big|\big|p_U(z)-p_V(z)\big|\dd z$. The bound for the other integral is identical and thus omitted. Let $f_a:\RR^d \to \RR$ be the PDF of $\mathcal{N}\left(0,\frac{1}{2a} \mathrm{I}_d\right)$, for $a > 0$ specified later.
The Cauchy-Schwarz inequality implies
\begin{align*}
&\left(\mathbb{E}_{P^{\otimes n}}\int\big|\log\tilde r_{S^n}(z)\big| \big|q(z) - r_{S^n}(z)\big| \dd z\right)^2 \leq
\\&\int\E_{P^{\otimes n}} \big(\log \tilde r_{S^n}(z)\big)^2f_a(z)\dd z\!\cdot\!\!\int \E_{P^{\otimes n}} \frac{\big(q(z)-r_{S^n}(z)\big)^2}{f_a(z)} \dd z.\numberthis\label{EQ:CS_bound}
\end{align*}
Using Lemma \ref{LEMMA:log_bound}, we bound the first integral as
\begin{align*}
\int\E_{P^{\otimes n}}&\big( \log \tilde r_{S^n}(z)\big)^2f_a(z)\dd z\\&\leq   \int\frac{\E \|z - S\|^4}{4 \sigma^4} \frac{\exp\big(-a \|z\|^2\big)}{\sqrt{ \pi^d a^{-d}}} \dd z  \\
&\stackrel{(a)} \leq \frac{2}{\sigma^4} \E \|S\|^4 + \frac{2}{\sigma^4}\int \|z\|^4 \frac{\exp\big(-a \|z\|^2\big)}{\sqrt{ \pi^d a^{-d}}} \dd z \\
&\stackrel{(b)} \leq  \frac{32K^4d^2}{\sigma^4} + \frac{1}{2 \sigma^4 a^2} d(d+2)
\end{align*}
where (a) follows from the triangle inequality, and (b) uses the $K$-subgaussianity of $S$ \cite[Lemma~5.5]{vershynin2010introduction}


To bound the second integral, we repeat steps \eqref{EQ:W1_1/f_bound}-\eqref{EQ:W1_int2_bound} from the proof of Proposition \ref{PROP:W1_parametric}. Specifically, we have $\mathbb{E}_{P^{\otimes n}}\big(q(z)-r_{S^n}(z)\big)^2\mspace{-3mu}\leq\mspace{-3mu}\frac{c_1^2}{n}\E e^{-\frac{1}{\sigma^2} \|z - S\|^2}$, because $r_{S^n}(z)$ is a sum of i.i.d. random variables with $\E_{P^{\otimes n}}r_{S^n}(z)=q(z)$. This gives
\begin{equation}
\int\E_{P^{\otimes n}} \frac{\big(q(z)-r_{S^n}(z)\big)^2}{f_a(z)} \dd z \leq \frac{c_1}{n2^{d/2}} \E \frac{1}{f_a(S+ Z/\sqrt{2} )},\label{EQ:1/f_bound}
\end{equation}
for independent $Z\sim\cN(0, \sigma^2\mathrm{I_d})$ and $S\sim P$. Recalling that $\big(f_a(z)\big)^{-1}= c_2 \exp\big(a \|z\|^2\big)$, for $c_2\triangleq\left(\frac{\pi}{a}\right)^{\frac{d}{2}}$, the subgaussianity of $S$ and $Z$ implies
\begin{align*}
    &\frac{c_1}{n2^{d/2}}\E \frac{1}{f_a(S + Z/\sqrt{2})}\leq\\& \frac{c_1 c_2}{n2^{d/2}} \exp\!\left(\!\big(K+\sigma/ \sqrt 2\big)^2 a d\!+\!\frac{(K+\sigma/ \sqrt 2)^4 a^2d}{1 -2(K+\sigma/ \sqrt 2)^2 a}\!\right)\!,\numberthis\label{EQ:int2_bound}
\end{align*}
where $0<a<\frac{1}{2(K+\sigma/\sqrt{2})^2}$.

Setting $a\mspace{-4mu}=\mspace{-4mu}\frac{1}{4(K+\sigma/\sqrt{2})^2}$, we combine \eqref{EQ:CS_bound}-\eqref{EQ:int2_bound} to obtain the result (recalling that the second integral from \eqref{EQ:ent_bound_lemmata} is bounded exactly as the first). 
For any $P\in\FSG$ we have
\begin{align*}
&\bigg(\E_{P^{\otimes n}}\big|h(P\ast\Gauss)-h(\hat{P}_{S^n}\ast\Gauss)\big|\bigg)^2
\leq\\& 
\frac{64\big(2 d^2K^4 + d(d+2)( K + \sigma/\sqrt{2})^4\big)}{\sigma^4}
\left(\!\left(\frac{1}{\sqrt{2}}+\frac{K}{\sigma}\right)e^{\frac{3}{8}}\right)^d\!\frac{1}{n}\ .\numberthis\label{EQ:Plugin_risk_bound_constant}
\end{align*}


\subsection{Proof of Theorem \ref{TM:SP_Bias_LB}}\label{SUBSEC:SP_Bias_LB}

First note that since $h(q)$ is concave in $q$ and because $\mathbb{E}_{P^{\otimes n}}\hat{P}_{S^n}=P$, we have
\begin{equation}
    \mathbb{E}_{P^{\otimes n}} h(\hat{P}_{S^n}\ast \gauss) \leq h(P\ast\Gauss),\label{EQ:SP_LB_truth}
\end{equation}
for all $P\in\mathcal{F}_d$. Now, let $W\sim\mathsf{Unif}([n])$ be independent of $(S^n,Z)$ and define $Y=S_W+Z$. We have the following lemma, whose proof is found in Appendix \ref{APPEN:SP_bias_MI_proof}.

\begin{lemma}\label{LEMMA:SP_bias_MI}
For any $P\in\mathcal{F}_d$, we have
\begin{equation}
    h(P\ast\Gauss)-\mathbb{E}_{P^{\otimes n}}h(\hat{P}_{S^n}\ast \Gauss)=I(S^n;Y).
\end{equation}
\end{lemma}

Using the lemma, we have
\begin{equation}
    \sup_{P\in\mathcal{F}_d}\big|h(P\ast\Gauss)-\mathbb{E}_{P^{\otimes n}}h(P\ast\Gauss)\big|=\sup_{P\in\mathcal{F}_d}I(S^n;Y),\label{EQ:Bias_MI_relation}
\end{equation}
where the right hand side is the mutual information between $n$ i.i.d. random samples $S_i$ from $P$ and the random vector $Y=S_W+Z$, formed by choosing one of the $S_i$'s at random and adding Gaussian noise. 

To obtain a lower bound on the supremum, we consider the following $P$. Partition the hypercube $[-1,1]^d$ into $k^d$ equal-sized smaller hypercubes, each of side length $k$. Denote these smaller hypercubes as $\mathsf{C}_1,\mathsf{C}_2,\ldots,\mathsf{C}_{k^d}$ (the order does not matter). For each $i\in[k^d]$ let $c_i\in\mathsf{C}_i$ be the centroid of the hypercube $\mathsf{C}_i$. Let $\mathcal{C}\triangleq\{c_i\}_{i=1}^{k^d}$ and choose $P$ as the uniform distribution over $\mathcal{C}$.

By the mutual information chain rule and the non-negativity of discrete entropy, we have
\begin{align*}
    I(S^n;Y)&=I(S^n;Y,S_W)-I(S^n;S_W|Y)\\
    &\stackrel{(a)}\geq I(S^n;S_W)-H(S_W|Y)\\
    &= H(S_W)-H(S_W|S^n)-H(S_W|Y),\numberthis\label{EQ:MI_LB_SP_bias}
\end{align*}
where step (a) uses the independence of $(S^n,W)$ and $Z$. Clearly $H(S_W)=\log|\mathcal{C}|$, while $H(S_W|S^n)\leq H(S_W,W|S^n)\leq H(W)=\log n$, via the independence of $W$ and $S^n$. For the last (subtracted) term in \eqref{EQ:MI_LB_SP_bias} we use Fano's inequality to obtain
\begin{equation}
    H(S_W|Y)\leq H\big(S_W\big|\psi_\mathcal{C}(Y)\big)\leq H_b\big(\mathsf{P}_\mathsf{e}(\mathcal{C})\big)+\mathsf{P}_\mathsf{e}(\mathcal{C})\cdot\log|\mathcal{C}|,
\end{equation}
where $\psi_\mathcal{C}:\mathbb{R}^d\to\mathcal{C}$ is a function for decoding $S_W$ from $Y$ and $\mathsf{P}_\mathsf{e}(\mathcal{C})\triangleq \mathbb{P}\big(S_W\neq\psi_\mathcal{C}(Y)\big)$ is the probability that $\psi_\mathcal{C}$ commits an error. 

Fano's inequality holds for any decoding function $\psi_\mathcal{C}$. We choose $\psi_\mathcal{C}$ as the maximum likelihood decoder, i.e., upon observing a $y\in\mathbb{R}^d$ it returns the closest point to $y$ in $\mathcal{C}$. Denote by $\mathcal{D}_i\triangleq\psi_\mathcal{C}^{-1}(c_i)$ the decoding region on $c_i$, i.e., the region $\cp{y \in \mathbb{R}^d \middle| \psi_\mathcal{C}(y) = c_i}$ that $\psi_\mathcal{C}$ maps to $c_i$. Note that $\mathcal{D}_i=\mathsf{C}_i$ for all $i\in[k^d]$ for which $\mathsf{C}_i$ doesn't intersect with the boundary of $[-1,1]^d$. The probability of error for the decoder $\psi_\mathcal{C}$ is bounded as:
\begin{align*}
    \mathsf{P}_{\mathsf{e}}(\mathcal{C})&=\frac{1}{k^d}\sum_{i=1}^{k^d}\mathbb{P}\Big(\psi_\mathcal{C}(c_i+Z)\neq c_i\Big|S_W=c_i\Big)\\
    &=\frac{1}{k^d}\sum_{i=1}^{k^d}\mathbb{P}\big(c_i+Z\notin\mathcal{D}_i\big)\\
    &\stackrel{(a)}\leq \mathbb{P}\left(\|Z\|_\infty > \frac{2/k}{2}\right)\\
    &\stackrel{(b)} = 1 - \left(1 - 2 Q\left(\frac{1}{k\sigma}\right)\right)^d,
    \numberthis\label{EQ:Error_Prob_UB_bias}
\end{align*}
where (a) holds since the $\mathsf{C}_i$ have sides of length $2/k$ and the error probability is largest for $i\in[k^d]$ such that $\mathsf{C}_i$ is in the interior of $[-1,1]^d$. Step (b) follows from independence and the definition of the Q-function.


Taking $k=k_\star$ in \eqref{EQ:Error_Prob_UB_bias} as given in the statement of the theorem
gives the desired bound $\mathsf{P}_{\mathsf{e}}(\mathcal{C})\leq \epsilon$. Collecting the pieces and inserting back to \eqref{EQ:MI_LB_SP_bias}, we obtain
\begin{equation}
    I(S^n;Y)\geq \log\left(\frac{k_\star^{d(1-\epsilon)}}{n}\right)-H_b(\epsilon).
\end{equation}
Together with \eqref{EQ:Bias_MI_relation} this concludes the proof.


\subsection{Proof of Theorem \ref{TM:MC_MSE}}\label{SUBSEC:MC_MSE_proof}

Denote the joint distribution of $(C,Z,V)$ by $P_{C,Z,V}$. Marginal or conditional distributions are denoted as usual by keeping only the relevant subscripts. Lowercase $p$ denotes a probability mass function (PMF) or a PDF depending on whether the random variable in the subscript is discrete or continuous. In particular, $p_C$ is the PMF of $C$, $p_{C|V}$ is the conditional PMF of $C$ given $V$, while $p_Z=\gauss$ and $p_V=g$ are the PDFs of $Z$ and $V$, respectively.

First observe that the estimator is unbiased:
\begin{equation}
    \mathbb{E}\hat{h}_\mathsf{MC}=-\frac{1}{n\cdot \nmc}\sum_{i=1}^n\sum_{j=1}^{\nmc}\mathbb{E}\log g\left(\mu_i+Z_j^{(i)}\right)=h(g). 
\end{equation}
Therefore, the MSE expands as
\begin{equation}
    \mathsf{MSE}\left(\hat{h}_\mathsf{MC}\right)=\frac{1}{n^2\cdot \nmc}\sum_{i=1}^n\var\Big(\log g(\mu_i+Z)\Big).\label{EQ:MC_MSE_expansion}
\end{equation}

We next bound the variance of $\log g(\mu_i+Z)$ via the Gaussian Poincar{\'e} inequality (with Poincar{\'e} constant $\sigma^2)$. For each $i\in[n]$, we have
\begin{equation}
    \var\Big(\log g(\mu_i+Z)\Big)\leq \sigma^2 \mathbb{E}\Big[\big\|\nabla\log g(\mu_i+Z)\big\|^2\Big].\label{EQ:Poincare}
\end{equation}
We proceed with separate derivations of \eqref{EQ:MC_MSE_Tanh} and \eqref{EQ:MC_MSE_ReLU}. 

\vspace{2mm}
\subsubsection{\underline{MSE for Bounded Support}} Since $\|C\|_2\leq \sqrt{d}$ almost surely, Proposition 3 from \cite{PolyWu_Wasserstein2016} implies
\begin{equation}
    \big\|\nabla\log g(v)\big\|_2\leq \frac{\|v\|+\sqrt{d}}{\sigma^2}.
\end{equation}
Inserting this into the Poincar{\'e} inequality and using $(a+b)^2\leq 2a^2+2b^2$ we have,
\begin{equation}
    \var\Big(\log g(\mu_i+Z)\Big)\leq\frac{2d(4+\sigma^2)}{\sigma^2},
\end{equation}
for each $i\in[n]$. Together with \eqref{EQ:MC_MSE_expansion}, this produces \eqref{EQ:MC_MSE_Tanh}.

\vspace{2mm}
\subsubsection{\underline{MSE for Bounded Second Moment}}

To prove \eqref{EQ:MC_MSE_ReLU}, we use Proposition 2 from \cite{PolyWu_Wasserstein2016} to obtain
\begin{equation}
    \big\|\nabla\log g(v)\big\|\leq \frac{1}{\sigma^2}\big(3\|v\|+4\mathbb{E}\|C\|\big).
\end{equation}
Using \eqref{EQ:Poincare}, the variance is bounded as
\begin{align*}
   & \var\Big(\log g(\mu_i+Z)\Big)\leq \frac{1}{\sigma^2}\mathbb{E}\Big[\big(3\|\mu_i+Z\|+4\mathbb{E}\|C\|\big)^2\Big]\leq\\
    & \frac{1}{\sigma^2}\!\left(9d\sigma^2\!+\!16m\!+\!24\sigma\sqrt{dm}\!+\!3\|\mu_i\|\left(3\!+\!9\sigma\sqrt{d}\!+\!8\sigma\sqrt{dm}\right)\right),\numberthis\label{EQ:MC_MSE_ReLU_var_bound}
\end{align*}
where the last step uses H{\"o}lder's inequality $\mathbb{E}\|C\|\leq \sqrt{\mathbb{E}\|C\|^2}$. The proof of \eqref{EQ:MC_MSE_ReLU} is concluded by plugging \eqref{EQ:MC_MSE_ReLU_var_bound} into the MSE expression from \eqref{EQ:MC_MSE_expansion} and noting that $\frac{1}{n}\sum_{i=1}^n\|\mu_i\|\leq \sqrt{m}$.

\subsection{Proof of Proposition \ref{PROP:MI_True_Data_Dist}}\label{SUBSEC:MI_True_Data_Dist_proof}

Fix $P_X$, define $g(x)\triangleq h(T|X=x)=h(P_{S|X=x}\ast\Gauss)$ and write
\begin{equation}
    I(X;T)=h(T)-h(T|X)=h(P_S\ast\Gauss)-\mathbb{E}g(X).
\end{equation}
Applying the triangle inequality to \eqref{EQ:MI_CONT} we obtain
\begin{align*}
\mathbb{E}\bigg|I(X;T)&-\hat{I}_{\mathsf{Input}}\p{X^n,\hat{h},\sigma}\bigg|\\
&\leq \mathbb{E} \left|\hat{h}(S^n,\sigma) - h(P_S\ast\Gauss) \right| \\&\hskip 15mm + \mathbb{E}\left|\frac{1}{n} \sum_{i=1}^n \hat{h}\big(S^n(X_i),\sigma\big) - \mathbb{E} g(X)\right|\\
&\leq \underbrace{\mathbb{E} \left|\hat{h}(S^n,\sigma) - h(P_S\ast\Gauss) \right|}_{(\mathrm{I})} \\
&\hskip 10mm + \underbrace{\frac{1}{n} \sum_{i=1}^n \mathbb{E} \left| \hat{h}\big(S^n(X_i),\sigma\big) -  g(X_i)\right|}_{(\mathrm{II})}\\
&\hskip 20mm + \underbrace{\mathbb{E}\left|\frac{1}{n} \sum_{i=1}^n g(X_i)- \mathbb{E} g(X)\right|}_{(\mathrm{III})}\numberthis\label{eq:BarError}
\end{align*}
By assumption \eqref{EQ:Data_Dist_Assumption} and because $P_S\in\mathcal{F}_d$, we have
\begin{equation}
\mathbb{E} \left|\hat{h}(S^n,\sigma) - h(P_S\ast\Gauss)\right| \leq \Delta_{\sigma,d}(n).\label{EQ:FIRSTTERM}
\end{equation}
Similarly, for any fixed $X^n=x^n$, $P_{S|X=x_i}\in\mathcal{F}_d$, for all $i\in[n]$, and hence
\begin{align*}
    \mathbb{E}\bigg[&\left| \hat{h}(S^n(X_i),\sigma)- g(X_i)\right|\bigg|X^n=x^n\bigg]\\&\stackrel{(a)}=\mathbb{E}\left| \hat{h}\big(S^n(x_i),\sigma\big)- h(P_{S|X=x_i}\ast\Gauss)\right|\\&\leq \Delta_{\sigma,d}(n),\numberthis\label{EQ:FIRSTTERM_temp}
\end{align*}
where (a) is because for a fixed $x_i$, sampling from $P_{S|X=x_i}$ corresponds to drawing multiple noise realization for the previous layers of the DNN. Since these noises are independent of $X$, we may remove the conditioning from the expectation. Taking an expectation on both sides of \eqref{EQ:FIRSTTERM_temp} and applying the law of total expectation, we have
\begin{equation}
(\mathrm{II}) = \frac{1}{n} \sum_{i=1}^n \mathbb{E} \left|\hat{h}\big(S^n(x_i),\sigma\big)- g(X_i)\right|\leq\Delta_{\sigma,d}(n).\label{EQ:SECONDTERM}
\end{equation}

Turning to term $(\mathrm{III})$, observe that $\big\{g(X_i)\big\}_{i=1}^n$ are i.i.d random variables. Hence 
\begin{equation}
\frac{1}{n} \sum_{i=1}^n g(X_i)- \mathbb{E}g(X)
\end{equation}
is the difference between an empirical average and the expectation. By monotonicity of moments we have
\begin{align*}
(\mathrm{III})^2 &= \left(\mathbb{E} \left|\frac{1}{n} \sum_{i=1}^n g(X_i)- \mathbb{E}g(X)\right|\right)^2\\
&\leq \mathbb{E}\left[\left(\frac{1}{n} \sum_{i=1}^n g(X_i)- \mathbb{E}g(X)\right)^2\right]\\
&= \frac{1}{n} \var \big(g(X)\big)\\
&\leq \frac{1}{4n}\left(\sup_x h(P_{T|X=x}) - \inf_x {h}(P_{T|X=x})\right)^2.\numberthis\label{EQ:THIRDTERM}
\end{align*}
The last inequality follows since $\var(A)\leq \frac{1}{4}(\sup A - \inf A)^2$ for any random variable $A$.

It remains to bound the supremum and infimum of $h(P_{T|X=x})$ uniformly in $x\in\mathbb{R}^{d_0}$. By definition $T = S + Z$, where $S$ and $Z$ are independent and $Z \sim \mathcal{N}(0,\sigma^2 \mathrm{I}_d)$. Therefore, for all $x\in\mathbb{R}^{d_0}$
\begin{align*}
h(P_{T|X=x})&=h(S+Z|X=x)\\&\geq h(S+Z|S,X=x)\\&=h(Z)\\&=\frac{d}{2}\log( 2\pi e \sigma^2),\numberthis\label{EQ:INFX}
\end{align*}
where we have used the independence of $Z$ and $(S,X)$ and the fact that conditioning cannot increase entropy. On the other hand, denoting the entries of $T$ by $T\triangleq\big(T(k)\big)_{k=1}^d$, we can obtain an upper bound as
\begin{equation}
    h(P_{T|X=x})=h(T|X=x)\leq \sum_{k=1}^d h\big(T(k)\big|X=x\big),
\end{equation}
since independent random variables maximize differential entropy. Now for any $k\in[d]$, we have
\begin{equation}
    \var\big(T(k)\big|X=x\big)\leq\mathbb{E}\big[T^2(k)\big|X=x\big]\leq 1+\sigma^2,
\end{equation}
because $S(k)\in[-1,1]$ almost surely. Since the Gaussian distribution maximizes differential entropy under a variance constraint, we have
\begin{equation}
h(P_{T|X =x}) \leq \frac{d}{2}\log\big(2\pi e (1+\sigma^2)\big).\label{EQ:SUPEX}
\end{equation}
for all $x\in\mathbb{R}^{d_0}$. Substituting the lower bound \eqref{EQ:INFX} and upper bound \eqref{EQ:SUPEX} into \eqref{EQ:THIRDTERM} gives
\begin{equation}
(\mathrm{III})^2 \leq \left(\frac{d\log\left(1+\frac{1}{\sigma^2}\right)}{4\sqrt{n}}\right)^2. \end{equation}
Inserting this along with \eqref{EQ:FIRSTTERM} and  \eqref{EQ:SECONDTERM} into the bound \eqref{eq:BarError} bounds the expected estimation error as
\begin{equation}
    \mathbb{E}\left|\hat{I}_{\mathsf{Input}}\p{X^n,\hat{h},\sigma} - I(X;T)\right|\leq 2\Delta_n+\frac{d\log\left(1+\frac{1}{\sigma^2}\right)}{4\sqrt{n}}.
\end{equation}
Taking the supremum over $P_X$ concludes the proof.


\section{Summary and Concluding Remarks}\label{SEC:summary}

This work first explored the problem of empirical approximation under Gaussian smoothing in high dimensions. To quantify the approximation error, we considered various statistical distances, such as 1-Wasserstein, squared 2-Wasserstein, TV, KL divergence and $\Chi$-divergence. It was shown that when $P$ is subgaussian, the 1-Wasserstein and the TV distances converge as $n^{-1/2}$. The parametric convergence rate is also attained by the KL divergence, squared 2-Wasserstein distance and $\Chi$-divergence, so long that the $\Chi$ mutual information $I_{\Chi}(S;Y)$, for $Y=S+Z$ with $S\sim P$ independent of $Z\sim\Gauss$, is finite. The latter condition is always satisfied by $K$-subgaussian $P$ distributions in the low SNR regime where $K<\frac{\sigma}{2}$. However, when SNR is high ($K>\sqrt{2}\sigma$), there exist $K$-subgaussian distributions $P$ for which $I_{\Chi}(S;Y)=\E_{P^{\otimes n}}\Chi\left(\hat{P}_{S^n}\ast\Gauss\middle\|P\ast\Gauss\right)=\infty$. Whenever this happens, it was further established that the KL divergence and the squared 2-Wasserstein distance are $\omega(n^{-1})$. Whenever the parametric convergence rate of the smooth empirical measure is attained, it strikingly contrasts classical (unconvolved) results, e.g., for the Wasserstein distance, which suffer from the curse of dimensionality (see, e.g., \cite[Theorem 1]{FournierGuillin2015}). 

The empirical approximation results were used to study differential entropy estimation under Gaussian smoothing. Specifically, we considered the estimation of $\mathsf{T}_{\sigma}(P)=h(P\ast\Gauss)$ based on i.i.d. samples from $P$ and knowledge of the noise distribution $\Gauss$. It was shown that the absolute-error risk of the plug-in estimator over the bounded support and subgaussian classes converges as $e^{O(d)}n^{-1/2}$ (with the prefactor explicitly characterized). This established the plug-in estimator as minimax-rate optimal. The exponential dependence of the sample complexity on dimension was shown to be necessary. These results were followed by a bias lower bound of order $\log\left(2^d n^{-1}\right)$, as well as an efficient and provably accurate MC integration method for computing the plug-in estimator. 

The considered differential entropy estimation framework enables studying information flows in DNNs \cite{ICML_Info_flow2019}. In Section \ref{SEC:applications} we showed how the mutual information between layers of a DNN reduces to estimating $h(P\ast\Gauss)$. An ad hoc estimator for $h(P\ast\Gauss)$ was important here because the general-purpose estimators (based on noisy samples from $P\ast\Gauss$) available in the literature are unsatisfactory for several (theoretical and/or practical) reasons. Most theoretical performance guarantees for such estimators are not valid in our setup, as they typically assume that the unknown density is positively lower bounded inside its compact support. 

To the best of our knowledge, the only two works that provide convergence results that apply here are \cite{han2017optimal} and \cite{berrett2019efficient}. The rate derived for the KDE-based estimator in \cite{han2017optimal}, however, effectively scales as $n^{-1/d}$ for large dimensions, which is too slow for practical purposes. Remarkably, \cite{berrett2019efficient} proposes a wKL estimator in the very smooth density regime that provably attains the parametric rate of estimation in our problem (e.g., when $P$ is compactly supported). This result, however, does not characterize the dependence of that rate on $d$. Understanding this dependence is crucial in practice. Indeed, in Section \ref{SEC:simulations} we show that, empirically, the performance of the wKL significantly deteriorates as $d$ grows. In all our experiments, the plug-in estimator outperforms the wKL method from \cite{berrett2019efficient} (as well as all other generic estimator we have tested), converging faster with $n$ and scaling better with $d$. 

For future work, open questions regarding the smooth empirical measure convergence were listed in Section \ref{SUBSEC:open_questions}. On top of that, there are appealing extensions of the differential estimation question to be considered. This includes non-Gaussian additive noise models or multiplicative Bernoulli noise (which corresponds to DNNs with dropout regularization). The question of estimating $h(P\ast\Gauss)$ when only samples from $P\ast\Gauss$ are available (yet the Gaussian convolution structure is known) is also attractive. This would, however, require a different technique to that employed herein. Our current method strongly relies on having `clean' samples from $P$. Beyond this work, we see considerable virtue in exploring additional ad hoc estimation setups with exploitable structure that might enable improved estimation results.

\section*{Acknowledgement}
The authors thank Yihong Wu for helpful discussions.


\appendices


\section{Proof of Proposition \ref{PROP:minimaxRate}}\label{app:rateVerification}


We lower bound the minimax risk in the nonparametric estimation risk by a reduction to a parametric setup. Without loss of generality, assume $d=1$ (the risk of the one-dimensional estimation problem trivially lower bounds that of its $d$-dimensional counterpart). Recall that $\FSGO$ is the class of $K$-subgaussian measures on $\RR$. Define $\mathcal{G}_K:=\{\cN_\nu\}_{\nu\in[K/2,K]}$ as the collection of all centered Gaussian measures, each with variance $\nu^2$. Noting that $\cG\subset\FSGO$, we obtain
\begin{equation}
    \mathcal{R}^\star\left(n,\sigma,\FSGO\right)\geq \inf\limits_{\hat{h}}\sup\limits_{P\in\cG}\mathbb{E}\left|h(P\ast\Gauss)-\hat{h}(S^n,\sigma)\right|,\label{EQ:risk_LB}
\end{equation}
and henceforth focus on lower bounding the RHS.


Note that $h(P\ast\Gauss)=\frac{1}{2}\log(\nu^2 + \sigma^2) + \frac{1}{2}\log(2\pi e)$, for $P=\cN_\nu\in\cG$. Thus, the estimation of $h(P\ast\Gauss)$, when $P\in\cG$, reduces to estimating $\frac{1}{2}\log (\nu^2 + \sigma^2)$ from samples of $\cN_\nu$. Recall that $\nu \in [K/2,K]$ is considered unknown and $\sigma$ known. This simple setting lands within the framework studied in \cite{chen1997general}, where lower bounds on the minimax absolute error in terms of the associated Hellinger modulus were derived. We follow the proof style of Corollary 3 therein.



Firstly, recall that the squared Hellinger distance between $\mathcal{N}(0,\nu_1^2)$ and $\mathcal{N}(0,\nu_0^2)$ is given by
\begin{equation}
\rho^2(\nu_1, \nu_0) = 2 \left(1 - \sqrt{\frac{2 \nu_1 \nu_0}{\nu_1^2 + \nu_0^2}}\right). \label{EQ:Hellinger_Gauss}
\end{equation}
Since we are estimating $\frac{1}{2}\log (\nu^2 + \sigma^2)$, for convenience we denote the parameter of interest as $\theta(\nu) := \frac{1}{2}\log (\nu^2 + \sigma^2)$.
The Hellinger modulus of $\theta(\nu_0)$ is defined as
\[
\omega_q\left(\sqrt{\frac{\alpha}{n}}; \theta(\nu_0)\right) = \sup_{\nu:\ \rho(\nu,\nu_0) \leq \sqrt{\frac{\alpha}{n}}}\left|\theta(\nu) - \theta(\nu_0) \right|.\numberthis \label{eq:hellinger}
\]
Based on Theorem 2 of \cite{chen1997general}, if $\omega_q\left(\sqrt{\frac{\alpha}{n}}; \theta(\nu_0)\right)=\Omega(n^{-1/2})$ then the RHS of \eqref{EQ:risk_LB} is also $\Omega(n^{-1/2})$. We thus seek to bound this modulus.

We start by characterizing the set of $\nu$ values that lie in the Hellinger ball of radius $\sqrt{\frac{\alpha}{n}}$ around~$\nu_0$. From \eqref{EQ:Hellinger_Gauss} and by defining $\xi = \left(1 -\frac{\alpha}{2n}\right)^{-2}$, one readily verifies that $\rho(\nu,\nu_0)\leq \sqrt{\frac{\alpha}{n}}$ if and only if
\[
\xi - \sqrt{\xi^2 - 1}\leq \frac{\nu}{\nu_0} \leq \xi + \sqrt{\xi^2 - 1}.
\]
Equivalently, the feasible set of $\nu$ values satisfies 
\begin{align*}
&\left|\log ({\nu}) - \log (\nu_0 )\right| \\
&\qquad\qquad\in \left[\log\left(\xi - \sqrt{\xi^2 - 1}\right),\log\left(\xi + \sqrt{\xi^2 - 1}\right) \right].
\end{align*}


One may check that $2\sqrt{\frac{\alpha}{2n}}$, for all positive $\alpha,n>0$ with $\alpha/n < 1$, belongs to the above interval. Hence, $\nu_\star$ such that $\log(\nu_\star)  = \log(\nu_0) + 2\sqrt{\frac{a}{2n}}$ is feasible and we may substitute it into \eqref{eq:hellinger} to lower bound the modulus as follows:
\begin{align*}
\omega_q\left(\sqrt{\frac{a}{n}}; \theta(\nu_0)\right)&\geq|\theta(\nu_\star) - \theta(\nu_0)|\\
&=\frac{1}{2}\log (\nu^2_\star + \sigma^2) - \frac{1}{2}\log (\nu^2_0 + \sigma^2) \\
&\stackrel{(a)}\geq \frac{1}{1 + \frac{\sigma^2}{\nu_0^2}} \log\left( \frac{\nu_\star}{\nu_0}\right)\\
&\stackrel{(b)}\geq \frac{1}{1 + \frac{4\sigma^2}{K^2}} \log\left( \frac{\nu_\star}{\nu_0}\right)\\
&= 2\frac{K^2}{K^2 + {4\sigma^2}}\sqrt{\frac{\alpha}{2n}},
\end{align*}
where (a) is because for any $a,b,c > 0$ with $b \geq a$ we have
\begin{align*}
\log (b + c) - \log (a + c)
&= \int_a^b \frac{1}{x+c} dx\\
&\geq \int_a^b \frac{1}{(1+c/a) x} dx\\
&= \frac{1}{1+\frac{c}{a}} \log \frac{b}{a},
\end{align*}
and (b) follows since $\nu_0 \geq K/2$.

Applying Theorem 2 of \cite{chen1997general} to this bound on the Hellinger modulus implies that the best estimator of $h(P\ast\Gauss)$ over the class $\cG_K$ achieves $\Omega(1/\sqrt{n})$ in absolute error.

\section{Label and Hidden Layer Mutual Information}\label{APPEN:label_MI}

Consider the estimation of $I(Y;T)$, where $Y$ is the true label and $T$ is a hidden layer in a noisy DNN. For completeness, we first describe the setup (repeating some parts of Remark \ref{REM:MI_label}). Afterwards, the proposed estimator for $I(Y;T)$ is presented and an upper bound on the estimation error is stated and proven. 

Let $(X,Y)\sim P_{X,Y}$ be a feature-label pair, whose distribution is unknown. Assume that $\mathcal{Y}\triangleq\supp(P_Y)$ is finite and known (as is the case in practice) and let $|\mathcal{Y}|=K$ be the cardinality of $\mathcal{Y}$. Let $\big\{(X_i,Y_i)\big\}_{i=1}^n$ be a set of $n$ i.i.d. samples from $P_{X,Y}$, and $T$ be a hidden layer in a noisy DNN with input $X$. Recall that $T=S+Z$, where $S$ is a deterministic map of the previous layer and $Z\sim\mathcal{N}(0,\sigma^2\mathrm{I}_d)$. The tuple $(X,Y,S,T)$ is jointly distributed according to $P_{X,Y}P_{S|X}P_{T|S}$, under which $Y-X-S-T$ forms a Markov chain. Our goal is to estimate the mutual information
\begin{equation}
    I(Y;T)=h(P_S\ast\Gauss)-\sum_{y\in\mathcal{Y}}p_Y(y)h(P_{S|Y=y}\ast\Gauss),\label{EQ:MI_label}
\end{equation}
based on a given estimator $\hat{h}$ of $h(P\ast\Gauss)$ that knows $\sigma$ and uses i.i.d. samples from $P\in\mathcal{F}_d$. In \eqref{EQ:MI_label}, $p_Y$ is the PMF associated with $P_Y$.

We first describe the sampling procedure for estimating each of the differential entropies from \eqref{EQ:MI_label}. For the unconditional entropy, $P_S$ is sampled in the same manner described in Section \ref{SUBSEC:MI_Noisy_DNNs_estimation} for the estimation of $I(X;T)$. Denote the obtained samples by $S^n$. To sample from $P_{S|Y=y}$, for a fixed label $y\in\mathcal{Y}$, fix a sample set  $\big\{(x_i,y_i)\big\}_{i=1}^n$ and consider the following. Define the set $\mathcal{I}_y\triangleq\big\{i\in[n]\big|y_i=y\big\}$ and let $\mathcal{X}_y\triangleq\{x_i\}_{i\in\mathcal{I}_y}$ be the subset of features whose label is $y$; the elements of $\mathcal{X}_y$ are conditionally i.i.d. samples from $P_{X|Y=y}$. Now, feed each $x\in\mathcal{X}_y$ into the noisy DNN and collect the values induced at the layer preceding $T$. By applying the appropriate deterministic function on each of these samples we get a set of $n_y\triangleq |\mathcal{I}_y|$ i.i.d. samples from $P_{S|Y=y}$. Denote this sample set by $S^{n_y}(\mathcal{X}_y)$.

Similarly to Section \ref{SUBSEC:MI_Noisy_DNNs_estimation}, suppose we are given an estimator $\hat{h}(A^m,\sigma)$ of $h(P\ast\Gauss)$, for $P\in\mathcal{F}_d$, based on $m$ i.i.d. samples $A^m=\{A_1,\ldots,A_m\}$ from $P$. Assume that $\hat{h}$ attains
\begin{equation}
    \sup_{P\in\mathcal{F}_d} \mathbb{E}_{P^{\otimes m}}\left|h(P\ast\Gauss)-\hat{h}(A^m,\sigma)\right|\leq \Delta_{\sigma,d}(m).\label{EQ:given_estimator}
\end{equation}
Further assume that $\Delta_{\sigma,d}(m)<\infty$, for all $m\in\mathbb{N}$, and that $\lim_{m\to\infty}\Delta_{\sigma,d}(m)=0$, for any fixed $\sigma$ and $d$ (otherwise, $\hat{h}$ is not a good estimator and there is no hope using it for estimating $I(Y;T)$). Without loss of generality we may also assume that $\Delta_{\sigma,d}(m)$ is monotonically decreasing in $m$. Our estimator of $I(Y;T)$~is
\begin{equation}
\hat{I}_\mathsf{Label}\p{X^n,Y^n,\hat{h},\sigma}\triangleq\hat{h}(S^n,\sigma)-\sum_{y\in\mathcal{Y}} \hat{p}_{Y^n}(y)\hat{h}\big(S^{n_y}(\mathcal{X}_y),\sigma\big),\label{EQ:MI_label_est}
\end{equation}
where $\hat{p}_{Y^n}(y)\triangleq \frac{1}{n}\sum_{i=1}^n\mathds{1}_{\{Y_i=y\}}$ is the empirical PMF associated with the labels $Y^n$. The following proposition bounds the expected absolute-error risk of $\hat{I}_\mathsf{Label}\p{X^n,Y^n,\hat{h},\sigma}$; the proof is given after the statement.

\begin{proposition}[Label-Hidden Layer Mutual Information]\label{PROP:MI_label}
For the above described estimation setting, we have
\begin{align*}
    &\sup_{P_{X,Y}:\ |\mathcal{Y}|=K}\mathbb{E} \left|I(Y;T)-\hat{I}_\mathsf{Label}\p{\mspace{-2mu}X^n,Y^n,\hat{h},\sigma\mspace{-2mu}}\right|\leq\\& \Delta_{\sigma,d}(n)+c^{\mathsf{(MI)}}_{\sigma,d}\!\sqrt{\frac{K\!-\!1}{n}}+K\!\p{\!\Delta^\star_{\sigma,d}\cdot e^{-\frac{np_l^2}{8p_u}}+\Delta_{\sigma,d}\p{\frac{np_l}{2}}\!}\mspace{-2mu}, 
\end{align*}
where 
\begin{subequations}
\begin{align}
    c^{\mathsf{(MI)}}_{\sigma,d}&\triangleq \frac{d}{2}\max\Big\{-\log(2\pi e\sigma^2),\log\big(2\pi e(1+\sigma^2)\big)\Big\}\\
    p_l&\triangleq \min_{y\in\mathcal{Y}}p_Y(y)\\
    p_u&\triangleq \max_{y\in\mathcal{Y}}p_Y(y)\\
    \Delta^\star_{\sigma,d}&\triangleq \max_{n\in\mathbb{N}}\Delta_{\sigma,d}(n).
\end{align}
\end{subequations}
\end{proposition}
The proof is reminiscent of that of Proposition \ref{PROP:MI_True_Data_Dist}, but with a few technical modifications accounting for $n_y$ being a random quantity (as it depends on the number of $Y_i$-s that equal to $y$). To control $n_y$ we use the concentration of the Binomial distribution about its mean. 

\begin{IEEEproof} Fix $P_{X,Y}$ with $|\mathcal{Y}|=K$, and use the triangle inequality to get
\begin{align*}
&\mathbb{E}\Big|I(Y;T)-\hat{I}_\mathsf{Label}\p{X^n,Y^n,\hat{h},\sigma}\Big|\\
&\leq \underbrace{\mathbb{E}\left|h(P_S\ast\Gauss)-\hat{h}(S^n,\sigma) \right|}_{(\mathrm{I})}\\
&\hskip 3mm+\underbrace{\sum_{y\in\mathcal{Y}}\big|h(P_{S|Y=y}\ast\Gauss)\big|\mathbb{E}\big|p_Y(y)-\hat{p}_{Y^n}(y)\big|}_{(\mathrm{II})}\\
&\hskip 6mm+\underbrace{\sum_{y\in\mathcal{Y}}\mathbb{E}\bigg|\hat{p}_{Y^n}(y)\Big(h(P_{S|Y=y}\ast\Gauss)-\hat{h}\big(S^{n_y}(\mathcal{X}_y),\sigma^2\big)\Big)\bigg|}_{(\mathrm{III})},\numberthis\label{EQ:label_MI_error_terms}
\end{align*}
where we have added and subtracted $\sum_{y\in\mathinner{Y}}\hat{p}_{Y^n}(y)h(P_{S|Y=y}\ast\Gauss)$ inside the original expectation.

Clearly, (I) is bounded by $\Delta_{\sigma,d}(n)$. For (II), we first bound the conditional differential entropies. For any $y\in\mathcal{Y}$, we have
\begin{align*}
    h(P_{S|Y=y}\ast\Gauss)&=h(S+Z|Y=y)\geq h(S+Z|S,Y=y)\\&=\frac{d}{2}\log(2\pi e \sigma^2),\numberthis\label{EQ:label_ent_bound1}
\end{align*}
where the last equality is since $(Y,S)$ is independent of $Z\sim\mathcal{N}(0,\sigma^2\mathrm{I}_d)$. Furthermore, 
\begin{align*}
    h(P_{S|Y=y}\ast\Gauss)&\leq \sum_{k=1}^d h\big(S(k)+Z(k)\big|Y=y\big)\\&\leq \frac{d}{2}\log\big(2\pi e (1+\sigma^2)\big),\numberthis\label{EQ:label_ent_bound2}
\end{align*}
where the first inequality is because independence maximizes differential entropy, while the second inequality uses $\var\big(S(k)+Z(k)\big|Y=y\big)\leq 1+\sigma^2$. Combining \eqref{EQ:label_ent_bound1} and \eqref{EQ:label_ent_bound2} we obtain
\begin{align*}
    \big|h(&P_{S|Y=y}\ast\Gauss)\big|\\&\leq c^{\mathsf{(MI)}}_{\sigma,d}\triangleq\frac{d}{2}\max\Big\{-\log(2\pi e\sigma^2),\log\big(2\pi e(1+\sigma^2)\big)\Big\}.\numberthis\label{EQ:label_termII_const}
\end{align*}
For the expected value in (II), monotonicity of moment gives
\begin{align*}
    \mathbb{E}\big|p_Y(y)-\hat{p}_{Y^n}(y)\big|&\leq \sqrt{\var\big(p_{Y^n}(y)\big)}\\&=\sqrt{\frac{1}{n}\var\p{\mathds{1}_{\{Y=y\}}}}\\&=\sqrt{\frac{p_Y(y)\big(1-p_Y(y)\big)}{n}}.\numberthis\label{EQ:label_termII_expectation}
\end{align*}
Using \eqref{EQ:label_termII_const} and \eqref{EQ:label_termII_expectation} we bound Term (II) as follows:
\begin{equation}
    (\mathrm{II})\leq \frac{c^{\mathsf{(MI)}}_{\sigma,d}}{\sqrt{n}}\sum_{y\in\mathcal{Y}}\sqrt{p_Y(y)\big(1-p_Y(y)\big)}\leq c^{\mathsf{(MI)}}_{\sigma,d}\sqrt{\frac{K-1}{n}},\label{EQ:label_termII}
\end{equation}
where the last step uses the Cauchy-Schwarz inequality. 

For Term (III), we first upper bound $\hat{p}_{Y^n}(y)\leq 1$, for all $y\in\mathcal{Y}$, which leaves us to deal with the sum of expected absolute errors in estimating the conditional entropies. Fix $y\in\mathcal{Y}$, and notice that $n_y\sim\mathsf{Binom}\big(p_Y(y),n\big)$. Define $p_l\triangleq\min_{y\in\mathcal{Y}}p_Y(y)$ and $p_u\triangleq\max_{y\in\mathcal{Y}}p_Y(y)$ as in the statement of Proposition \ref{PROP:MI_label}. Using a Chernoff bound for the Binomial distribution we have that for any $k\leq np_Y(y)$,
\begin{align*}
    \mathbb{P}\Big(n_y\leq k\Big)&\leq \exp\p{-\frac{1}{2p_Y(y)}\cdot\frac{\big(np_Y(y)-k\big)^2}{n}}\\&\leq \exp\p{-\frac{1}{2p_u}\cdot\frac{\big(np_Y(y)-k\big)^2}{n}}.
\end{align*}
Set $k^\star_y=n\p{p_Y(y)-\frac{1}{2}p_l}\in\big(0,np_Y(y)\big)$ into the above to get
\begin{equation}
    \mathbb{P}\Big(n_y\leq k^\star_y\Big)\leq \exp\p{-\frac{np_l^2}{8p_u}}.\label{EQ:label_termIII_Chernoff}
\end{equation}

Setting $\Delta^\star_{\sigma,d}\triangleq\max_{n\in\mathbb{N}}\Delta_{\sigma,d}(n)$, we note that $\Delta^\star_{\sigma,d}<\infty$ by hypothesis, and bound (III) as follows:
\begin{align*}
    &(\mathrm{III})\leq\sum_{y\in\mathcal{Y}}\mathbb{E}\Big|h(P_{S|Y=y}\ast\Gauss)-\hat{h}\big(S^{n_y}(\mathcal{X}_y),\sigma^2\big)\Big|\\
    &\stackrel{(a)}=\!\!\sum_{y\in\mathcal{Y}}\!\mathbb{E}_{n_y}\!\Bigg[\!\mathbb{E}\bigg[\!\Big|\hat{p}_{Y^n}(y)\Big(h(P_{S|Y\!=y}\!\ast\!\Gauss)\!-\!\hat{h}\big(S^{n_y}\!(\mathcal{X}_y),\!\sigma^2\big)\Big|\bigg|n_y\bigg]\!\Bigg]\\
    &\stackrel{(b)}=\sum_{y\in\mathcal{Y}}\mathbb{E}_{n_y}\Delta_{\sigma,d}(n_y)\\
    &\stackrel{(c)}=\sum_{y\in\mathcal{Y}}\mathbb{P}\Big(n_y\leq k^\star_y\Big)\mathbb{E}\big[\Delta_{\sigma,d}(n_y)\big|n_y\leq k^\star_y \big]\\&\qquad\qquad+\mathbb{P}\Big(n_y>k^\star_y\Big)\mathbb{E}\big[\Delta_{\sigma,d}(n_y)\big|n_y>k^\star_y \big]\\
    &\stackrel{(d)}\leq K\p{\Delta^\star_{\sigma,d}\cdot e^{-\frac{np_l^2}{8p_u}}+\Delta_{\sigma,d}\p{\frac{np_l}{2}}}, \numberthis\label{EQ:label_termIII}
\end{align*}
where (a) and (c) use the law of total expectation, (b) is since for each fixed $n_y=k$, the expected differential entropy estimation error (inner expectation) is bounded by $\Delta_{\sigma,d}(k)$, while (d) relies on \eqref{EQ:label_termIII_Chernoff}, the definition of $\Delta_{\sigma,d}^\star$ and the fact that $\Delta_{\sigma,d}(n)$ is monotonically decreasing with $n$ along with $k^\star_y\geq \frac{np_l}{2}$, for all $y\in\mathcal{Y}$. Inserting (I) $\leq \Delta_{\sigma,d}(n)$ together with the bounds from \eqref{EQ:label_termII} and \eqref{EQ:label_termIII} back into \eqref{EQ:label_MI_error_terms} and taking the supremum over all $P_{X,Y}$ with $|\mathcal{Y}|=K$ concludes the proof.

\end{IEEEproof}

\section{Proof of Lemma \ref{LEMMA:SP_bias_MI}}\label{APPEN:SP_bias_MI_proof}

We expand $I(S^n;Y)=h(Y)-h(Y|S^n)$. Let $T=S+Z\sim P\ast\Gauss$ and first note that for any measurable set~$\mathcal{A}$,
\begin{align*}
    \mathbb{P}(Y \in \mathcal{A})&=\mathbb{P}\big(S_W+Z \in \mathcal{A} \big)\\&=\frac{1}{n}\sum_{i=1}^n\mathbb{P}(S_i+Z \in  \mathcal{A})= \mathbb{P}(T \in \mathcal{A}).
\end{align*}
Thus, $h(Y)=h(P\ast\Gauss)$. It remains to show that $h(Y|S^n)=\mathbb{E}_{P^{\otimes n}}h(\hat{P}_{S^n}\ast\Gauss)$. Fix $S^n=s^n$ and consider
\begin{align*}
    \mathbb{P}\big(Y \in \mathcal{A} \big|S^n=s^n\big)&=\mathbb{P}\big(S_W+Z \in \mathcal{A}\big|S^n=s^n\big)\\&=\frac{1}{n}\mathbb{P}\big(s_i+Z \in \mathcal{A} \big),
\end{align*}
which implies that the density $p_{Y|S^n=s^n}$ equals the density of $\hat{P}_{s^n}\ast\Gauss$. Consequently, $h(Y|S^n=s^n)=h(\hat{P}_{s^n}\ast\Gauss)$, and by definition of conditional entropy $h(Y|S^n)=\mathbb{E}_{P^{\otimes n}}h(\hat{P}_{S^n}\ast\Gauss)$.

\section{Proof of Proposition \ref{prop:div}}
\label{app:counterex}

We start from the derivation of \eqref{EQ:Chi_empirical_to_MI}, which shows that
\begin{align*}
    \E_{P^{\otimes n}} \Chi &\left(\hat{P}_{S^n} \ast\Gauss\middle\|P\ast\Gauss\right)= \frac{1}{n}I_{\Chi}(S;Y)\\&=\frac{1}{n}\left( \int_{\RR^d}\frac{\E_P\gauss^2(z - S)}{q(z)}\dd z-1\right)
\end{align*}
for $S\sim P$ and $Y=S+Z$, where $Z\sim \Gauss$ is independent of $S$. 
Recalling that without loss of generality $\sigma =1$ and that $q(z) = \mathbb{E}_{P} \gaussI(z - S)$, a sufficient condition for divergence in Proposition \ref{prop:div} is
 \begin{equation}\label{eq:Chi2Int}
 \int_{\RR} \frac{\E_P \gausssI(z - S)}{\E_P \gaussI(z - S)}\dd z=\infty .
 \end{equation}

Under the $P$ from \eqref{eq:Pcounter}, the left-hand side (LHS) of \eqref{eq:Chi2Int} becomes
\begin{align*}
    &\int_{\mathbb{R}} \frac{\sum_{k = 0}^\infty p_k \gausssI(z - r_k)}{\sum_{k = 0}^\infty p_k \gaussI(z - r_k)}\dd z 
    \\&= \sum_{k = 0}^\infty \int_{\mathbb{R}} \frac{ p_k \gausssI(z - r_k)}{ p_k \gaussI(z - r_k) + \sum_{j \neq k} p_j \gaussI(z - r_j)}\dd z \\
    &= \sum_{k = 0}^\infty \int_{\mathbb{R}} \frac{ \gausssI(z - r_k)}{  \gaussI(z - r_k)} \frac{1}{1 + \sum_{j \neq k} \frac{p_j}{p_k} \frac{\gaussI(z - r_j)}{\gaussI(z - r_k)}}\dd z\\
    &\geq \sum_{k = 1}^\infty \int_{r_k - \frac{1}{100}}^{r_k + \frac{1}{100}} \Bigg[\frac{ \gausssI(z - r_k)}{  \gaussI(z - r_k)} \\&\qquad\qquad\times\left.\frac{1}{1 + \sum\limits_{j =0}^{k-1} \frac{p_j}{p_k} \frac{\gaussI(z - r_j)}{\gaussI(z - r_k)}+ \sum\limits_{j =k + 1}^\infty \frac{p_j}{p_k} \frac{\gaussI(z - r_j)}{\gaussI(z - r_k)}}\right]\dd z,\numberthis \label{eq:ratiobd}
\end{align*}
where the inequality follows since the integrands are all nonnegative and the domain of integration has been reduced.

We now bound the sums in the second denominator of \eqref{eq:ratiobd} for $k > 0$ and $z \in \left\{r_k-\frac{1}{100},r_k+\frac{1}{100}\right\}$ (as indicated by the support of the outer sum and integral).
First, consider the ratio $\frac{p_j}{p_k} \frac{\gaussI(z - r_j)}{\gaussI(z - r_k)}$. For $j = 0$ and $\epsilon \leq \frac{1}{4}$, we have
\begin{align*}
    \frac{p_0}{p_k} \frac{\gaussI(z )}{\gaussI(z - r_k)}&\leq C\cdot\exp\left\{{\epsilon r_k^2} +{\frac{r_k^2-2\left(r_k -\frac{1}{100}\right) r_k}{2}}\right\}\\&= C\cdot \exp\left\{\left(\epsilon - \frac{1}{2}\right) r_k^2 + \frac{r_k}{100}\right\}\leq C,\numberthis\label{EQ:demon_bound1}
\end{align*}
where $C$ is a constant depending on $K$ only, and in the last inequality is because $r_k \geq 1$, for $k \geq 1$. For $j > 0$, denoting $\alpha_\epsilon\triangleq1-\sqrt{2\epsilon}$, the bound becomes
\begin{align*}
    &\frac{p_j}{p_k} \frac{\gaussI(z - r_j)}{\gaussI(z - r_k)} 
    = \exp\left\{\epsilon (r_k^2-r_j^2)+\frac{r_k^2 - r_j^2-2z(r_k - r_j)}{2}\right\}\\
    &= \exp\Bigg\{{\epsilon r_k^2\left(1 - \alpha_\epsilon^{-2(j-k)}\right)} \\&\hskip 15mm +\left.{\frac{r_k^2\left(1 - \alpha_\epsilon^{-2(j-k)}\right)-2z r_k\left(1 - \alpha_\epsilon^{-(j-k)}\right)}{2}}\right\}.\numberthis\label{EQ:demon_bound2}
\end{align*}

Using \eqref{EQ:demon_bound1} and \eqref{EQ:demon_bound2}, for $\epsilon=\frac{1}{4}$ and $\alpha\triangleq\alpha_{\frac{1}{4}}=1-\frac{1}{\sqrt{2}}$, we have
\begin{align*}
    &\sum_{j =0}^{k-1} \frac{p_j}{p_k} \frac{\gaussI(z - r_j)}{\gaussI(z - r_k)}\\
    &\leq C + \sum_{j =1}^{k-1} \exp\Bigg\{{\epsilon r_k^2\left(1 - \alpha^{-2(j-k)}\right)} +\frac{r_k^2}{2}\left(1 - \alpha^{-2(j-k)}\right)\\&\hskip 35mm \left.-r_k\left(r_k - \frac{1}{100}\right) \left(1 - \alpha^{-(j-k)}\right)\right\}\\
    &\leq C + \sum_{j =1}^{k-1} \exp\left\{{\epsilon r_k^2\left(1 - \alpha^{-2(j-k)}\right)} +\frac{r_k^2}{2}\left(1 - \alpha^{-2(j-k)}\right)\right.\\
    &\left.-r_k\left(r_k \!-\! \frac{1}{100}\right)\! \left(1 - \alpha^{-2(j-k)}\right)\min_{j}\left(\left(1 + \alpha^{-(j-k)}\right)^{-1}\!\right)\!\right\}\\
    &= C + \sum_{j =1}^{k-1} \exp\left\{{\left(\epsilon + \frac 1 2 -\frac{1}{2 - \sqrt{1/2}}\right) r_k^2\left(1 - \alpha^{-2(j-k)}\right)}\right.\\&\hskip 52 mm +\frac{r_k}{100}\left(1 - \alpha^{-2(j-k)}\right)\Bigg\}\\
    &\leq C + k - 1,\numberthis
\end{align*}
where the last inequality follows since $\epsilon = \frac{1}{4}$, $r_k^2 \geq r_k$, for $k \geq 1$, and $\frac{1}{4} + \frac 1 2 - \frac{1}{2 - \sqrt{1/2}} + \frac{1}{100} < 0$. 

Proceeding onto the series for $j\geq k+1$, we have
\begin{align*}
    &\sum_{j =k+1}^{\infty} \frac{p_j}{p_k} \frac{\gaussI(z - r_j)}{\gaussI(z - r_k)} 
    \\&\stackrel{(a)}\leq \sum_{j =k+1}^{\infty} \exp\left\{{\frac{3}{4} r_k^2\left(1 - \alpha^{-2(j-k)}\right)}\right. \\&\hskip 35mm -r_k\left(r_k + \frac{1}{100}\right) \left(1 - \alpha^{-(j-k)}\right)\Bigg\}\\
    &\leq \sum_{j =k+1}^{\infty} \exp\left\{{-\frac{3}{4} r_k^2 \alpha^{-2(j-k)}} +r_k\left(r_k + \frac{1}{100}\right) \alpha^{-(j-k)}\right\}\\
    &\leq \sum_{j=k+1}^\infty \exp\left\{-r_k^2\alpha^{-(j-k)}  \left( \frac{3 \alpha^{-(j-k)}}{4}- \frac{101}{100}\right) \right\}\\
    &\stackrel{(b)}\leq \sum_{j=k+1}^\infty \exp\left\{-\frac{1}{4}r_k^2\alpha^{-(j-k)}\right\}\\
   &\leq \sum_{\ell=1}^\infty \exp\left\{-\frac{1}{4}\alpha^{-\ell}\right\}\\
   &\stackrel{(c)}\leq \frac{1}{4},
\end{align*}
where (a) uses \eqref{EQ:demon_bound2} and the fact that $\left(1 - \alpha^{-(j-k)}\right)$ is negative for $j > k$, (b) is since $\frac{3\alpha^{-t}}{4}- \frac{101}{100} \geq 1/4$ for all $t \geq 1$, and (c) follows by numerical computation and because the series converges by the ratio test.

Substituting these bounds into the LHS of \eqref{eq:Chi2Int}, we get 
\begin{align*}
    \int_{\mathbb{R}} &\frac{\sum_{k = 0}^\infty p_k \gausssI(z - r_k)}{\sum_{k = 0}^\infty p_k \gaussI(z - r_k)}\dd z \\&\geq \sum_{k = 1}^\infty \int_{r_k -  \frac{1}{100}}^{r_k +  \frac{1}{100}} \frac{ \gausssI(z - r_k)}{  \gaussI(z - r_k)} \frac{1}{1 + C + k - 1 +  \frac{1}{4}}\dd z\\
    &= \left(\int_{ - \frac{1}{4}}^{ \frac{1}{4}} \frac{ \gausssI(z)}{  \gaussI(z)}\dd z\right)\sum_{k = 1}^\infty \frac{1}{k + C + \frac{1}{4}}.\numberthis
\end{align*}
The RHS above diverges because the integral is nonzero and $\sum_{k=1}^\infty \frac{1}{k + C + 1/4}$ is a harmonic series.

\bibliographystyle{IEEEtran}
\bibliography{ref}

\end{document}